\documentclass[12pt]{amsart}
\usepackage{amscd}
\usepackage{amssymb}
\usepackage{epic,eepic}

\allowdisplaybreaks[1]

\newtheorem{thm}{Theorem}[section]
\newtheorem{lem}[thm]{Lemma}
\newtheorem{cor}[thm]{Corollary}
\newtheorem{prop}[thm]{Proposition}

\theoremstyle{definition}

\newtheorem{rem}[thm]{Remark}
\newtheorem{defn}[thm]{Definition}

\newtheorem{ex}[thm]{Example}

\theoremstyle{remark}

\numberwithin{equation}{section}

\newcommand{\C}{{\mathbb{C}}}
\newcommand{\R}{{\mathbb{R}}}
\newcommand{\Z}{{\mathbb{Z}}}
\newcommand{\N}{{\mathbb{N}}}
\newcommand{\T}{{\mathbb{T}}}
\newcommand{\bw}{{\boldsymbol{w}}}

\DeclareMathOperator{\Ad}{Ad}
\DeclareMathOperator{\Tr}{Tr}

\DeclareMathOperator{\id}{id}

\DeclareMathOperator{\spe}{Sp}

\DeclareMathOperator{\spa}{span}
\DeclareMathOperator{\Hom}{Hom}
\DeclareMathOperator{\End}{End}
\DeclareMathOperator{\ev}{ev}
\DeclareMathOperator{\pr}{pr}
\DeclareMathOperator{\diag}{diag}

\def\ups{\upsilon}

\def\ovl{\overline}
\def\wdh{\widehat}
\def\wdt{\widetilde}
\def\op{\mbox{\tiny{op}}}

\def\suq{SU_q(2)}
\def\sumq{SU_{-q}(2)}
\def\soq{SO_q(3)}
\def\somq{SO_{-q}(3)}
\def\som{SO_{-1}(3)}
\def\summ{SU_{-1}(2)}
\def\csuq{C(SU_q(2))}
\def\asuq{A(SU_q(2))}
\def\csoq{C(SO_q(3))}
\def\csomq{C(SO_{-q}(3))}
\def\csum{C(SU_{-1}(2))}
\def\csom{C(SO_{-1}(3))}
\def\cso{C(SO_{1}(3))}

\def\uqsu{U_q(su_2)}
\def\umsu{U_{-1}(su_2)}

\def\CSq{C(S_{q,\lambda}^2)}
\begin{document}

\title[Compact quantum ergodic systems]
{Compact quantum ergodic systems}

\author[Reiji Tomatsu]{Reiji Tomatsu}
\address{Department of Mathematical Sciences,
University of Tokyo, 3-8-1 Komaba, Meguro, Tokyo 153-8914, Japan}
\email{tomatsu\char`\@ms.u-tokyo.ac.jp}

\begin{abstract}
We develop theory of multiplicity maps 
for compact quantum groups, as an application, 
we obtain a complete classification of 
right coideal $C^*$-algebras of $\csuq$ 
for $q\in [-1,1]\setminus \{0\}$. 
They are labeled with Dynkin diagrams, 
but classification results for positive and negative cases 
of $q$ are different. 
Many of the coideals are quantum spheres or quotient spaces by 
quantum subgroups, but we do have other ones in our classification list. 
\end{abstract}

\maketitle

\section{Introduction}

The aim of this paper is studying compact quantum ergodic 
systems and classifying right coideals of $\csuq$ for 
$q\in [-1,1]\setminus \{0\}$. 
Compact quantum groups have been introduced by S. L. Woronowicz 
\cite{Woronowicz2}. 
They are considered as the continuous function algebra on 
a non commutative space of a compact quantum group. 
The continuous function algebra on a compact quantum group 
is denoted by $C(G)$ and we have the coproduct 
$\delta:C(G)\longrightarrow C(G)\otimes C(G)$. 
In this paper, a quantum group discovered by V. G. Drinfel'd and 
M. Jimbo is called a quantum universal enveloping algebra 
because it is considered as a deformation of the universal 
enveloping algebra of a Lie algebra 
\cite{Drinfeld}, \cite{Jimbo}. 
A $C^*$-subalgebra $A\subset C(G)$ is called a right coideal $C^*$-algebra 
or often simply a right coideal 
if it is invariant by the right action of $G$, that is, we have 
$\delta(A)\subset A\otimes C(G)$. 
Let $A$ be a unital $C^*$-algebra and $G$ be a compact quantum 
group. 
A triple 
$\{A,G,\alpha\}$ is a compact quantum covariant system if a 
unital faithful $*$-homomorphism 
$\alpha:A\longrightarrow A\otimes C(G)$ satisfies 
$(\id\otimes\delta)\circ\alpha=(\alpha\otimes\id)\circ\alpha$. 
Moreover if it satisfies $A^\alpha:=\{a\in A\mid \alpha(a)=a\otimes1\}=\C1_A$, 
then it is called a compact quantum ergodic system. 
When we consider the classical compact group case, 
the invariant state on a $C^*$-algebra $A$ must be a faithful tracial state 
and 
multiplicities of irreducible representations are bounded by 
dimensions of their representation spaces \cite{KLS}. 
However in the quantum case the invariant state has no longer tracial property 
and the multiplicities are bounded by quantum dimensions which 
are larger than or equal to usual dimensions. 
Keeping in mind these phenomena, one shall notice there are 
rather differences between classical compact 
group case and compact quantum group case. 
It is still a fascinating problem to construct an example 
of a compact quantum ergodic system which has a multiplicity of 
an irreducible representation strictly larger than 
its dimension. 
In \cite{Wassermann1}, he has developed the theory of multiplicity maps 
and it has worked very well in the classification program 
of ergodic systems of $SU(2)$ \cite{Wassermann3}. 
Since multiplicity maps are defined via equivariant $K$-theory, 
one easily has its quantum version by using the work due to, 
for example, \cite{BaajSkandalis1,Boca,Vergnioux}. 
The most nontrivial important problem is whether there exists 
a common eigenvector of multiplicity maps or not. 
In compact group case its existence is guaranteed by the tracial 
invariant state. 
For that problem we show that such a vector exists if 
$A$ has a (not necessarily faithful) 
tracial state and $G=\suq$ 
in Theorem \ref{multiplicity vector}. 
If we consider a right coideal $A\subset C(G)$, 
then we get a compact quantum ergodic system $\{A,\delta,C(G)\}$. 
Fortunately 
we can show there always exists a desired common eigenvector 
in Theorem \ref{multiplicity vector for right coideal} in this case. 
With this result we can proceed 
to the classification of right coideals of 
$\csuq$. 
In the classical case we know the continuous function algebra on 
the homogeneous space $H\setminus G$ by a closed subgroup $H\subset G$ gives 
a right coideal and its converse holds via 
the Gelfand-Naimark theorem on abelian $C^*$-algebras. 
However in the quantum case, such correspondence does not hold in general. 
For example, in \cite{Podles1} he has constructed one-parameter family of, 
so called, quantum spheres $C(S_{q,\lambda}^2)$ which are 
right coideals and most of them are not obtained by taking 
quotient by subgroups. 
Therefore it is interesting to investigate the other non-quotient type 
right coideals. 
The aim of the latter half of this paper is 
giving the complete classification of right coideals of $\csuq$. 
Most important information of a right coideal is 
its spectral pattern, that is, multiplicities of irreducible representations. 
We will see 
this information is remembered by equivariant $K$-theory and 
we classify right coideals into 
the several types by connected graphs of norm $2$ 
as is used in the classification of ergodic systems of $SU(2)$. 
They are labeled by closed subgroups of $SU(2)$ or $\summ$ as 
McKay diagrams with respect to the fundamental two-dimensional 
representation. 
When we work on the negative $q$ case, 
it is needed to treat the graphs with a single loop at a 
vertex. 
Then we carry out case-by-case study of them. 
In that procedure one has to be careful of the essential effect of 
$|q|\neq1$. 
In fact, right coideals of some types such as 
the regular polyhedrons do not appear in those case. 
We state the main result on this classification. 

$(1)$ The case $0<q<1$: A right coideal must be one of type 
$1$, $SU(2)$, $\T_n$, $\T$ and $D_\infty^*$. 
When it is of type $\T$, then it is one of series of the quantum spheres. 
Otherwise it is uniquely determined by the type.

$(2)$ The case $-1<q<0$: A right coideal must be one of type 
$1$, $SU(2)$, $\T_n$, $\T$, $D_\infty^*$ and $D_1$. 
When it is of type $\T$, then it is one of series of the quantum spheres. 
Otherwise it is uniquely determined by the type.

$(3)$ The case $q=-1$: 
If a right coideal $A$ is not of type $\T_n$ (odd $n\geq3$) 
or $D_n$ 
(odd $n\geq1$), 
there exists a closed subgroup $H$ in $SO_{-1}(3)$ such that 
$A$ is $C(H\setminus \som)$. 
If a right coideal $A$ is of type $\T_n$ (odd $n\geq3$), 
$A$ is conjugated to 
$C(\T_n\setminus \summ)$ or 
$C^*(\eta^{\frac{n}{2}},\wdh{\eta}^{\frac{n}{2}})$. 
If a right coideal $A$ is of type $D_1$, then 
$A$ is conjugated to $C(D_1\setminus \summ)$. 
If a right coideal $A$ is of type $D_n$ (odd $n\geq3$), 
$A$ is conjugated to 
$C(D_n\setminus \summ)$ or $C^*(\eta^{\frac{n}{2}})$. 
Here conjugation is given by the left action 
$\beta_z^L$ of the maximal torus for some $z\in\T$. 

In each of the above cases, uniqueness means conjugation by 
the left action of the maximal torus $\T$. 
We find right coideals which are 
of type $D_\infty^*$ in the case $0<q<1$, 
of type $D_\infty^*$ and $D_1$ in the case $-1<q<0$ 
and 
$C^*(\eta^{\frac{n}{2}},\wdh{\eta}^{\frac{n}{2}})$, 
$C^*(\eta^{\frac{n}{2}})$ in the case $q=-1$. 
They are not the quantum spheres nor the quotient spaces. 
In the case $q=-1$, a right coideal which is not of type $\T_n$ (odd $n\geq3$) 
is $\summ$-isomorphic to each other in the same type, however, 
$C(\T_n\setminus \summ)$ is not $\summ$-isomorphic to 
$C^*(\eta^{\frac{n}{2}},\wdh{\eta}^{\frac{n}{2}})$. 

We briefly explain the contents of study of each section. 
In section $2$ we collect required facts on compact quantum groups 
and their actions. 
In section $3$ we study general compact quantum ergodic systems, 
especially hom-spaces of them. 
In section $4$ we consider quantum version of the theory of 
multiplicity maps. 
In section $5$ we give a summary for representation theory of $\suq$. 
In sections $6,7$ and $8$ 
we carry out classification of right coideals of $\csuq$ for 
$0<q<1$, $-1<q<0$ and $q=-1$, respectively. 
In Appendix we list all the connected graphs of norm $2$. 
If one has only interest in classification of 
right coideals, 
he or she can skip results in sections $1$, $2$, $3$ and $4$ except for 
Corollary \ref{multiplicity vector for right coideal}. 

\textbf{Acknowledgements}. 
The author is grateful to Yasuyuki Kawahigashi, 
Shigeru Yamagami and 
Junichi Shiraishi for permanent encouragement and 
fruitful discussions. 
This work was partially supported by Research Fellowship for 
Young Scientists of the Japan Society for the Promotion of 
Science.

\section{Preliminaries}

We collect basic notions on compact quantum groups and 
their actions. 
Our standard references are 
\cite{Woronowicz2} for theory of compact quantum groups 
and \cite{BaajSkandalis2} for their actions. 
In this paper we only treat minimal tensor products for $C^*$-algebras 
and use the notation simply $\otimes$. 
Although it is not essential, 
separability of compact quantum groups is always assumed. 

\begin{defn}\cite[p.853]{Woronowicz2}
Let $A$ be a unital $C^*$-algebra and 
$\delta:A \longrightarrow A\otimes A$ be a 
$*$-homomorphism. If they satisfy the following conditions, the pair 
$(A,\delta)$ is called a compact quantum group.

\begin{enumerate}

\item The map $\delta$ satisfies coassociativity condition:
\[(\delta\otimes\id_{A})\circ \delta
=(\id_{A}\otimes\delta)\circ\delta.\]

\item The vector spaces $\delta(A)(\C\otimes A)$ and $\delta(A)(A\otimes\C)$ 
are dense in $A\otimes A$.

\end{enumerate}

\end{defn}

Let $(A,\delta)$ be a compact quantum group. Then there exists the unique 
state $h$ which satisfies the following invariance condition:
\[
(h\otimes\id_{A})\circ\delta(a)=(\id_{A}\otimes h)\circ\delta(a)=h(a)1_{A},
\mbox{for all} \  a\in A.
\]

This state is called the \textit{Haar state}. We always consider a compact 
quantum group whose Haar state is faithful. A compact quantum group 
$(A,\delta)$ is often regarded as the continuous function algebra 
on a non-commutative space $G$ and we write $G=(A,\delta)$ and $A=C(G)$.

Let  $H$ be a Hilbert space and $v$ be a unitary 
in $M(\mathbb{K}(H)\otimes A)$. 
A unitary $v$ is called a 
\textit{unitary representation} of $G$ 
if it satisfies the following equality:
\[
(\id_{H}\otimes\delta)(v)=v_{12}v_{13}.
\]

Let $\{v_{i}\}_{i\in I}$ 
be a family of unitary representations on Hilbert spaces 
$\{H_{v_{i}}\}_{i\in I}$, respectively. 
The direct sum representation $\prod_{i \in I}v_{i}$ is defined as 
the direct sum of unitary operators via the natural inclusion 
$\prod_{i\in I}M(\mathbb{K}(H_{v_{i}})\otimes C(G))
\subset M(\mathbb{K}(\bigoplus_{i\in I}H_{v})
\otimes C(G))$. 
Let $v,w$ be unitary representations on Hilbert spaces $H_{v},H_{w}$ 
respectively.
The tensor product representation $v\dot{\otimes}w$ is defined as a 
unitary operator $v_{13}w_{23}$ which is an element of 
$M(\mathbb{K}(H_{v}\otimes H_{w})\otimes C(G))$. An operator $T$ 
in $\mathbb{B}(H_{v},H_{w})$ is called an \textit{intertwiner} of $v$ 
and $w$ if it satisfies $(T\otimes 1)v=w(T\otimes 1)$. 
The set of intertwiners between $v$ and $w$ is written as $\Hom(v,w)$ and 
it is called a \textit{hom-space}.
For one unitary representation $v$ its hom-space $\Hom(v,v)$ becomes a 
$C^*$-subalgebra of $\mathbb{B}(H_{v})$. If the space is trivial 
(that is, $\Hom(v,v)=\C1_{H_{v}}$), 
we say that the unitary representation $v$ is 
\textit{irreducible}. 
Let $v$ be a unitary representation of a compact quantum group 
$G$ and $v$ is a direct sum of finite dimensional irreducible 
representations (\cite[p.864]{Woronowicz2}). 
In particular any irreducible unitary representation is finite dimensional.

In the set of all irreducible unitary representations of $G$, 
we define unitary equivalence relation naturally 
and denote its quotient set by $\widehat{G}$. 
We can select one unitary representation $w(\pi)\in 
\mathbb{B}(H_{\pi})\otimes C(G)$ 
for each equivalence class $\pi\in \widehat{G}$. 
Let $\{\xi_{i}\}_{i\in I_{\pi}}$ be an orthonormal basis for $H_{\pi}$. 
We can identify $w(\pi)$ 
with the matrix $\{w(\pi)_{i,j}\}_{i,j\in I_{\pi}}$ 
where $w(\pi)_{i,j}$ are elements of $C(G)$. 
From the definition of unitary representation we obtain the following 
equality:
\[
\delta(w(\pi)_{i,j})=\sum_{k\in I_{\pi}}w(\pi)_{i,k}\otimes w(\pi)_{k,j}
,\ \mbox{for all}\  i,j\in I_{\pi}.
\]
Let us define the \textit{smooth function algebra} 
$A(G)=\spa\{w(\pi)_{i,j}\mid i,j\in I_{\pi},\  \pi\in \widehat{G}\}$. 
It is in fact a unital $*$-subalgebra dense in $C(G)$ and 
the set $\{w(\pi)_{i,j}\mid i,j\in I_{\pi}, \pi\in \widehat{G}\}$ is a linear 
basis for $A(G)$. 
Introducing the maps $\varepsilon:A(G)\longrightarrow \C$ and 
$\kappa:A(G)\longrightarrow A(G)$ by
\[
\varepsilon(w(\pi)_{i,j})=\delta_{i,j},
\  
\kappa(w(\pi)_{i,j})=w(\pi)_{j,i}^*
\]
 for each $\pi \in \widehat{G}$ and $i,j\in I_{\pi}$ gives the algebra 
 $A(G)$ a \textit{Hopf}$*$-\textit{algebra} structure:

\begin{enumerate}

\item The map $\delta:A(G)\longrightarrow A(G)\otimes A(G)$ satisfies 
coassociativity:
\[
(\delta\otimes\id)\circ\delta=(\id\otimes\delta)\circ\delta.
\]

\item The map (called the \textit{counit}) 
$\varepsilon: A(G) \longrightarrow \C$ is a 
unital $*$-homomorphism and satisfies
\[
(\varepsilon\otimes\id)\circ\delta=(\id\otimes\varepsilon)\circ\delta=\id.
\]

\item The map (called the \textit{antipode})
$\kappa:A(G)\longrightarrow A(G)$ is a linear anti-multiplicative map and 
satisfies
\[
m(\kappa\otimes\id)\circ\delta=m(\id\otimes\kappa)\circ\delta=\varepsilon
\]
and
\[
\kappa(\kappa(a^*)^*)=a \  \mbox{for all}\  a\in A(G)
.\]
\end{enumerate} 

Let $A(G)^*$ be an algebraic dual space of $A(G)$. 
For $a\in A(G)$ and $\theta,\omega_{1},\omega_{2}\in A(G)^*$, 
we define their convolution products and involution:
\[
\theta*a=(\id\otimes\theta)(\delta(a)),
\]
\[
a*\theta=(\theta\otimes\id)(\delta(a)),
\]
\[
\omega_{1}*\omega_{2}=(\omega_{1}\otimes\omega_{2})\circ\delta
,\]
\[
\theta^*(a)=\overline{\theta(\kappa(a)^*)}.
\]
With these operations the dual space 
$A(G)^*$ becomes a unital $*$-algebra 
(its unit is the counit) and acts on $A(G)$ from the left and right.
On the smooth function algebra $A(G)$ there exists a unique 
family of characters 
$\{f_z\}_{z\in\C}$which are called \textit{Woronowicz characters} 
with following properties:

\begin{enumerate}

\item $f_z(1)=1$ for all $z\in \C$.

\item For any element $a\in A(G)$, the mapping $z\in\C \mapsto f_z(a)\in\C$ 
is an entire holomorphic function.

\item $f_{z_{1}}*f_{z_{2}}=f_{z_{1}+z_{2}}$ for all $z_{1},z_{2}\in\C$.

\item $f_0=\varepsilon$.

\item $f_z(\kappa(a))=f_{-z}(a)$ for all $z\in \C, a\in A(G)$.

\item $f_{z}(a^*)=\overline{f_{-\bar{z}}(a)}$ for all $z\in \C, a\in A(G)$.

\end{enumerate}

We can define one-parameter automorphism groups 
$\{\sigma^{h}_{t}\}_{t\in\R}$ and 
$\{\tau_{t}\}_{t\in \R}$ and $*$-anti-multiplicative linear map $R$ by
\[
\sigma^{h}_{t}(a)=f_{it}*a*f_{it} \  \mbox{for all}\  a\in A(G),
\]
\[
\tau_{t}(a)=f_{-it}*a*f_{it} \  \mbox{for all}\  a\in A(G),
\]
\[
R(a)=f_{\frac{1}{2}}*\kappa(x)*f_{-\frac{1}{2}} 
\  \mbox{for all}\  a\in A(G).
\]
They are called the \textit{modular automorphism group}, 
the \textit{scaling automorphism group} and the \textit{unitary antipode} 
respectively. 
They are norm continuously extendable to the continuous function algebra 
$C(G)$. The Haar state $h$ is a $\sigma^{h}$-KMS state: 
for any $a,b\in A(G)$, we have
\[
h(ab)=h(\sigma^{h}_{i}(b)a).
\]
Relations of these maps are as follows.
\[
\sigma^{h}_{s}\circ\tau_{t}=\tau_{t}\circ\sigma^{h}_{s}
\  \mbox{for all}\  s,t\in\R,
\]
\[
R\circ\sigma^{h}_{s}=\sigma^{h}_{s}\circ R,\  
R\circ\tau_{t}=\tau_{t}\circ R,
\]
\[
\kappa=R\circ\tau_{-\frac{i}{2}}=\tau_{-\frac{i}{2}}\circ R.
\]

Let $v$ be a unitary representation on a finite dimensional 
Hilbert space $H_v$ and $j_v :H_v\longrightarrow \overline{H_v}$ be 
a conjugate unitary map where $\overline{H_v}$ is the conjugate Hilbert 
space. The transpose map 
$t_v :\mathbb{B}(H_v)\longrightarrow \mathbb{B}(\overline{H_v})$ is defined 
by $t_v(x)=j_v x^* j_v^{-1}$. 
We define the \textit{contragradient representation} 
(which is not necessarily a unitary representation) 
$v^c$ and the conjugate unitary representation $\bar{v}$
of a finite dimensional unitary representation $v$ by
\[
v^c=(t_v\otimes \kappa)(v),\  \bar{v}=(t_v\otimes R)(v).
\]

The \textit{F-matrix} of a finite dimensional unitary representation 
$v\in\mathbb{B}(H_v)\otimes C(G)$ is defined by
\[
F_v=(\id\otimes f_1)(v)\ \in \mathbb{B}(H_v).
\]
This matrix is strictly positive definite. It is known that 
$F_{v}^z=(\id\otimes f_{z})(v)$ for any $z\in\C$ 
and $\Tr(F_v)=\Tr(F_v^{-1})$. 
For a unitary irreducible representation $v$, its second contragradient 
(not necessarily unitary) representation $v^{cc}$ is also irreducible and 
equivalent to $v$. In fact we have $\Hom(v,v^{cc})=\C F_v$. 
The value $\Tr(F_v)$ is called the 
\textit{quantum dimension} and denoted by $D_{v}$. 
Since we have $F_{\bar{v}}=t_{v}(F_{v}^{-1})$, 
the quantum dimension of $v$ is equal to that of $\bar{v}$. 
We write $d_v$ for the usual dimension of the vector space $H_v$. 
The quantum dimension is larger than or equal to the usual dimension. 
When we fix a selection of $\{w(\pi)\}_{\pi\in\widehat{G}}$, we write 
$H_\pi$, $F_\pi$, $D_{\pi}$ and $d_{\pi}$ 
for $H_{w(\pi)}$, $F_{w(\pi)}$, $D_{w(\pi)}$ and $d_{w(\pi)}$, 
respectively. 

Using $F$-matrices, we have the following equalities 
about the Haar state:
\[
h(w(\pi)_{i,j}w(\rho)_{r,s}^*)
=D_{\pi}^{-1}(F_{\pi})_{s,j}\delta_{\pi,\rho}\delta_{i,r}
\]
\[
h(w(\pi)_{i,j}^{*}w(\rho)_{r,s})
=D_{\pi}^{-1}(F_{\pi}^{-1})_{r,i}\delta_{\pi,\rho}\delta_{j,s}
\]
for any $\pi,\rho\in\widehat{G}$, $i,j\in I_\pi$ and $r,s \in I_\rho$. 

For the tensor product $H\otimes K$ of two Hilbert spaces $H$ and $K$, 
we define its conjugate unitary map by 
$j_{H\otimes K}:H\otimes K \longrightarrow 
\overline{K}\otimes\overline{H},
j_{H\otimes K}(\xi\otimes\eta)=j_{K}\eta\otimes j_{H}\xi$. 
Then for two finite dimensional unitary representation $v$ and $w$, 
we obtain $\overline{v\dot{\otimes}w}=\bar{w}\dot{\otimes}\bar{v}$. 

Let $L^2(G)$ be the GNS-representation space of $C(G)$ associated to 
the Haar state $h$. 
Its cyclic vector is denoted by $\hat{1}_h$. 
We define the \textit{modular conjugation} $J$ and $\hat{J}$ by
\begin{align*}
J \,x \hat{1}_h =&\,\sigma_{\frac{i}{2}}^h(x)^* \hat{1}_h \, ,\\
\hat{J}\, x \hat{1}_h=&\, R(x)^* \hat{1}_h\, 
\end{align*}
where $x$ is an element of $A(G)$. 
We also define a unitary $U=J\hat{J}=\hat{J}J$. 
When we consider $L^\infty(G)$ which is the $\sigma$-weak closure 
of $C(G)$ in $\mathbb{B}(L^2(G))$, it becomes 
a \textit{von Neumann algebraic compact quantum group}. 
Precisely speaking, the coproduct $\delta$ and the Haar state $h$ 
extend to $L^\infty(G)$ as a normal $*$-homomorphism 
and the invariant faithful normal state respectively. 
By the invariance of $h$, we define the following two unitaries 
on $L^2(G)\otimes L^2(G)$: for all $x,y\in C(G)$, 
\[
V_{\ell}^*(x\hat{1}_h\otimes y\hat{1}_h)=\delta(y)(x\hat{1}_h\otimes\hat{1}_h)
\, ,
\]
\[
V(x\hat{1}_h\otimes y\hat{1}_h)=\delta(x)(\hat{1}_h\otimes y\hat{1}_h)
\, .
\]
They satisfy the following \textit{pentagonal equality} and therefore 
are called \textit{multiplicative unitaries}:
\[
V_{12}V_{13}V_{23}=V_{23}V_{12}
\, .
\]
We also use other unitaries:
\begin{align*}
\widetilde{V}=&\,\Sigma (1\otimes U) V (1\otimes U) \Sigma\, , \\
W=&\, (U\otimes1)V(U\otimes1)\, ,
\end{align*}
where $\Sigma:L^2(G)\otimes L^2(G)\longrightarrow L^2(G)\otimes L^2(G)$ 
is a flipping unitary $\xi\otimes\eta\mapsto \eta\otimes\xi$. 
The unitary $\widetilde{V}$ satisfies the above pentagonal equality and 
$W$ satisfies $W_{23}W_{13}W_{12}=W_{12}W_{23}$. 
It is easy to see 
$V_{\ell}
=\Sigma (\hat{J}\otimes\hat{J})V^*(\hat{J}\otimes\hat{J})\Sigma
=(1\otimes U)\Sigma V\Sigma (1\otimes U)
$. 
We define a left $G$-action on $H_\pi$ by 
$\xi_i^\pi\mapsto \sum_{j\in I_\pi}w(\pi)_{i,j}\otimes \xi_j^\pi$ and 
denote it by $H_\pi^{\ell}$. 
Let $\Theta_\pi$ be a unitary map 
$H_\pi^{\ell}\otimes H_\pi\longrightarrow L^2(G)$ 
defined by 
$\Theta_\pi(\xi_i\otimes\xi_j)=\sqrt{D_\pi (F_\pi)_{i,i}}w(\pi)_{i,j}$. 
It intertwines the left and right $G$-actions. 
Then we have the \textit{Peter-Weyl} decomposition 
$\Theta:\oplus_{\pi\in\wdh{G}}
H_\pi^{\ell}\otimes H_\pi \longrightarrow L^2(G)$ 
with $\Theta=\oplus_{\pi\in\wdh{G}}\Theta_\pi$. 
The \textit{left, right reduced group $C^*$-algebras} are 
defined by
\[
C_{\ell}^*(G)=\ovl{\spa}
\{(\omega\otimes\id)(V_{\ell})\mid \omega\in \mathbb{B}(L^2(G))_{*}\}
\, ,
\]
\[
C_r^*(G)=\ovl{\spa}
\{(\id\otimes\omega)(V)\mid \omega\in \mathbb{B}(L^2(G))_{*}\}
\, ,
\]
where the closure is taken with respect to the operator norm of 
$\mathbb{B}(L^2(G))$. 
Note that $C_r^*(G)$ is contained in the commutant algebra of 
$C_{\ell}^*(G)=\hat{J}C_r^*(G)\hat{J}$. 
There is a distinguished projection $p_0=(\id\otimes h)(V)$ in 
$C_r^*(G)$. 
It is also included in $C_{\ell}^*(G)$ and hence 
it is a central projection of $C_r^*(G)$. 
Moreover it is a minimal projection of $\mathbb{K}(L^2(G))$, 
in fact, we have $p_0 L^2(G)=\C \hat{1}_h$. 
We now have a map $\rho:\mathbb{B}(L^2(G))_*\longrightarrow C_r^*(G)$ 
by $\rho(\omega)=(\id\otimes \omega)(V)$. 
Let us define a normal functional 
$\theta(\pi)_{i,j}(x)=D_\pi (F_\pi)_{j,j}^{-1}h(xw(\pi)_{i,j}^*)$ 
for $x\in \mathbb{B}(L^2(G))$. 
Then the linear space 
$\spa\{\rho(\theta(\pi)_{i,j})\mid i,j\in I_\pi,\pi\in \widehat{G}\}$ is 
dense in $C_r^*(G)$ and moreover we have 
$\rho(\theta(\pi)_{i,j})\rho(\theta(\sigma)_{p,q})
=
\delta_{\pi,\sigma}\delta_{j,p}
\rho(\theta(\pi)_{i,q})\, .
$
Hence we obtain a $*$-isomorphism: 
\[
C_r^*(G)\longrightarrow \bigoplus_{\pi\in \widehat{G}}\End(H_\pi)
\]
defined by $\rho(\theta(\pi)_{i,j})\mapsto E_{i,j}^{\pi}$ which 
is a matrix unit. 
With this identification, the action of $E_{i,j}^{\pi}$ on $L^2(G)$ is 
given by
\[
E_{i,j}^{\pi} w(\sigma)_{p,q}\hat{1}_h
=
\delta_{\pi,\sigma}\delta_{j,q}w(\pi)_{p,i}\hat{1}_h
\, .
\]
Using this fact and the Peter-Weyl decomposition, we see that 
the $C^*$-algebras $C_{\ell}^*(G)$ and $C_r^*(G)$ act on the first and 
the second tensor components of the Hilbert space 
$\oplus_{\pi\in\wdh{G}}H_\pi^{\ell}\otimes H_\pi$ respectively. 
This is just the differential representation of 
$\mathbb{B}(L^2(G))_*\subset A(G)^*$. 
Notice that $E_{i,j}^\pi$ is an operator of rank $n$ and 
hence $C_{\ell}^*(G)$ and $C_r^*(G)$ are contained in $\mathbb{K}(L^2(G))$. 
The projection $*$-homomorphism 
$C_r^*(G)\longrightarrow \End(H_\pi)$ is denoted by 
$\pr_\pi$. 

The \textit{dual coproduct} is defined by
$\hat{\delta}:C_r^*(G)\longrightarrow M(C_r^*(G)\otimes C_r^*(G))$ by 
\[
\hat{\delta}(x)=V^*(1\otimes x)V \quad \mbox{for\ all\ }x\in C_r^*(G)
\, 
\]
and it follows that linear subspaces 
$\hat{\delta}(C_r^*(G))(\C\otimes C_r^*(G))$ and 
$\hat{\delta}(C_r^*(G))(C_r^*(G)\otimes \C)$ are dense in 
$C_r^*(G)\otimes C_r^*(G)$. 
On $C_{\ell}^*(G)$, we define 
$\wdh{\delta_{\ell}}(x)=\Sigma V_{\ell}^*(x\otimes1)V_{\ell} \Sigma$. 

Let $R(G)$ be a $\Z$-module 
$\bigoplus_{\pi\in\widehat{G}}\Z \pi$. 
For two $\pi,\sigma\in \widehat{G}$, we define their product via 
irreducible decomposition of the tensor product representation 
$w(\pi)\dot{\otimes}w(\sigma)$. 
More precisely, let $N^{\pi\,\sigma}_{\tau}$ be the dimension of 
the \textit{hom-space} 
$\Hom(w(\tau),w(\pi)\dot{\otimes}w(\sigma))$, which does not 
depend on the choice of representatives of 
$\{w(\pi)\}_{\pi\in\widehat{G}}$. Then we have
\[
\pi\cdot\sigma=\sum_{\tau\in\widehat{G}} N^{\pi\,\sigma}_{\tau}\tau.
\]
The module $R(G)$ has an involution, that is, 
$\ovl{\pi}$ is an equivalence class of $\overline{w(\pi)}$, 
which also does not depend on the choice of representatives of 
$\{w(\pi)\}_{\pi\in\widehat{G}}$. 
By these operations, $R(G)$ has an involutive $\Z$-ring structure 
and it is called the \textit{representation ring} of $G$. 
We define the positive cone $R(G)_+=\bigoplus_{\pi\in\wdh{G}}\Z_{\geq0}\pi$. 

Since a hom-space $\Hom(w(\tau),w(\pi)\dot{\otimes}w(\sigma))$ is 
naturally 
isomorphic to a hom-space 
$\Hom(w(\pi),w(\tau)\dot{\otimes}\overline{w(\sigma)})$ or 
$\Hom(w(\sigma),\overline{w(\pi)}\dot{\otimes}w(\tau))$, 
we have symmetry of 
$N^{\pi\,\sigma}_{\tau}=N^{\tau\,\bar{\sigma}}_{\pi}
=N^{\bar{\pi}\,\tau}_{\sigma}$. 

From now on, 
we use the letter $A$ not for a continuous function algebra 
on a quantum group but for an arbitrary $C^*$-algebra. 
We give a brief summary of compact quantum group actions. 
The following definition is standard. 

\begin{defn}
Let $A$ be a (not necessarily unital) $C^*$-algebra 
and $(C(G),\delta)$ be a compact quantum group and 
$\alpha:A\longrightarrow A\otimes C(G)$ be a $*$-homomorphism. 
The triple $\{A,G,\alpha \}$ is called compact quantum covariant system 
(or simply covariant system) if the following statements hold:

\begin{enumerate}

\item The map $\alpha$ is injective and satisfies 
$(\alpha\otimes\id)\circ\alpha=(\id\otimes\delta)\circ\alpha$. 

\item The vector space $\alpha(A)(\C\otimes C(G))$ 
is dense in $A\otimes C(G)$.

\end{enumerate}

\end{defn}
If we consider a compact quantum covariant system $\{A,G,\alpha\}$ 
where $A$ is unital, 
we always assume that $\alpha$ is a unital $*$-homomorphism. 
Let $\{A,G,\alpha\}$ be a compact quantum covariant system. 
We can select representatives $\{w(\pi)\}_{\pi\in\widehat{G}}$ 
of $\widehat{G}$ whose $F$-matrices are diagonal 
and fix them in this section. 
We give comodule structure 
$\Gamma_{\pi}:H_{\pi}\longrightarrow H_{\pi}\otimes A(G)$ 
to a representation space $H_{\pi}$ by 
fixing an orthonormal basis $\{\xi_{j}^{\pi}\}_{j\in I_{\pi}}$: 
$\Gamma_{\pi}(\xi^{\pi}_{j})=\sum_{k\in I_{\pi}}\xi^{\pi}_{k}\otimes w(\pi)_{k,j}$ 
for any $j\in I_{\pi}$. 
For $\pi\in\widehat{G}$, define the functional 
$\theta_{\pi}$ on $C(G)$ by 
$\theta_{\pi}(x)
=D_{\pi}h(x(\sum_{i\in I_\pi}(F_{\pi}^{-1})_{i,i}w(\pi)_{i,i}))$ for any 
$x\in C(G)$. 
As in the case of a compact group \cite[Theorem 2]{Shiga}, 
$A$ is completely decomposable. 
Define the 
linear map $P_{\pi}=(\id\otimes\theta_{\pi})\circ\alpha$ on $A$. 
From simple calculations, we have 
$P_{\pi}P_{\rho}=\delta_{\pi,\rho}P_{\pi}$. 
Put $A_\pi$ for the range space of the projection $P_{\pi}$. 
Then the linear subspace $\mathcal{A}=\bigoplus_{\pi\in\widehat{G}}A_{\pi}$ 
is norm dense in $A$. 
Let $\Hom_{G}(H_{\pi}, A)$ be a set of intertwiners, 
that is, $\Hom_{G}(H_{\pi}, A)\ni S$ 
if and only if $\alpha\circ S=(S\otimes\id)\circ\Gamma_{\pi}$. 
Then we have $A_{\pi}=\overline{\spa}
\{ S\xi \mid \xi\in H_{\pi}, S\in \Hom_{G}(H_{\pi}, A)\}$. 
The dimension of $\Hom_{G}(H_{\pi}, A)$ is called 
the \textit{multiplicity} of $\pi$. 
Especially for the trivial representation $\pi_{0}$, 
the closed linear subspace $A_{\pi_0}$ coincides with the 
fixed point algebra $A^\alpha=\{x\in A\mid \alpha(x)=x\otimes1\}$. 
Let $m(\pi)$ be the multiplicity of $\pi$. 
We often use a formal symbol 
$\oplus_{\pi\in\widehat{G}}m(\pi)\pi$ which is 
called the \textit{spectral pattern} of $\{A,G,\alpha\}$ (or simply $A$). 

Next we recall the notion of the \textit{crossed product}. 
Let $\{A,G,\alpha\}$ be a covariant system. 
Consider the Hilbert $A$-module $A\otimes L^2(G)$. 
In a $C^*$-algebra 
$\mathbb{K}(A\otimes L^2(G))=A\otimes \mathbb{K}(L^2(G))$ 
define the crossed product $C^*$-algebra: 
\[A\rtimes_\alpha G=\overline{\alpha(A)(\C\otimes C_r^*(G))}\, ,\]
where the closure of linear space is taken with respect to the operator norm. 
It is also characterized as a fixed point algebra of 
$\{A\otimes\mathbb{K}(L^2(G)),G,\widetilde{\alpha}\}$: 
$A\rtimes_{\alpha}G=(A\otimes\mathbb{K}(L^2(G)))^{\widetilde{\alpha}}$, 
where $\widetilde{\alpha}(a\otimes k)
=W_{23}\alpha(a)(1\otimes k\otimes 1)W_{23}^*$ 
for all $a\in A$ and $k\in \mathbb{K}(L^2(G))$. 
Note that $\alpha(A)$, ${C_{r}}^*(G)$ are $C^*$-subalgebras 
of $M(A\rtimes_\alpha G)$ naturally. 
When we consider $A=C(G)$ and $\alpha=\delta$, 
we have a natural isomorphism:
\[
C(G)\rtimes_\delta G \longrightarrow \mathbb{K}(L^2(G))
\]
defined by $\delta(x)(1\otimes y)\mapsto xy$. 
The \textit{dual coaction} 
$\hat{\alpha}:A\rtimes_\alpha G \longrightarrow 
M(A\rtimes_\alpha G\otimes C_r^*(G))$ is defined by 
\[
\hat{\alpha}(x)=\widetilde{V}_{23}(x\otimes1)\widetilde{V}_{23}^* \, 
\]
which satisfies the density condition that a linear space 
$\hat{\alpha}(A\rtimes_\alpha G)(\C\otimes C_r^*(G))$ is 
dense in $A\rtimes_\alpha G\otimes C_r^*(G)$.  
We define a $*$-homomorphism 
$\hat{\alpha}_\pi:A\rtimes_\alpha G\longrightarrow 
A\rtimes_\alpha G\otimes \End(H_\pi)$ 
by the composition 
$(\id\otimes\pr_\pi)\circ \hat{\alpha}$. 
Now we introduce the famous duality theorem. 
Let $\{A,G,\alpha\}$ be a covariant system. 
We consider the $*$-homomorphism 
$\hat{\alpha}^o=
(\id\otimes \Ad U)\circ \hat{\alpha}:A\rtimes_\alpha G\longrightarrow 
M(A\rtimes_\alpha G\otimes C_{\ell}^*(G))$. 
It satisfies 
$(\id\otimes \wdh{\delta_{\ell}}^{\op})\circ \hat{\alpha}^o
=
(\hat{\alpha}^o\otimes\id)\circ\hat{\alpha}^o
$, that is, $\hat{\alpha}^o$ is a coaction of the opposite 
discrete quantum group $C_{\ell}^*(G)$. 
Define the crossed product 
$A\rtimes_\alpha G \rtimes_{\hat{\alpha}^o} (\wdh{G})'{}^{\op}$
 by the norm closure of the linear 
space $\hat{\alpha}^o(A\rtimes_\alpha G)(\C\otimes\C\otimes C(G))$. 
Its dual action $\wdh{\hat{\alpha}^o}$ is given by 
$\wdh{\hat{\alpha}^o}(x)=V_{34}(x\otimes1)V_{34}^*$. 
The duality theorem says there exists an isomorphism 
between $A\rtimes_\alpha G \rtimes_{\alpha^o} (\wdh{G})'{}^{\op}$ and 
$A\otimes \mathbb{K}(L^2(G))$ such that 
$\hat{\alpha}^o\big{(}\alpha(a)(1\otimes x)\big{)}(1\otimes1\otimes y)$ 
is mapped to $\alpha(a)(1\otimes x)(1\otimes UyU)$ 
for all $a\in A$, $x\in C_r^*(G)$ and $y\in C(G)$. 
This isomorphism conjugates $\wdh{\hat{\alpha}^o}$ and 
$\widetilde{\alpha}$. 

Let $G=(C(G),\delta_G)$ and $H=(C(H),\delta_H)$ be 
compact quantum groups. We say $H$ is a 
\textit{quantum subgroup} (or sometimes simply called a subgroup) 
if there exists a surjective $*$-homomorphism 
$r:C(G)\longrightarrow C(H)$ such that $r$ satisfies 
$\delta_H\circ r=(r\otimes r)\circ \delta_G$. 
This surjection is often called the restriction homomorphism. 
In general a choice of $r$ is not unique. 
The \textit{(left) quotient space} is defined by 
\[C(H\setminus G)=\{ x\in C(G)
\mid (r\otimes\id)\circ \delta_G(x)=1\otimes x\}\, .\] 
Actually it depends on $r$ and we shall denote it by 
${}^H C(G)$, however, we abuse the notation 
of $H\setminus G$. 
We have to be careful of $r$ when we consider 
a quotient space. 

\section{Ergodic systems}

If a compact quantum covariant system $\{A,G,\alpha\}$ 
with a unital $C^*$-algebra $A$ satisfies $A^\alpha=\C1$, 
then we call it a \textit{compact quantum ergodic system} 
(or simply an ergodic system). 
We investigate the ergodic system in this section. 
The faithful conditional expectation $\varphi=(\id\otimes h)\circ\alpha$ onto 
$A^\alpha=\C1$ is the unique invariant state in this case, where invariance 
means: $(\varphi\otimes\id)(\alpha(x))=\varphi(x)1$ for any $x\in A$. 
Note that $\varphi$ is a trace if $G$ is a compact group 
\cite[Theorem 4.1]{KLS}. 
In \cite[Theorem 17]{Boca}, 
it has been proved that the dimension of $A_{\pi}$ is 
less than or equal to $D_{\pi}^2$ for any $\pi\in\widehat{G}$. 
We will make a sharper estimate of these dimensions. 
In order to do this, we have to explain the existence of 
the modular automorphism group with respect to the state $\varphi$. 
By \cite[Proposition 18]{Boca}, there exists the unital 
multiplicative linear map $\Theta:\mathcal{A}\longrightarrow \mathcal{A}$ 
such that $\varphi(ab)=\varphi(\Theta(b)a)$ for any $a,b\in \mathcal{A}$. 
Let $M$ be the von Neumann algebra $\pi_{\varphi}(A)''$ associated to 
the state $\varphi$ via GNS-representation 
$\{H_{\varphi},\pi_{\varphi},\xi_{\varphi}\}$ 
and the extension of $\varphi$ to $M$ is also denoted by $\varphi$. 
We often identify $A$ with $\pi_{\varphi}(A)$. 
Let $p\in M$ be a projection with $\varphi(p)=0$. 
Then for any $a\in \mathcal{A}$ we have 
$\varphi(a^{*}pa)=\varphi(\Theta(a)a^{*}p)$. 
This is in fact equal to $0$ by the Cauchy-Schwarz inequality. 
From this we see $p\mathcal{A}\xi_{\varphi}=0$ and $p=0$, 
because the linear subspace $\mathcal{A}\xi_{\varphi}$ 
is dense in $H_{\varphi}$. 
Hence $\varphi$ is a faithful normal state on $M$ and there exists the 
modular automorphism group $\{\sigma^{\varphi}\}_{t\in\R}$ on $M$. 
Since $\varphi$ is an invariant state on $A$ for the action of $G$, 
the action $\alpha$ extends to the action on $M$. 
An action of a compact quantum group on a von Neumann algebra is 
defined similarly to the case of a $C^*$-algebra 
by using the $\sigma$-weak topology. 
Note the following useful equality about the modular automorphism group 
and the scaling automorphism group:
\[
\alpha\circ\sigma^{\varphi}_{t}
=(\sigma^{\varphi}_{t}\otimes\tau_{-t})\circ\alpha
\  \mbox{for all}\  t\in\R
.\]
For its proof, readers are referred to 
\cite[Th\'{e}or\`{e}me 2.9]{Enock1}. 
The spectral subspace $M_\pi$ 
is the $\sigma$-weak closure of $A_{\pi}$, however, they coincide 
because of the finite dimensionality of $A_{\pi}$. 

\begin{lem}\label{R-action}
Let $S$ be an intertwiner of $H_{\pi}$ and $A$. 
For any $t\in\R$, 
we define the map $S_{t}: H_{\pi}\longrightarrow A$ by 
$S_{t}\xi^{\pi}_{j}=(F_{\pi})_{j,j}^{-it}\sigma^{\varphi}_{t}(S\xi^{\pi}_{j})$ 
for any $j\in I_{\pi}$. Then $S_{t}$ is also an intertwiner of $H_{\pi}$ 
and $A$ for any $t\in\R$, in particular we have 
$\sigma^{\varphi}_{t}(A_{\pi})=A_{\pi}$ for any $t\in\R$ 
and we see any element of $A_{\pi}$ is analytic for 
$\{\sigma^{\varphi}_{t}\}_{t\in\R}$.
\end{lem}

\begin{proof}
We shall show that 
$\alpha(S_{t}\xi^{\pi}_{j})
=\sum_{k\in I_{\pi}}S_{t}\xi^{\pi}_{k}\otimes w(\pi)_{k,j}$ 
for any $j\in I_{\pi}$ and $t\in\R$. In fact, we have 
\begin{align*}
\alpha(S_{t}\xi^{\pi}_{j})=&(F_{\pi})_{j,j}^{-it}\alpha(\sigma^{\varphi}_{t}
(S\xi^{\pi}_{j}))\\
=&(F_{\pi})_{j,j}^{-it}
(\sigma^{\varphi}_{t}\otimes\tau_{-t})(\alpha(S\xi^{\pi}_{j})))\\
=&(F_{\pi})_{j,j}^{-it}\sum_{k\in I_{\pi}}
\sigma^{\varphi}_{t}(S\xi^{\pi}_{k})\otimes\tau_{-t}(w(\pi)_{k,j})\\
=&(F_{\pi})_{j,j}^{-it}\sum_{k\in I_{\pi}}
\sigma^{\varphi}_{t}(S\xi^{\pi}_{k})\otimes 
(F_{\pi})^{-it}_{k,k}(F_{\pi})^{it}_{j,j}w(\pi)_{k,j}\\
=&\sum_{k\in I_{\pi}}
S_{t}\xi^{\pi}_{k}\otimes w(\pi)_{k,j}.
\end{align*}

\end{proof}

From this lemma we see that $\sigma^{\varphi}_{t}(A)=A$ 
for $t\in \R$ and the following proposition holds. 

\begin{prop}
Let $\{A,G,\alpha\}$ be an ergodic system and $\varphi$ be the 
invariant state on $A$. Then there exists 
the modular automorphism group $\{\sigma^{\varphi}_t\}_{t\in\R}$ on A. 
\end{prop}

\begin{defn}
Let $\{A,G,\alpha\}$ be an ergodic system and $\pi$ 
be an element of $\widehat{G}$. 
A vector 
$\xi=(\xi_{j})_{j\in I_{\pi}}\in \bigoplus_{j\in I_{\pi}}A_{\pi}$ is 
called a $\pi$-eigenvector if 
$\alpha(\xi_{j})=\sum_{k\in I_{\pi}}\xi_{k}\otimes w(\pi)_{k,j}$ 
for any $j\in I_{\pi}$. The set of $\pi$-eigenvectors is called 
the $\pi$-eigenvector space and denoted by $X_{w(\pi)}$.
\end{defn}

Let $\{A,G,\alpha\}$ and $\{\xi^\lambda\}_{\lambda\in \Lambda}$ be 
a covariant system and a set of its 
$\pi_\lambda$-eigenvectors, respectively. 
A $C^*$-subalgebra generated by $\xi_r^\lambda$ for all $r\in I_{\pi_\lambda}$ 
and $\lambda \in \Lambda$ is denoted by 
$C^*(\{\xi^\lambda\}_{\lambda\in \Lambda})$. 
The quantum group $G$ acts on it invariantly. 
We say that it is a $G$-invariant $C^*$-subalgebra generated 
by $\{\xi^\lambda\}_{\lambda\in \Lambda}$. 
We make a $\pi$-eigenvector space to a Hilbert space 
by defining its inner product with
\[
(\xi\mid\eta)=\sum_{k\in I_{\pi}} \xi_{k} \eta_{k}^*,
\  \mbox{for all}\  \xi,\eta\in X_{w(\pi)}.
\]
In fact, the above right hand side is an element of the fixed point algebra 
and therefore a scalar by the ergodicity of the action. 

We remark eigenvector spaces do not essentially 
depend on the selection of 
unitaries $\{w(\pi)\}_{\pi\in\widehat{G}}$. 
In fact, let us take 
two selections from $\pi$: $w(\pi)$, $w'(\pi)$ 
and let $X_{w(\pi)}$ and $X_{w'(\pi)}$ be the $\pi$-eigenvector spaces 
corresponding $w(\pi)$ and $w'(\pi)$, respectively. 
Then there exists a unitary $u: H_{w(\pi)}\longrightarrow H_{w'(\pi)}$ 
satisfying $u w(\pi)=w'(\pi) u$. 
The unitary $u$ gives a unitary isomorphism 
between $X_{w(\pi)}$ and $X_{w'(\pi)}$ 
by sending $\xi\in X_{w(\pi)}$ to $\bar{u}\xi\in X_{w(\pi)}$, 
where $(\bar{u}\xi)_{k}=\sum_{\ell\in I_\pi}\overline{u_{k,\ell}}\xi_{\ell}$. 

\begin{defn}
Let $\{A,G,\alpha\}$ be an ergodic system and $\pi$ 
be an element of $\widehat{G}$. 
We define the following two operations on eigenvector spaces. 
\begin{enumerate}

\item For any $t\in\R$, we define the linear map 
$U^{w(\pi)}_t:X_{w(\pi)} \longrightarrow X_{w(\pi)}$ by
\[
(U^{w(\pi)}_t\xi)_k=(F_{\pi})^{-it}_{k,k} \sigma^{\varphi}_{t}(\xi_{k})
\quad \mbox{for all} \quad \xi \in X_{\pi}, k\in I_{\pi}.
\]

\item We define the conjugate linear map 
$T_{w(\pi)}:X_{w(\pi)}\longrightarrow X_{\ovl{w(\pi)}}$. 
\[
(T_{w(\pi)}\xi)_{k}=(F_{\pi})^{-\frac{1}{2}}_{k,k}\xi^*_{k}
\quad \mbox{for all} \quad \xi\in X_{\pi}, k\in I_{\pi}.
\]

\end{enumerate}

\end{defn}

The well-definedness of the above first operation 
has been already checked in Lemma \ref{R-action}. 
For the second one we shall justify the equality 
$\alpha((T_{w(\pi)}\xi)_{k})=\sum_{\ell\in I_{\pi}}
(T_{w(\pi)}\xi)_{\ell}\otimes \overline{w(\pi)}_{\ell,k}$ 
for any $\xi\in X_{w(\pi)}$ and $k\in I_{\pi}$. 
In fact, 
\begin{align*}
\alpha((T_{w(\pi)}\xi)_{k})=&(F_{\pi})^{-\frac{1}{2}}_{k,k}\alpha(\xi^*_{k})\\
=&(F_{\pi})^{-\frac{1}{2}}_{k,k}
\sum_{\ell\in I_{\pi}}\xi_{\ell}^*\otimes w(\pi)_{\ell,k}^*\\
=&(F_{\pi})^{-\frac{1}{2}}_{k,k}
\sum_{\ell\in I_{\pi}}\xi_{\ell}^*
\otimes(F_{\pi})^{\frac{1}{2}}_{k,k}
(F_{\pi})^{-\frac{1}{2}}_{\ell,\ell}\overline{w(\pi)}_{\ell,k}\\
=&\sum_{\ell\in I_{\pi}}(T_{w(\pi)}\xi)_{\ell}\otimes 
\overline{w(\pi)}_{\ell,k}.
\end{align*}

\begin{prop}\label{intertwiner}
Let $\{A,G,\alpha\}$ be an ergodic system and $\pi$ 
be an element of $\widehat{G}$. 
Then the map 
$\Hom_{G}(H_{\pi},A)\ni S\longrightarrow 
(S(\xi^{\pi}_k))_{k\in I_{\pi}}\in X_{w(\pi)}$ is a linear isomorphism.
\end{prop}

\begin{proof}
For $\xi\in X_{w(\pi)}$, define the element $S\in \Hom_{G}(H_{\pi},A)$ by 
$S(\xi^{\pi}_{k})=\xi_{k}$ for any $k\in I_{\pi}$. 
It is easy to see that this map is well-defined and gives the 
desired inverse map. 
\end{proof}

From this proposition we see the spectral subspace $A_{\pi}$ is 
spanned by the entries of $\pi$-eigenvectors. 

The next lemma has already appeared in the proof of 
\cite[Proposition 18]{Boca}, 
however, we give a proof for readers' convenience. 

\begin{lem}\label{innerproduct}
Let $\xi$ and $\eta$ be $\pi$-eigenvectors of an ergodic system 
$\{A,G,\alpha\}$. Then the following equalities hold

\begin{enumerate}

\item $\varphi(\xi_{k}\eta_{\ell}^*)=
D_{\pi}^{-1}(F_{\pi})_{k,k}\delta_{k,\ell}(\xi\mid\eta)
$, for all $k,\ell \in I_\pi$,

\item $\varphi(\xi_{k}^*\eta_{\ell})=0$, \ if $k\neq \ell \in I_{\pi}$,
\  and \  
$\varphi(\xi_{k}^*\eta_{k})=\varphi(\xi_{\ell}^*\eta_{\ell})$, 
\  for all $k, \ell \in I_{\pi}$.
\end{enumerate}

\end{lem}

\begin{proof}
(1)It is a straightforward calculation 
\begin{align*}
\varphi(\xi_{k}\eta_{\ell}^*)
=&\sum_{r,s\in I_{\pi}}
\xi_{r}\eta_{s}^* h(w(\pi)_{r,k}w(\pi)_{s,\ell}^*)\\
=&\sum_{r,s\in I_{\pi}}
\xi_{r}\eta_{s}^*  D_{\pi}^{-1}(F_{\pi})_{k,k}\delta_{r,s}\delta_{k,\ell}\\
=&D_{\pi}^{-1}(F_{\pi})_{k,k}\delta_{k,\ell}
\sum_{r\in I_{\pi}}
\xi_{r}\eta_{r}^*\\
=&D_{\pi}^{-1}(F_{\pi})_{k,k}\delta_{k,\ell}(\xi\mid\eta).
\end{align*}
(2) 
Let $Z=(\varphi(\xi_{k}^*\eta_{\ell}))_{k,\ell\in I_{\pi}}$ be 
a matrix on the representation space $H_{\pi}$. Then we have
\begin{align*}
(w(\pi)^* Zw(\pi))_{i,j}
=&\sum_{k,\ell\in I_{\pi}}
w(\pi)_{k,i}^* Z_{k,\ell}w(\pi)_{\ell,j}\\
=&\sum_{k,\ell\in I_{\pi}}
w(\pi)_{k,i}^* \varphi(\xi_{k}^*\eta_{\ell})w(\pi)_{\ell,j}\\
=&\sum_{k,\ell\in I_{\pi}}
(\varphi\otimes\id )
((\xi_{k}^*\otimes w(\pi)_{k,i}^*)( \eta_{\ell}\otimes w(\pi)_{\ell,j}))\\
=&(\varphi\otimes\id )
\alpha(\xi_{i}^*\eta_{j})\\
=&\varphi(\xi_{i}^*\eta_{j})\\
=&Z_{i,j}.
\end{align*}

Therefore the operator $Z$ commutes with the irreducible 
unitary representation $w(\pi)$, so we have $Z\in\C$.  
\end{proof}

\begin{lem}\label{T*T=U}
Let $\{A,G,\alpha\}$ be an ergodic system and 
$\pi$ be an element of $\widehat{G}$. Then it follows 
$T_{w(\pi)}^*T_{w(\pi)}=U_{i}^{w(\pi)}$.
\end{lem}

\begin{proof}
Take $\pi$-eigenvectors $\xi$ and $\eta$. Then we have
\begin{align*}
(T_{w(\pi)}\xi\mid T_{w(\pi)}\eta)
=&\sum_{k\in I_{\pi}}(T_{w(\pi)}\xi)_k (T_{w(\pi)}\eta)^*_k\\
=&\sum_{k\in I_{\pi}}
(F_{\pi})^{-\frac{1}{2}}_{k,k}(F_{\pi})^{-\frac{1}{2}}_{k,k}
\xi_{k}^* \eta_{k}\\
=&\sum_{k\in I_{\pi}}
(F_{\pi})^{-1}_{k,k}\varphi(\xi_{k}^* \eta_{k})\\
=&\sum_{k\in I_{\pi}}
(F_{\pi})^{-1}_{k,k}\varphi(\sigma^{\varphi}_{i}(\eta_{k})\xi_{k}^* )\\
=&\sum_{k\in I_{\pi}}
(F_{\pi})^{-1}_{k,k}
\varphi((F_{\pi})^{-1}_{k,k}(U^{w(\pi)}_{i}\eta)_{k}\xi_{k}^* )\\
=&\sum_{k\in I_{\pi}}
(F_{\pi})^{-2}_{k,k}\cdot D_{\pi}^{-1}(F_{\pi})_{k,k}
(U^{w(\pi)}_{i} \eta \mid \xi)\\
=&(U^{w(\pi)}_{i} \eta \mid \xi),
\end{align*}
where we use the result of the previous lemma 
and $\sum_{k\in I_{\pi}}(F_{\pi}^{-1})_{k,k}=D_{\pi}$. Hence we obtain 
the desired equality.
\end{proof}

Let $T_{w(\pi)}=J_{w(\pi)}|T_{w(\pi)}|$ be the polar decomposition of 
the conjugate linear map 
$T_{w(\pi)}:X_{w(\pi)}\longrightarrow X_{\overline{w(\pi)}}$. 
Since it is easy to see 
$T_{\overline{w(\pi)}}\circ T_{w(\pi)}=\id_{X_{w(\pi)}}$, 
the map $J_{w(\pi)}$ is a conjugate unitary map and satisfies
\[
J_{\overline{w(\pi)}}=J_{w(\pi)}^{-1}, 
\  \mbox{and}\  
J_{\overline{w(\pi)}}|T_{\ovl{w(\pi)}}|J_{w(\pi)}
=\big{|}T_{w(\pi)}\big{|}^{-1}
.\]

In particular we obtain 
$\spe\Big{(}U^{w(\pi)}_{\frac{i}{2}}\Big{)}
=\spe\Big{(}U^{\overline{w(\pi)}}_{\frac{i}{2}}\Big{)}^{-1}$. 
With these preparations we show the following result about 
bounds of 
multiplicities. 

\begin{thm}\label{multiplicity bound theorem}
Let $\{A,G,\alpha\}$ be an ergodic system and 
$\pi$ be an element of $\widehat{G}$. Then we have
\[
\dim\Hom_{G}(H_{\pi},A)\leq D_{\pi}.
\]
\end{thm}

\begin{proof}
The proof is essentially due to \cite[Theorem 17]{Boca} 
or \cite[Theorems 1, 2]{Wassermann1}. 
As we have proved, in Proposition \ref{intertwiner}, 
the hom-space $\Hom_{G}(H_{\pi},A)$ is naturally 
isomorphic to $X_{w(\pi)}$. 
Hence it suffices to show $d:=\dim(X_{w(\pi)})\leq D_{\pi}$. 
The unitary $\R$-action $\{U^{w(\pi)}_{t}\}_{t\in\R}$ enables us 
to take an orthonormal basis of 
$X_{w(\pi)}$ $\{\xi^{p}\}_{1\leq p\leq d}$ 
as $U^{w(\pi)}_{t}\xi^{p}=\lambda_{p}^{it}\xi^{p}$ for all $t\in\R$ 
where $\lambda_{p}$ is a positive real number. 
Then we define an operator entry matrix $M$ by 
\[
M=\begin{pmatrix}
\xi^{1}\\

\xi^{2}\\

\vdots\\

\xi^{d}\\
\end{pmatrix}
\in 
\mathbb{M}_{d,d_{\pi}}(A_{\pi}),
\]
where each $\xi^{p}$ is treated as a row vector. 
Adding $0$ entries, we embed the matrix $M$ 
to a larger square 
matrix (still denoted by $M$) of size 
$n$ $(\geq d, d_{\pi})$, if necessary. 
By orthonormality, we have $MM^*=1_{d}$. 
Let $k,\tilde{F}$ be positive matrices 
$\diag(\lambda_{1}^{-1},\lambda_{2}^{-1},\ldots,\lambda_{d}^{-1},
1,\ldots,1)$ and $\diag(F_{w(\pi)},1,\ldots,1)$ 
in $\mathbb{M}_{n}(\C)_{+}$, respectively. 
Then we have
\begin{align*}
\sum_{p=1}^{d}\lambda_{p}^{-1}
=&\Tr_{k}(1_{d})\\
=&(\varphi\otimes\Tr_{k})(MM^*)\\
=&(\varphi\otimes\Tr_{k})(M^*(\sigma^{\varphi}_{-i}\otimes \Ad(k))(M))\\
=&(\varphi\otimes\Tr_{k})(M^*M\tilde{F}k^{-1})\\
=&(\varphi\otimes\Tr)(M^*M\tilde{F})\\
\leq&(\varphi\otimes\Tr)(1\otimes \tilde{F}1_{d_{\pi}})\\
=&D_{\pi},
\end{align*}
where we use 
$(\sigma^{\varphi}_{t}\otimes \Ad(k^{it}))(M)=M\tilde{F}^{it}k^{-it}$ 
for $t\in\R$. 
Hence we obtain $\Tr_{X_{w(\pi)}}(U^{w(\pi)}_{i})
\leq D_{\pi}$ 
and similarly $\Tr_{X_{\ovl{w(\pi)}}}(U^{\ovl{w(\pi)}}_{i})
\leq D_{\ovl{\pi}}=D_{\pi}$. 
Since we have already known 
$\spe\Big{(}U^{w(\pi)}_{i}\Big{)}
=\spe\Big{(}U^{\overline{w(\pi)}}_{i}\Big{)}^{-1}$, 
it follows 
$2d\leq \Tr_{X_{w(\pi)}}(U^{w(\pi)}_{i})+
\Tr_{X_{\overline{w(\pi)}}}(U^{\overline{w(\pi)}}_{i})$. 
Therefore we obtain $2d\leq 2 D_{\pi}$. 
\end{proof}

From this result, the dimension of the $\pi$-eigenspace $A_{\pi}$ is 
less than or equal to $d_{\pi}D_{\pi}$. 
We can derive 
the quantum version of \cite[Theorem 2]{Wassermann1}. 

\begin{thm}\label{pAq multiplicity}
Let $\{A,G,\alpha\}$ be a compact quantum covariant system 
and $\pi$ be an element of $\widehat{G}$. 
If $p$ and $q$ are  minimal projections in $A^\alpha$, then we have:
\[
\dim\Hom_{G}(H_{\pi},pAq)\leq D_{\pi}.
\]
\end{thm}

\begin{proof}
This proof is also essentially due to \cite[Theorem 17]{Boca} 
or \cite[Theorems 1, 2]{Wassermann1}, 
however, we cannot use the trace property as in the case of 
a compact group, we have to prove the finite dimensionality of 
$\Hom_{G}(H_{\pi},pAq)$ at first. 
If the linear space $pA^{\alpha}q$ is not $0$, 
we can take non-zero norm $1$ element $x$ from $pA^{\alpha}q$. 
Then we see $xx^*\in pA^{\alpha}p=\C p$ and $x^*x\in qA^{\alpha}q=\C q$ 
and have $xx^*=p$ and $x^*x=q$. 
It follows $(pAq)_{\pi}=(pAp)_{\pi}x$. 
Therefore the assertion of this theorem 
follows from applying the previous theorem to $pAp$. 
Next we assume $pA^\alpha q=0$. 
Because $p\perp q$, we may assume that $A$ is unital 
and $1=p+q$ by considering $(p+q)A(p+q)$. 
Let $\varphi(x)=\varphi_p(pxp) + \varphi_q(qxq)$ 
where $\varphi_p$ and $\varphi_q$ are the invariant states on $pAp$ 
and $qAq$ respectively. 
Since the conditional expectation with respect to $\varphi$ 
$E_{\alpha}:A\longrightarrow A^\alpha$ is 
given by $E_{\alpha}(x)=pxp+qxq$, 
$\varphi$ is the invariant positive functional on $A$. 
Let $d\leq\infty$ be 
the dimension of the linear space $\Hom_{G}(H_{\pi},A)$ and 
fix a natural number $0\leq d'\leq d$. 
Note that we can take $\pi$-eigenvectors 
$\{\xi^{j}\}_{1\leq j\leq d'}$ 
from $(pAq)_{\pi}=pA_{\pi}q$, 
which satisfy 
$\varphi((\xi^{j}_{r})^* \xi^{k}_{s})=\delta_{j,k}\delta_{r,s}$ 
for $1\leq j,k \leq d'$ and $r,s \in I_{\pi}$. 
In fact, if we have $\{\xi^{j}\}_{1\leq j \leq k-1}$ as before, 
we can take a $\pi$-eigenvector $\xi^{k}$ which satisfies 
$\varphi((\xi^{j}_{1})^*\xi^{k}_{1})=\delta_{j,k}$ for $1\leq j \leq p-1$. 
Then modifying Lemma \ref{innerproduct} to the case of 
$pAq$, 
we achieve to take desired vectors. 
Let $M$ be the matrix
\[
M
=
\begin{pmatrix}
T_{w(\pi)}\xi^{1} \\
T_{w(\pi)}\xi^{2} \\
\vdots\\
T_{w(\pi)}\xi^{d'} 
\end{pmatrix}.
\]
Adding $0$-entries, we embed $M$ to a larger square matrix 
(still denoted by $M$) of size $n$ $(\geq d',d_{\pi})$. 
Note that $T_{w(\pi)}\xi^{j}$ is a $\bar{\pi}$-eigenvector for $qAp$:
$T_{w(\pi)}\xi^{j}=\big{(}
(F_{\pi})^{-\frac{1}{2}}_{r,r}{\xi^{j}_{r}}^{*}
\big{)}_{r\in I_{\pi}}$. 
By the computation in the proof of Lemma \ref{T*T=U}, 
it follows $MM^*=D_{\pi}q\otimes1_{d'}$. 
Hence we have $M^*M\leq D_{\pi}p\otimes1_{d_{\pi}}$ 
in $pAp\otimes \mathbb{B}(H_{\pi})$. 
Then it yields:
\begin{align*}
D_{\pi}^2
\geq&(\varphi\otimes\Tr_{F_{\pi}})(M^*M)\\
=&\sum_{r\in I_{\pi}}
\big{(}
(F_{\pi})_{r,r}(F_{\pi}^{-\frac{1}{2}})_{r,r}
(F_{\pi}^{-\frac{1}{2}})_{r,r}
\sum_{1\leq j\leq d'}
\varphi(\xi^{j}_{r}{\xi^{j}_{r}}^*)
\big{)}\\
=&\sum_{1\leq j\leq d'}(\xi^{j}\mid \xi^{j}), 
\end{align*}
where $(\xi^{j}\mid \xi^{j})$ means the inner product of $\pi$-eigenvectors 
for $pAq$. 
Let us define another inner product $(\cdot\mid\cdot)_{2}$ 
of $\pi$-eigenvectors for $pAq$ by 
$(\xi\mid\eta)_{2}=\varphi(\eta_{r}^*\xi_{r})$ 
for any $\pi$-eigenvectors $\xi,\eta$. 
Note that this value does not depend on the choice of $r$ 
by Lemma \ref{innerproduct}. 
Let $W$ be a linear space spanned by $\{\xi^{j}\}_{1\leq j\leq d'}$. 
There are two inner products $(\cdot\mid\cdot)$ 
and $(\cdot\mid\cdot)_{2}$ on $W$. 
Let $A_{W}$ be the matrix which satisfies 
$(A_{W}\xi\mid\eta)_{2}=D_{\pi}^{-1}(\xi\mid\eta)$. 
Then from the above inequality we have 
$\Tr_{W}(A_{W})\leq D_{\pi}$, 
where $\Tr_{W}$ is the non-normalized trace 
associated to $W$. 
Similarly considering $T_{w(\pi)}(W)$ for $W$, 
we also obtain 
$\Tr_{T_{w(\pi)}(W)}(A_{T_{w(\pi)}(W)})\leq D_{\pi}$. 
Moreover, we see $A_{W}=T_{w(\pi)}^*T_{w(\pi)}$, 
where the involution $*$ comes from 
the conjugate linear map 
$T_{w(\pi)}:W\longrightarrow T_{w(\pi)}(W)$ 
between $(\cdot\mid\cdot)_{2}$-inner product spaces.  
In fact, we have
\begin{align*}
(T_{w(\pi)}\xi\mid T_{w(\pi)}\eta)_{2}
=&\varphi({(T\eta)_{1}}^*(T\xi)_{1})\\
=&(F_{\pi}^{-1})_{1,1}\varphi(\eta_{1} \xi_{1}^*)\\
=&(F_{\pi}^{-1})_{1,1}D_{\pi}^{-1}(F_{\pi})_{1,1}(\eta\mid\xi)\\
=&D_{\pi}^{-1}(\eta\mid\xi).
\end{align*}
Let us denote $T_{w(\pi)}(W)$ by $\overline{W}$.
Let $T_{w(\pi)}=J_{W}|T_{w(\pi)}|$ 
and $T_{\overline{w(\pi)}}
=J_{\overline{W}}|T_{\overline{w(\pi)}}|$ be 
the polar decomposition of the conjugate linear maps 
$T_{w(\pi)}$ and $T_{\overline{w(\pi)}}$ 
between $W$ and $T_{w(\pi)}(W)$. 
The equality $T_{w(\pi)}
\circ T_{\overline{w(\pi)}}=\id_{\overline{W}}$ 
yields 
$J_{W}A_{W}^{\frac{1}{2}}J_{\overline{W}}=A_{\overline{W}}^{-\frac{1}{2}}$. 
Then we have $J_{W}A_{W}J_{W}^*=A_{\overline{W}}^{-1}$ 
because $J_{W}$ and $J_{\overline{W}}$ are conjugate unitary. 
Hence we obtain $\Tr(A_{W})=\Tr(A_{\overline{W}}^{-1})$ 
and it follows $2d'\leq2\Tr(A_{W})\leq 2D_{\pi}$. 
\end{proof}

Finally we state the modified version of 
Theorem \ref{multiplicity bound theorem} and Theorem \ref{pAq multiplicity}. 

\begin{prop}
Let $\{A,G,\alpha\}$ be an ergodic system 
and $\pi$ be an element of $\widehat{G}$. 
Assume that $A$ has a tracial state. 
Then we have:
\[
\dim\Hom(H_{\pi},A)\leq d_{\pi}.
\]
\end{prop}

\begin{proof}
We have proved finite dimensionality of $X_{w(\pi)}$. 
Let $d$ be its dimension and 
$\{\xi^{p}\}_{1\leq p \leq d}$ be an orthonormal basis of $X_{w(\pi)}$ 
for the inner product $(\cdot\mid\cdot)$. 
Take a tracial state $\tau$ on $A$. 
Define the matrix $M$ by
\[
M
=
\begin{pmatrix}
\xi^{1}\\
\xi^{2}\\
\vdots\\
\xi^{d}
\end{pmatrix}. 
\]
As the proof of Theorem \ref{multiplicity bound theorem}, 
we look at $M$ in a bigger square matrix algebra. 
Then we have $MM^*=1_{d}$ and it yields
\begin{align*}
d_{\pi}
\geq&(\tau\otimes\Tr)(M^*M)\\
=&(\tau\otimes\Tr)(MM^*)\\
=&d.
\end{align*}

\end{proof}

\begin{prop}
Let $\{A,G,\alpha\}$ be a compact quantum covariant system 
and $\pi$ be an element of $\widehat{G}$. 
Assume that $A$ has a tracial state. 
If $p$ and $q$ are minimal projections in $A^\alpha$ 
and they are not annihilated by a trace on $A$, 
then we have
\[
\dim\Hom(H_{\pi},pAq)\leq d_{\pi}.
\]
\end{prop}

\begin{proof}
Let $d (\leq D_{\pi})$ be the dimension of 
$\dim\Hom(H_{\pi},pAq)$. 
Let $\{\xi^{p}\}_{1\leq p \leq d}$ be an orthonormal $\pi$-eigenvector 
for $pAq$ where we define the inner product by
\[
(\xi\mid\eta)p=\sum_{r\in I_{\pi}}\xi_{r}\eta_{r}^*
\]
for all $\pi$-eigenvector $\xi$ and $\eta$. 
Take a tracial state $\tau$ on A with $\tau(p)\tau(q)\neq 0$. 
Define the matrix $M$ by
\[
M
=
\begin{pmatrix}
\xi^{1}\\
\xi^{2}\\
\vdots\\
\xi^{d}
\end{pmatrix}. 
\]
and embed it into the larger square matrix as in previous theorems. 
Then we have $MM^*=p\otimes1_{d}$. This yields
\begin{align*}
d\,\tau(p)
=&(\tau\otimes\Tr)(MM^*)\\
=&(\tau\otimes\Tr)(M^*M)\\
\leq&(\tau\otimes\Tr)(q\otimes 1_{d_{\pi}})\\
=&\tau(q)\,d_{\pi}.
\end{align*}
Hence we have $d\leq d_{\pi} \frac{\tau(q)}{\tau(p)}$. 
Considering $qAp=(pAq)^*$ and $d=\dim\Hom(H_{\pi},qAp)$, 
we also have 
$d\leq d_{\ovl{\pi}} \frac{\tau(p)}{\tau(q)}
=d_{\pi} \frac{\tau(p)}{\tau(q)}$. 
Finally we obtain 
$d\leq d_{\pi} \min\big{(}\frac{\tau(q)}{\tau(p)}
,\frac{\tau(p)}{\tau(q)}\big{)}\leq d_{\pi}$. 

\end{proof}

\begin{rem}\label{tau(p)>0}
Assume that the linear space $pAq$ is not $0$. 
Then $\tau(p)\neq0$ if and only if $\tau(q)\neq0$. 
Indeed, take $\pi\in\widehat{G}$ as $pA_{\pi}q\neq0$ and  
assume $\tau(q)=0$, 
then 
we obtain $\tau(p)=0$ by the proof of the above proposition. 
The converse assertion holds by applying the involution $*$ from 
$pAq$ to $qAp$. 
\end{rem}

\section{Equivariant $K$-theory}

We follow \cite{Wassermann1} and use 
equivariant $K$-theory to obtain 
a multiplicity map. 
We briefly recall the notions. 
Readers are referred to, for example, \cite{BaajSkandalis1} and \cite{Vergnioux} for the theory of equivariant $K$-theory. 
Let $\{A,G,\alpha\}$ be a compact quantum covariant system 
with $A$ unital. 
Let $E$ be a finitely generated projective Hilbert $A$-module 
and $\delta_E:E\longrightarrow E\otimes C(G)$ be a 
linear map. 
The pair $\{E,\delta_E\}$ is called a $G$-equivariant $A$-module 
if they satisfy the following
\begin{enumerate}

\item $(\id\otimes\delta)\circ\delta_E=(\delta_E\otimes\id)\circ\delta_E\, ,$

\item $\delta_E(ea)=\delta_E(e)\alpha(a)$ for all $e\in E$ and $a\in a\, ,$

\item $\langle\delta_E(e)\mid \delta_E(e')\rangle
=\alpha(\langle e\mid e'\rangle)$ for all $e,e'\in E\, ,$

\item the linear subspace $\delta_E(E)(\C\otimes C(G))$ is 
      dense in $E\otimes C(G)$. 
\end{enumerate}

Actually the second equality follows from the third one. 
Inner products of Hilbert modules are conjugate-linear for 
the first variable. 
Two $G$-equivariant $A$-modules are equivalent if 
there exists 
a bijective linear map intertwining 
of the actions of $G$ and $A$. 
The set of these equivalence classes becomes an 
abelian semigroup by the direct sum 
and its Grothendieck group is denoted by 
$K_0^G(A)$ and 
this is called an \textit{equivariant K-group}. 
For a $G$-equivariant $A$-module $\{E,\delta_E\}$ 
we define its crossed product Hilbert $A\rtimes_\alpha G$-module 
by $E\rtimes_\alpha G=E\otimes_A (A\rtimes_\alpha G)$. 
This module becomes a (left) $C_r^*(G)$-module by an appropriate way 
(\cite[Lemme 5.2]{Vergnioux}). 
Here we take a non-rigorous picture for its action. 
Let $x$ be an element of $C_r^*(G)$ and we assume its dual coproduct 
has an expansion $\hat{\delta}(x)=\sum x_{(0)}\otimes x_{(1)}$. 
It is considered that $C_r^*(G)$ acts on $E$ through the 
differential representation with respect to $\delta_E$. 
Hence we have an action $x\cdot(e\otimes_A y)=\sum x_{(0)}(e)
\otimes_A x_{(1)}y$. 
This action is compatible with the right action of $A$, 
because in $A\rtimes_\alpha G$ 
we have $z (\alpha(y))=\sum (z_{(0)}\cdot y) z_{(1)}$ 
for $z\in C_r^*(G)$ and $y\in A$ 
where $C_r^*(G)$ acts on $A$ with the differential representation 
about $\alpha$. 
We next recall the Green-Julg isomorphism. 
If $E_\pi$ is a tensor product module $H_\pi\otimes A$ with 
a finite dimensional irreducible $G$-module $H_\pi$, 
then we have 
$E_\pi\otimes_A (A\rtimes_\alpha G)=H_\pi\otimes (A\rtimes_\alpha G)$. 
We apply $p_0\in C_r^*(G)$ to this space. 
Since this action is obtained by the dual coproduct, we have 
$p_0(E_\pi\otimes_A A\rtimes_\alpha G)
=p_0(H_\pi\otimes C_r^*(G))\otimes_{C_r^*(G)}A\rtimes_\alpha G
\cong p_{\overline{\pi}}C_r^*(G)\otimes_{C_r^*(G)}A\rtimes_\alpha G
=(1\otimes p_{\overline{\pi}})A\rtimes_\alpha G\, .$ 
In general case $E$ is a direct summand of a direct sum of $E_\pi$. 
Hence the module $p_0(E\rtimes_\alpha G)$ corresponds to 
a projection in $M_n(A\rtimes_\alpha G)$ for some $n$. 
We also notice that $A\rtimes_\alpha G$ is stably unital. 
In fact we have an approximate unit of $\mathbb{K}(L^2(G))$ 
which consists of projections in $C_r^*(G)$. 
In this way we obtain the Green-Julg isomorphism 
(see \cite[Th\'{e}or\`{e}me 5.10]{Vergnioux} for its proof of $KK$-version):
\[
\Phi_1:K_0^G(A)\longrightarrow K_0(A\rtimes_\alpha G) 
\] 
defined by $[E]\mapsto [p_0(E\rtimes_\alpha G)]$. 
The inverse map of $\Phi_1$ is given by
\[
\Phi_2:K_0(A\rtimes_\alpha G)\longrightarrow K_0^G(A)
\] 
which sends $[q]$ to $[q (A\otimes L^2(G))^n]$, 
where $q\in M_n(A\rtimes_\alpha G)$ is a projection and 
$A\otimes L^2(G)$ is a $G$-equivariant $A$-module with 
$\delta_{A\otimes L^2(G)}(a\otimes \xi)
=(1\otimes W(\xi\otimes1)) \alpha(a)_{13}$. 
We define the usual $R(G)$-module structures on $K_0^G(A)$ and 
$K_0(A\rtimes_\alpha G)$. 
Let $\pi$ be an element of $\widehat{G}$ and $E$ be a 
$G$-equivariant $A$-module 
and $q\in M_n(A\rtimes_\alpha G)$ be a projection. 
We define the action of $\pi$
\begin{align*}
\pi\cdot [E]=&\,[H_\pi\otimes E]\, ,\\
\pi\cdot [q]=&\,[(\id_{M_n(\C)}\otimes\hat{\alpha}_\pi)(q)]\, ,
\end{align*}
where we use the isomorphism 
$K_0(A\rtimes_\alpha G)\cong K_0(A\rtimes_\alpha G \otimes \End(H_\pi))$. 
We observe the Green-Julg isomorphism is an $R(G)$-module map in the 
compact group case. 

\begin{lem}
Let $B$ be a $C^*$-algebra and 
$E$ be a Hilbert $B$-module which has 
a $C_r^*(G)$-action 
$\phi:C_r^*(G)\longrightarrow \mathbb{B}(E)$. 
Then for any $\pi\in \widehat{G}$ there is an isomorphism as 
Hilbert $B$-modules 
between $p_0 (H_\pi \otimes E)$ and $p_0(E \otimes H_\pi)$. 
\end{lem}

\begin{proof}
First we observe the difference of $\hat{\delta}(p_0)$ and its 
opposite $\hat{\delta}^{op}(p_0)$. 
Using $p_0=(\id\otimes h)(V)$ and the pentagonal 
identity, we obtain 
$\hat{\delta}(p_0)=(\id\otimes\id\otimes h)(V_{13}V_{23})$ and 
$\hat{\delta}^{op}(p_0)=(\id\otimes\id\otimes h)(V_{23}V_{13})$. 
Since $V$ has the expansion 
$\sum_{\pi\in \widehat{G},i,j\in I_\pi}E_{i,j}^\pi\otimes w(\pi)_{i,j}$, 
we have the following equalities:
\begin{align*}
\hat{\delta}(p_0)
=&\, 
\sum_{\pi\in \widehat{G},i,j\in I_\pi,k,\ell\in I_{\ovl{\pi}}}
h(w(\pi)_{i,j}w(\overline{\pi})_{k,\ell})
E_{i,j}^{\pi}\otimes E_{k,\ell}^{\overline{\pi}} \, ,\\
\hat{\delta}^{op}(p_0)
=&\, 
\sum_{\pi\in \widehat{G},i,j\in I_\pi,k,\ell\in I_{\ovl{\pi}}}
h(w(\ovl{\pi})_{k,\ell}w(\pi)_{i,j})
E_{i,j}^{\pi}\otimes E_{k,\ell}^{\overline{\pi}}\, .\\
\end{align*}
We use 
$\sigma_{-i}^h(w(\ovl{\pi})_{k,\ell})
=(F_{\ovl{\pi}})_{k,k}(F_{\ovl{\pi}})_{\ell,\ell}w(\ovl{\pi})_{k,\ell}$ 
in order to get 
$h(w(\ovl{\pi})_{k,\ell}w(\pi)_{i,j})
=
(F_{\ovl{\pi}})_{k,k}(F_{\ovl{\pi}})_{\ell,\ell}
h(w(\pi)_{i,j}w(\overline{\pi})_{k,\ell})\, .
$
Define a densely defined unbounded positive operator 
$F=\sum_{\pi\in \widehat{G}} F_\pi\, .$ 
With this operator we get the identity
\[
\hat{\delta}^{op}(p_0)
=
(1\otimes F)\hat{\delta}(p_0)(1\otimes F)\, .
\]
Then we define the map 
$\chi:H_\pi \otimes E\longrightarrow E\otimes H_\pi$ 
by 
$\chi(\xi\otimes e)=e\otimes F_\pi^{-1}\xi\, .$ 
This is a bijective adjointable map for Hilbert $B$-modules. 
Then we get
\begin{align*}
\chi\circ (\id\otimes\phi)(\hat{\delta}(p_0))(\xi\otimes e)
=&\,
(1\otimes F_\pi^{-1})(\phi\otimes\id)
(\hat{\delta}^{op}(p_0))(e\otimes \xi)\\
=&\,
(\phi\otimes\id)(\hat{\delta}(p_0))(e\otimes F_\pi\xi)\, .
\end{align*}
Hence we see 
$\chi$ maps $p_0(H_\pi\otimes E)$ onto $p_0(E\otimes H_\pi)$. 
\end{proof}

\begin{lem}\label{R(G)-module map}
The Green-Julg isomorphism $\Phi_1$ is an $R(G)$-module map. 
\end{lem}

\begin{proof}
Take $\pi\in \widehat{G}$ and a $G$-equivariant $A$-module $E$. 
We shall prove 
$\Phi_1(\pi\cdot[E])=(\hat{\alpha}_\pi)_*(\Phi_1([E]))$. 
For the left side we have 
\[
[p_0((H_\pi\otimes E)\otimes_A A\rtimes_\alpha G )]
=
[p_0(H_\pi\otimes E\rtimes_\alpha G)]\, .
\]
For the right one we have 
\begin{align*}
(\hat{\alpha}_\pi)_*([p_0(E\rtimes_\alpha G)])
=
&\,
[p_0(E\rtimes_\alpha G)
\otimes_{\hat{\alpha}_\pi}
\big{(}A\rtimes_\alpha G \otimes \End(H_\pi)\big{)}]\\
=
&\,
[p_0(E\rtimes_\alpha G\otimes \End(H_\pi))]\, .
\end{align*}
Moving it to $K_0(A\rtimes_\alpha G)$ by 
tensoring $A\rtimes_\alpha G\otimes H_\pi$ over 
$A\rtimes_\alpha G\otimes \End(H_\pi)$, 
we get: 
$(\hat{\alpha}_\pi)_*(\Phi_1([E]))
=
[p_0(E\rtimes_\alpha G\otimes H_\pi)]\, .$ 
Hence we obtain the desired equality by applying 
the previous lemma to $E\rtimes_\alpha G$ and $B=A\rtimes_\alpha G$. 
\end{proof}

We consider a compact quantum ergodic system $\{A,G,\alpha\}$. 
Then there exists an index set $I$ and Hilbert spaces $H_i$ for 
each $i\in I$ such that 
the crossed product $A\rtimes_{\alpha}G$ is isomorphic to 
$\bigoplus_{i\in I}\mathbb{K}(H_i)$, where $\{H_i\}_{i\in I}$ 
are Hilbert spaces (see \cite[Theorem 19]{Boca}, 
\cite[Corollary 2]{Wassermann1} for its proof). 
Let us fix minimal projections $e_i \in \mathbb{K}(H_i)$ for all $i\in I$. 
Hence we have 
$K_0^G(A)\cong \bigoplus_{i\in I}\Z [e_i]$ by the Green-Julg theorem.  
Using the isomorphism 
$\Phi_2 :K_0(A\rtimes_\alpha G)\longrightarrow K_0^G(A)$, 
we can easily check that 
in $K_0^G(A)$, $[e_i]$ becomes a $G$-equivariant $A$-module 
$[e_i(A\otimes L^2(G))]\, .$ 
Now we consider the $R(G)$-module structure of $K_0^G(A)$. 
Let $\pi$ be an element of $\widehat{G}$. 
We define the (not necessarily finite size) matrix 
$\mathbb{M}(\pi)=(\mathbb{M}(\pi)_{i,j})_{i,j\in I}$ by
\[
\pi\cdot [e_j]=\sum_{i\in I}\mathbb{M}(\pi)_{i,j}[e_i].
\] 
This equality holds in $K_0^G(A)\cong K_0(A\rtimes_\alpha G)$. 
If we treat $[e_i]=[e_i (A\otimes L^2(G))]$ in $K_0^G(A)$, 
we get the equality 
$\mathbb{M}(\pi)_{i,j}
=
\dim\Hom_{(G,A)}(H_{\pi}\otimes e_j(A\otimes L^2(G)),
e_i(A\otimes L^2(G)))$. 

\begin{lem}
We have an isomorphism between linear spaces 
$\Hom_{(G,A)}(H_{\pi}\otimes e_j(A\otimes L^2(G)),
e_i(A\otimes L^2(G)))$ 
and 
$
\Hom_{G}(H_\pi,e_i(A\otimes \mathbb{K}(L^2(G)))e_j)\,. 
$
\end{lem}

\begin{proof}
We may assume projections $e_i$ and $e_j$ are 
majored by $1\otimes p$ where $p$ is a projection in 
$C_r^*(G)\subset \mathbb{K}(L^2(G))$. 
Let $H_0$ be a closed subspace $p L^2(G)$. 
It is finite dimensional and $G$-invariant, that is, 
$W(H\otimes \C)\subset H\otimes C(G)$. 
Let $\{\eta_p\}_{p\in J}$ be an orthonormal basis of $H_0$. 
With this basis 
we give a matrix representation $(w_{p,q})_{p,q\in J}$ of $W$ by 
$W(\eta_p\otimes1)=\sum_{q\in J}\eta_q\otimes w_{q,p}\, .$
From now we write $\Gamma_0$ and $\Gamma_1$ for 
$C(G)$-comodule maps 
$H_{\pi}\otimes A\otimes H_0\longrightarrow 
H_{\pi}\otimes A\otimes H_0\otimes C(G)$ and 
$A\otimes H_0\longrightarrow A\otimes H_0\otimes C(G)$, 
respectively. 
Take an intertwiner 
$S$ from 
$\Hom_{(G,A)}(H_{\pi}\otimes A\otimes H_0, A\otimes H_0)$. 
Then we can choose $a_{k, q p} \in A$ for all $k\in I_\pi$ 
and $p,q\in J$ such that they satisfy 
$S(\xi_k\otimes a \otimes \eta_p)=\sum_{q\in J}a_{k, q p}a\otimes \eta_q$ 
for all $a\in A$ by $A$-linearity. 
Since $S$ is a $G$-homomorphism, we have 
$\Gamma_1S(\xi_k\otimes 1 \otimes \eta_p)
=
(S\otimes1)(\Gamma_0(\xi_k\otimes 1 \otimes \eta_p))$. 
The the left hand side is equal to
\begin{align*}
\sum_{q\in J}\Gamma_1(a_{k, q p}\otimes \eta_q)
=&\,
\sum_{q,r\in J}(1\otimes\eta_r \otimes w_{rq})
\alpha(a_{k, q p})_{13}\, .
\end{align*}
The right hand side is equal to
\begin{align*}
\sum_{\ell\in I_\pi, s\in J}
(S\otimes1)(\xi_{\ell}\otimes1\otimes \eta_s\otimes w(\pi)_{\ell k}w_{sp})
=
\sum_{\ell \in I_\pi, r,s\in J}
a_{\ell,rs}\otimes \eta_r\otimes w(\pi)_{\ell k}w_{sp}\, .
\end{align*}
Hence we get
\begin{equation}\label{equiv01}
\sum_{q\in J}(1\otimes w_{rq})
\alpha(a_{k, q p})
=
\sum_{\ell\in I_\pi, s\in J}
a_{\ell,rs}\otimes w(\pi)_{\ell k}w_{sp}
\, .
\end{equation}
Next take an intertwiner $T$ from 
$\Hom_{G}(H_\pi,A\otimes \mathbb{B}(H_0))$. 
Then we can choose $b_{k,q p}\in A$ for all 
$k\in I_\pi$ and $p,q\in J$ such that they satisfy
$T\xi_k=\sum_{q,p\in J}b_{k,q p}\otimes \eta_q \odot \eta_p\, .$
Since $T$ is a $G$-homomorphism, we have 
$\tilde{\alpha}\circ T(\xi_k)=(T\otimes1)\circ w(\pi)(\xi_k\otimes1)$ for 
all $k\in I_\pi$. 
The left hand side is equal to
\begin{align*}
\sum_{q,t\in J}\tilde{\alpha}(b_{k,q t}\otimes \eta_q \odot \eta_t)
=
&\,
\sum_{q,t\in J}
W_{23}\alpha(b_{k,q t})_{13}(1\otimes1\otimes \eta_q \odot \eta_t)W_{23}^*\\
=
&
\sum_{q,r,s,t\in J}
(1\otimes\eta_{r}\odot\eta_{s}\otimes1)
(1\otimes1\otimes w_{rq})\alpha(b_{k,qt})_{13}(1\otimes1\otimes w_{st}^*)
\, .
\end{align*}
The right hand side is equal to
\begin{align*}
\sum_{\ell\in I_\pi}
T\xi_{\ell} \otimes w(\pi)_{\ell k}
=
&\,
\sum_{\ell,r,s}
b_{\ell,rs}\otimes \eta_{r}\odot\eta_{s} \otimes w(\pi)_{\ell k}
\, .
\end{align*}
Hence we obtain the following equality
\begin{equation*}
\sum_{q,t\in J}
(1\otimes w_{rq})\alpha(b_{k,qt})(1\otimes w_{st}^*)
=
\sum_{\ell}
b_{\ell,rs}\otimes w(\pi)_{\ell k}\, .
\end{equation*}
Multiplying $w_{sp}$ and summing up with $s$, 
we get
\begin{equation}\label{equiv02}
\sum_{q\in J}
(1\otimes w_{rq})\alpha(b_{k,qp})
=
\sum_{\ell\in I_\pi,s\in J}
b_{\ell,rs}\otimes w(\pi)_{\ell k}w_{sp}\, .
\end{equation}
Now we construct maps 
$\nu:\Hom_{(G,A)}(H_{\pi}\otimes A\otimes H_0,A\otimes H_0)
\longrightarrow
\Hom_{G}(H_\pi,A\otimes \mathbb{B}(H_0) )
$ 
and 
$\nu^{-1}:\Hom_{G}(H_\pi,A\otimes \mathbb{B}(H_0) )
\longrightarrow
\Hom_{(G,A)}(H_{\pi}\otimes A\otimes H_0,A\otimes H_0)
$ 
by the following equalities. For all $k\in I_\pi$ and $p\in J$,
\begin{align*}
\nu(S)(\xi_k)
=
&\,
\sum_{q,p\in J}
a_{k,qp}\otimes \eta_{q}\odot\eta_{p}\, ,\\
\nu^{-1}(T)(\xi_k\otimes a\otimes \eta_p)
=
&\,
\sum_{q\in J} b_{k,qp}a\otimes \eta_q
\, .
\end{align*}
They are actually linear operators 
in each hom-space because equalities (\ref{equiv01}) and 
(\ref{equiv02}) have the same form. 
Since clearly we have 
$\nu\circ\nu^{-1}=\id$ and $\nu^{-1}\circ\nu=\id$, 
they give an isomorphism between 
$\Hom_{(G,A)}(H_{\pi}\otimes A\otimes H_0,A\otimes H_0)$
and
$
\Hom_{G}(H_\pi,A\otimes \mathbb{B}(H_0) )
$. 
Then we also easily check that $\nu$ maps the subspace 
$
\Hom_{(G,A)}(H_{\pi}\otimes e_j(A\otimes L^2(G)),
e_i(A\otimes L^2(G)))
\subset
\Hom_{(G,A)}(H_{\pi}\otimes A\otimes H_0,A\otimes H_0)
$ 
to
the subspace
$
\Hom_{G}(H_\pi,e_i(A\otimes \mathbb{K}(L^2(G)))e_j)
\subset
\Hom_{G}(H_\pi,A\otimes \mathbb{B}(H_0) )\, .
$
\end{proof}

Therefore, as in the classical case 
we can derive the following result by A. Wassermann 
(\cite[p.304]{Wassermann1}). 
\begin{cor}\label{reduced space characterization}
Let $\{A,G,\alpha\}$ be a compact quantum ergodic system and 
$\mathbb{M}:R(G)\longrightarrow M_{|I|}(\Z)$ be the 
multiplicity map. 
For all $i,j$ in the index set $I$, we have
\[\mathbb{M}(\pi)_{i,j}=\dim 
\Hom(H_\pi,e_i (A\otimes\mathbb{K}(L^2(G)))e_j).
\] 
\end{cor}

It yields the identity 
$\mathbb{M}(\ovl{\pi})_{i,j}=\mathbb{M}(\pi)_{j,i}$ 
for all $\pi\in \widehat{G}$ and $i,j\in I$ 
and therefore the multiplicity map 
$\mathbb{M}:R(G)\longrightarrow M_{|I|}(\Z)$ 
is a $*$-homomorphism. 
As we have seen, some 
properties of multiplicity maps also hold in the quantum case. 
However we have to be much careful about the existence of 
its common eigenvector (see, \cite[Theorem17]{Wassermann1} 
for classical case). 
In order to study this problem we need several lemmas. 

\begin{lem}\label{not 0}
For all $i,j\in I$, the reduced $G$-space 
$e_i (A\otimes\mathbb{K}(L^2(G)))e_j$ is not $0$. 
\end{lem}

\begin{proof}
Let 
$M$ be the $\sigma$-weak closure of $A$ in $\mathbb{B}(H_\varphi)$ 
where $H_\varphi$ is the GNS Hilbert space of the invariant state $\varphi$. 
Define the relation $\sim$ in the index set $I$ as 
$i\sim j$ if and only if $e_i (A\otimes\mathbb{K}(L^2(G)))e_j \neq 0$. 
We show this is an equivalence relation. 
A non-trivial part is transitivity of $\sim$. 
Assume $i\sim j$ and $j\sim k$. 
Then there exists $\pi,\rho \in \widehat{G}$ which satisfy 
$\mathbb{M}(\pi)_{i,j}>0$ and $\mathbb{M}(\rho)_{j,k}>0$. 
Then we have 
$\mathbb{M}(\pi\cdot\rho)_{i,k}\geq 
\mathbb{M}(\pi)_{i,j}\,\mathbb{M}(\rho)_{j,k}>0$ 
and 
the left hand side is equal to 
$\sum_{\tau\in\widehat{G}}N_\tau^{\pi\,\rho}\,\mathbb{M}(\tau)_{i,k}$. 
Hence there exists $\tau\in\widehat{G}$ such that 
$\mathbb{M}(\tau)_{i,k}>0$. 
Therefore the reduced space 
$e_i (A\otimes\mathbb{K}(L^2(G)))e_k$ is 
not $0$, that is, $i\sim k$. 
Let $I=\bigsqcup_{\lambda\in\Lambda}I_\lambda$ be 
the decomposition with respect to the equivalence relation $\sim$. 
Let $z_i$ be the central projection in $Z(M)$ with 
$z(e_i)=z_i\otimes1$ for all $i\in I$ 
where $z(e_i)$ means the central support projection of $e_i$ in 
$M\otimes \mathbb{B}(L^2(G))$. 
Since for $\lambda\neq\mu\in \Lambda$ we have 
$e_{i} (A\otimes\mathbb{K}(L^2(G)))e_{j}=0$ 
for all $i\in I_\lambda$ and $j\in I_\mu$, 
we see $z_i \, z_j=0$ 
for all $i\in I_\lambda$ and $j\in I_\mu$. 
We claim 
$1=\bigvee_{i\in I}z_i
=\sum_{\lambda\in\Lambda}\bigvee_{i\in I_\lambda}z_i$. 
This follows because 
$A\rtimes_\alpha G \cong \bigoplus_{i\in I}\mathbb{K}(H_i)$ 
and the sum of all diagonal atoms is equal to $1$ 
in strong operator topology. 
By definition of a central support projection 
we obtain 
$Mz_i \otimes \mathbb{B}(L^2(G))
=
M \otimes \mathbb{B}(L^2(G))\, e_i \,M \otimes \mathbb{B}(L^2(G))
$
for all $i\in I$. 
Then we see 
$\tilde{\alpha}(Mz_i \otimes \mathbb{B}(L^2(G)))
\subset
(M \otimes \mathbb{B}(L^2(G))e_iM \otimes \mathbb{B}(L^2(G)))
\otimes L^{\infty}(G)
=
Mz_i \otimes \mathbb{B}(L^2(G))
\otimes L^{\infty}(G)
$ 
and hence 
$\tilde{\alpha}(z_i\otimes1)\leq z_i\otimes 1\otimes1$ for all $i\in I$. 
Let $z_\lambda$ be the projection $\bigvee_{i\in I_\lambda}z_i$. 
Then a family $\{z_\lambda\}_{\lambda\in \Lambda}$ is 
a partition of unity and satisfies 
$\tilde{\alpha}(z_\lambda\otimes1)\leq z_\lambda\otimes1\otimes1$ for all 
$\lambda\in \Lambda$. 
Therefore we obtain 
$\tilde{\alpha}(z_\lambda\otimes1)=z_\lambda\otimes1\otimes1$ for all 
$\lambda\in \Lambda$. 
Since the left hand side is equal to 
$W_{23}\alpha(z_\lambda)_{13}W_{23}^*$, 
we have $\alpha(z_\lambda)=z_\lambda\otimes1$ for all 
$\lambda\in \Lambda$. 
Ergodicity of $\alpha$ yields $z_\lambda=1$ or $0$ for all 
$\lambda\in \Lambda$. 
Hence $\Lambda$ is a singleton and it completes the proof.
\end{proof}

From this lemma we have the following irreducibility criterion 
which has been studied in \cite[Section 10]{Wassermann1}. 

\begin{lem} \label{irreducibility}
Let $\pi$ be an element of $\widehat{G}$. 
If for any $\rho\in\wdh{G}$ there exists $k\in\Z_{\geq0}$ such that 
$\rho$ is contained in $\pi^k$, 
then for any $i,j\in I$ there exists $\ell\in\Z_{\geq0}$ with 
$(\mathbb{M}(\pi)^\ell)_{i,j}>0$.
\end{lem}

\begin{proof}
Since $G$-space $e_i (A\otimes\mathbb{K}(L^2(G)))e_j$ is not $0$, 
there exists an element $\rho$ in $I$ such that 
$\mathbb{M}(\rho)_{i,j}=\dim 
\Hom(H_\rho,e_i (A\otimes\mathbb{K}(L^2(G)))e_j)
>0$. 
By assumption there exists $\ell\in\N$ such that $\pi^\ell$ contains $\rho$. 
Then we have 
$(\mathbb{M}(\pi)^\ell)_{i,j}=\mathbb{M}(\pi^\ell)_{i,j}\geq 
\mathbb{M}(\rho)_{i,j}>0$. 
\end{proof}

In that case we want to show the existence of a Perron-Frobenius eigenvector, 
however, we cannot use the dual weight $\hat{\varphi}$ as 
in the compact group case 
because it is not necessarily a trace. 
We will give an eigenvector under some assumptions later. 

\begin{lem}\label{m(pi)}
Let $m:R(G)\longrightarrow \Z$ be a $*$-homomorphism such that 
it satisfies $m(\pi_0)=1$ and $m(\widehat{G})\subset \Z_{>0}$. 
Assume $\mathbb{M}(\pi)_{i,j}\leq m(\pi)$ for all $\pi\in \widehat{G}$ 
and $i,j \in I$. 
Then the matrix $\mathbb{M}(\rho)$ is a bounded operator on $l^2(I)$ and 
its norm is less than or equal to $m(\rho)^2$ for all $\rho\in R(G)_+$. 
\end{lem}

\begin{proof}
Take $\rho$ from $R(G)_+$ and 
we have 
$\sum_{j\in I}\mathbb{M}(\rho)_{i,j}^2 \leq m(\rho)^2$. 
In fact we can directly check it by the following calculation
\begin{align*}
\sum_{j\in I}\mathbb{M}(\rho)_{i,j}^2
=&(\mathbb{M}(\rho)\mathbb{M}(\ovl{\rho}))_{i,i}\\
=&\mathbb{M}(\rho\cdot\ovl{\rho})_{i,i}\\
=&\sum_{\sigma\in \widehat{G}}
N_{\sigma}^{\rho \ovl{\rho}} \mathbb{M}(\sigma)_{i,i}\\
\leq&\sum_{\sigma\in \widehat{G}}
N_{\sigma}^{\rho \bar{\rho}}m(\sigma)\\
=&m(\rho)m(\bar{\rho})\\
=&m(\rho)^2\, .
\end{align*}
Take a finitely supported vector $\xi$  in $l^2(\Z)$. 
Let $C(i,\rho)$ be a finite set 
$\{j\in I\mid \mathbb{M}(\rho)_{i,j}\neq0\}$. 
From the above inequality we particularly obtain 
$|C(i,\rho)|\leq m(\rho)^2$. 
We assume $\rho=\ovl{\rho}\in R(G)_+$. 
Then we have the desired inequality
\begin{align*}
||\mathbb{M}(\rho)\xi||^2
=&\sum_{i\in I}
\Big{|}\sum_{j\in I}\mathbb{M}(\rho)_{i,j}\xi_j\Big{|}^2\\
=
&\,
\sum_{i\in I}
\Big{|}\sum_{j\in C(i,\rho)}
\mathbb{M}(\rho)_{i,j}\xi_j\Big{|}^2\\
=
&\,
\sum_{i\in I}\sum_{j\in C(i,\rho)}
\mathbb{M}(\rho)_{i,j}^2
\cdot
\sum_{j\in C(i,\rho)}|\xi_j|^2\\
\leq
&\,
m(\rho)^2 \sum_{i\in I}\sum_{j\in C(i,\rho)}
|\xi_j|^2\\
=
&\,
m(\rho)^2 \sum_{j\in I}
|C(j,\rho)|
|\xi_j|^2\\
\leq&\,
m(\rho)^4 \sum_{j\in I}
|\xi_j|^2\, .\\ 
\end{align*}
For general $\rho\in R(G)_+$ we apply the above result to 
a positive self-conjugate 
$\ovl{\rho}\rho$ and get 
\begin{align*}
||\mathbb{M}(\rho)||^2
=
&\,
||\mathbb{M}(\ovl{\rho}\rho)||\\
\leq
&\, 
m(\ovl{\rho}\rho)^2\\
=
&\,
m(\rho)^4\, .
\end{align*}

\end{proof}

\begin{lem}\label{eigenvector}
If a self-conjugate element $\pi\in \widehat{G}$ 
has the property in Lemma \ref{irreducibility}, 
then $\mathbb{M}(\pi)$ has an eigenvector with the eigenvalue 
$||\mathbb{M}(\pi)||(\leq D_\pi^2)$. 
Moreover if $A$ has a tracial state, then 
$||\mathbb{M}(\pi)||(\leq d_\pi^2)$. 
\end{lem}

\begin{proof}
Since we know the matrix $\mathbb{M}(\pi)$ is irreducible and bounded, 
there exists an eigenvector with the eigenvalue $||\mathbb{M}(\pi)||$. 
For all $i,j\in I$ we have proved 
$e_i (A\otimes\mathbb{K}(L^2(G)))e_j \neq 0$. 
By Lemma \ref{pAq multiplicity} we have $\mathbb{M}(\pi)_{i,j}\leq D_\pi$. 
Hence we obtain the first assertion by Lemma \ref{m(pi)}. 
If $A$ has a tracial state $\tau$, then the tracial weight 
$\tilde{\tau}=\tau\otimes \Tr$ is semifinite on $A\rtimes_\alpha G$ and 
$\tilde{\tau}(e_i)>0$ for all $i\in I$. 
In fact if $\tilde{\tau}(e_i)$ is equal to 0 for some $i\in I$, 
we have $\tilde{\tau}(e_j)=0$ for all $j\in I$ 
by Remark \ref{tau(p)>0} 
and Lemma \ref{not 0}. 
It shows $\tilde{\tau}=0$ and this is contradiction. 
Now we apply Lemma \ref{m(pi)} with $m(\pi)=d_\pi$ 
and obtain the second assertion. 
\end{proof}

Now we consider the special case $G=SU_q(2)$. 
We refer to \cite{Woronowicz1,Masuda1} 
or the next section 
for its representation theory of $SU_q(2)$. 
Its irreducible representations are represented by 
$\{\pi_\nu\}_{\nu\in \frac{1}{2}\Z_{>0}}$ 
and its fusion rules are determined by
\[
\pi_{\frac{1}{2}}\cdot\pi_{\nu}=\pi_{\nu+\frac{1}{2}}+\pi_{\nu-\frac{1}{2}}
\]
for all $\nu\in  \frac{1}{2}$. 

\begin{lem}
Let $m:R(SU_q(2))\longrightarrow \Z$ be a dimension function with 
$m(\pi_\nu)\leq d_{\pi_\nu}^2$ for all $\nu\in \frac{1}{2}$. 
Then we obtain $m(\pi_\nu)=d_{\pi_\nu}$ for all $\nu\in \frac{1}{2}$.
\end{lem}

\begin{proof}
Write $t_0$ for $m(\pi_{\frac{1}{2}})$. 
Define polynomials $\{p_{\nu}\}_{\nu\in \frac{1}{2}\Z}$ 
recursively by $p_{\nu+\frac{1}{2}}(s)=s \,p_{\nu}(s)-p_{\nu-\frac{1}{2}}(s)$ 
and $p_0=1, \, p_{\frac{1}{2}}=s$. 
Then we get $m(\pi_\nu)=p_\nu(t_0)$. 
From the positivity of $m(\pi_\nu)$ we have $t_0\geq2$. 
Assume $t_0>2$. 
For $s>2$, the polynomials $p_\nu(s)$ are strictly increasing 
and we get a lower bound by an affine line $p_{\nu}(2)'(s-2)+p_{\nu}(2)$ 
where $p_{\nu}(s)'$ means the derivative at $s$. 
Putting $s=t_0$, we obtain an equality 
$m(\pi_\nu)\geq (t-2) p_{\nu}(2)'+2\nu+1$. 
The left hand side has upper bound by $d_{\pi_\nu}^2=(2\nu+1)^2$. 
In order to derive contradiction, it suffices to show the right hand side 
has a polynomial degree 3 with respect to $\nu$. 
By the definition of $p_\nu$, we get 
$p_{\nu+\frac{1}{2}}(2)'-p_\nu(2)'=p_\nu(2)'-p_{\nu-\frac{1}{2}}(2)'+2\nu+1$. 
This immediately gives $p_\nu(2)'=\frac{2}{3}\nu(\nu+1)(2\nu+1)$. 
\end{proof}

We show the following main result in this section 
which asserts the existence of 
multiplicity vector $\boldsymbol{c}$ under a tracial condition 
on $A$. 
\begin{thm}\label{multiplicity vector}
Let $\{A,SU_q(2),\alpha\}$ be a compact quantum ergodic system. 
Assume that $A$ has a tracial state. 
Then there exists a positive entry 
vector $\boldsymbol{c}=(c_i)_{i\in I}$ such that 
$\mathbb{M}(\pi_\nu)\,\boldsymbol{c}=d_{\pi_\nu}\,\boldsymbol{c}$ 
for all $\nu\in \frac{1}{2}$. 
\end{thm}

\begin{proof}
By Lemma \ref{eigenvector} 
there exists an eigenvector $\boldsymbol{c}$ 
such that 
$\mathbb{M}(\pi_{\frac{1}{2}})\,\boldsymbol{c}=t\, \boldsymbol{c}$ 
where $t=\big{|}\big{|}\mathbb{M}(\pi_{\frac{1}{2}})\big{|}\big{|}$. 
Since any $\pi_\nu$ is written by the polynomial $P_\nu$ of 
$\pi_\frac{1}{2}$, 
we can define the dimension function 
$m:R(SU_q(2))\longrightarrow \Z$ with 
$\mathbb{M}(\pi_\nu)\,\boldsymbol{c}=m(\pi_\nu)\,\boldsymbol{c}$ for all 
$\nu\in \frac{1}{2}$. 
We show $m(\pi_\nu)=||\mathbb{M}(\pi_\nu)||$ for all 
$\nu\in \frac{1}{2}$. 
The self-adjoint operator $\mathbb{M}(\pi_\nu)$ is written by 
the polynomial $P_\nu$ 
of $\mathbb{M}(\pi_\frac{1}{2})$. 
Hence we have 
$||\mathbb{M}(\pi_\nu)||\geq 
P_\nu(\big{|}\big{|}\mathbb{M}(\pi_{\frac{1}{2}})\big{|}\big{|})=
P_\nu(t)=m(\pi_\nu)$. 
The converse inequality is obtained by applying 
the Schur test to 
$\mathbb{M}(\pi_\nu)\,\boldsymbol{c}=m(\pi_\nu)\,\boldsymbol{c}$. 
Therefore we have $m(\pi_\nu)=||\mathbb{M}(\pi_\nu)||\leq d_{\pi_\nu}^2$. 
By previous lemma, we have $m(\pi_\nu)=d_{\pi_\nu}$ for all 
$\nu\in \frac{1}{2}$.
\end{proof}

The next lemma 
has already been proved in \cite[Lemma 1, Theorem 2]{Wassermann3} 
for $q=1$, that is, the $SU(2)$ case. 

\begin{lem}
Let $\{A,SU(2),\alpha\}$ be an ergodic system. 
Then its $\pi_{\frac{1}{2}}$-eigenvector space 
$X_{\pi_{\frac{1}{2}}}$ is not one-dimensional. 
\end{lem}

\begin{proof}
Let $\xi=(a,b)$ be a $\pi_{\frac{1}{2}}$-eigenvector. Then we have
\[
\alpha(a)=a\otimes x + b\otimes v\, , \quad
\alpha(b)=a\otimes u + b\otimes y.
\]
Consider a vector $\eta=(b^*,-q^{-1}a^*)$. 
By an easy calculation, 
we see $\eta$ is also a $\pi_{\frac{1}{2}}$-eigenvector. 
Assume that $X_{\pi_{\frac{1}{2}}}$ is one-dimensional. 
Then there exists a complex number $\mu$ such that 
$\eta=\mu\xi$. 
So we get $b^*=\mu a$ and $-q^{-1}a^*=\mu b$. 
It follows $-q^{-1}a=|\mu|^2 a$. 
Hence we obtain $0>-q^{-1}=|\mu|^2>0$. 
This is a contradiction. 
\end{proof}

Hence we obtain the result corresponding to \cite[Theorem 1]{Wassermann3} 
for positive $q$. 

\begin{cor}
Let $\{A,SU(2),\alpha\}$ be an ergodic system. 
Assume that $A$ has a tracial state and $q$ is positive. 
Then its multiplicity diagram is one of type 
$1,\T_n\,(n\geq2),\T,SU(2),D_n^*\,(n\geq 2),D_\infty^*,A_4^*,S_4^*$ 
and $A_5^*$. 
\end{cor}

\begin{ex}
We consider the multiplicity diagrams of quantum spheres $\CSq$. 
If $\lambda_0=c(n)$, then its spectral pattern 
(or finite dimensionality) allows only 
the diagram of type $SU(2)$ in Figure \ref{SU(2)} 
(see Appendix). 
If $0\leq\lambda_0\leq 1$, then its spectral pattern derives the 
diagram of type $\T$ in Figure \ref{Torus}. 
If $\{A,\suq,\alpha\}$ is an ergodic system and $A$ is finite dimensional, 
then its multiplicity diagram of $\pi_{\frac{1}{2}}$ must be of type $SU(2)$ 
by its finite dimensionality. 
Hence the classification by Podle\'{s} (\cite{Podles1}) 
shows that $A$ is $G$-isomorphic to $\End(H_{\pi_{\nu}})$ for some 
$\nu \in \Z_{\geq \frac{1}{2}}$. 
\end{ex}

We recall the definition of the McKay diagrams. 
Let $H\subset G$ be a quantum subgroup with the restriction map $r_H$ and 
$w\in \mathbb{B}(H_w)\otimes C(G)$ be a unitary representation 
of $G$. 
We denote the restricted representation $(\id\otimes r_H)(w)$ by 
$w|_H$. 
Prepare the vertices $\{\sigma\}_{\sigma\in \wdh{H}}$. 
Let $\sigma_0$ be the one-dimensional trivial representation of $H$. 
First we consider the irreducible decomposition 
$w|_H\cdot \sigma_0
=\oplus_{\sigma\in \wdh{H}}N_\sigma^{w|_H}\sigma$ where 
the scalar $N_\sigma^{w|_H}$ means a multiplicity. 
Then we draw arrows from $\sigma_0$ to the above arising 
irreducible representations $\{\sigma\}_{\sigma\in \wdh{H}}$ 
with $N_\sigma^{w|_H}$-times, respectively. 
Second for each $\sigma$ in the above decomposition we 
consider the irreducible decomposition of $w|_H\cdot \sigma$ 
and we draw arrows from $\sigma$ to the irreducible representations 
in a similar way. 
Continue this procedure we get the McKay diagram for $H\subset G$ 
with respect to $w$. 
If $w$ is self-conjugate, then the diagram is non-oriented 
because of the symmetry $N_\sigma^{w|_H\,\tau}=N_{w|_H\,\sigma}^\tau$. 
When we treat quantum subgroups of $\suq$, the McKay diagrams 
are drawn with respect to the fundamental representation 
$\pi_{\frac{1}{2}}$. 
Now we consider $G$-covariant system $\{C(H\setminus G),\delta\}$ 
where $H$ is a quantum subgroup with the restriction map $r_H$. 
It is well-known that the multiplicity diagram and the McKay diagram 
coincide in the classical case, however, we give the proof of the 
general quantum group case for readers' convenience. 
We denote $(r_H\otimes\id)(V_{\ell})$ by $V_{\ell}|_H$. 
The $C^*$-algebra $C(G)\otimes \mathbb{K}(L^2(G))$ 
has the left and right actions $\alpha_{\ell}$ and $\beta_r$ defined by 
$\alpha_{\ell}(x)=\Ad {V_{\ell}}_{13}^*{V_{\ell}}_{12}^*(1\otimes x)$ 
and $\beta_r(x)=\Ad V_{13}(x\otimes 1)$ for all 
$x\in C(G)\otimes \mathbb{K}(L^2(G))$, respectively. 
We define the left $H$-action $\alpha_{\ell}^H$ by 
the composition $(r_H\otimes \id)\circ \alpha_{\ell}$. 
Consider the map 
$\Ad V_{\ell}: C(G)\otimes \mathbb{K}(L^2(G)) \longrightarrow 
C(G)\otimes \mathbb{K}(L^2(G))$. 

\begin{lem}
\begin{enumerate}

\item 

The map $\Ad V_{\ell}$ gives an isomorphism between 
$C(H\setminus G)\otimes \mathbb{K}(L^2(G))$ and 
${}^H (C(G)\otimes \mathbb{K}(L^2(G)))$ where 
the latter one is the fixed point algebra 
of the left $H$-action $\alpha_{\ell}^H$. 
Moreover it maps $C(H\setminus G)\rtimes_\delta G$ to 
$\C\otimes {}^H \mathbb{K}(L^2(G)))$. 

\item 

$\Ad V_{\ell}$ intertwines the right $G$-actions 
$\wdt{\alpha}$ and $\beta_r$. 

\end{enumerate}
\end{lem}

\begin{proof}
(1) It suffices to show that the map 
$\Ad V_{\ell}$ intertwines left $H$-actions 
$\Ad V_{\ell}|_H \otimes \id$ and $\alpha_{\ell}^H$. 
This is immediately verified by the equality 
${V_{\ell}}_{\,23}{V_{\ell}|_H}_{12}^*={V_{\ell}|_H}_{13}^*{V_{\ell}|_H}_{12}^*{V_{\ell}}_{\,23}$. 
(2) Similarly the equality 
${V_{\ell}}_{\,12}W_{23}V_{13}=V_{13}{V_{\ell}}_{\,12}$ gives the desired 
intertwining property. 
\end{proof}

Now we choose a left irreducible $H$-module $K_\sigma^\ell$ for each 
$\sigma\in\wdh{H}$. 
Then we have the irreducible decomposition 
$L^2(G)=\oplus_{\sigma\in \wdh{H}}K_\sigma^\ell \otimes L_\sigma$ 
where $L_\sigma$ is the multiplicity space for $\sigma$. 
Note that all the $L_\sigma$ are non-zero, 
because there exists $\pi\in \wdh{G}$ such that $\sigma$ is 
contained in $\pi|_H$. 
Then we get 
${}^H \mathbb{K}(L^2(G)))
=\oplus_{\sigma\in \wdh{H}}\C1_\sigma\otimes \mathbb{K}(L_\sigma)$. 
Let $e_\sigma=1_\sigma\otimes p_\sigma$ be a minimal projection 
of $\C1_\sigma\otimes \mathbb{K}(L_\sigma)$. 
By the previous lemma we have 
$K_0(C(H\setminus G)\rtimes_\delta G)
=\oplus_{\sigma\in \wdh{H}}\Z [e_\sigma]$. 
Hence the vertices of multiplicity maps are represented by 
the elements of $\wdh{H}$. 
Then the following proposition holds. 

\begin{prop}
Let $\pi$ be an element of $\wdh{G}$. 
Then we have $\mathbb{M}(\pi)_{\rho,\sigma}=N_\sigma^{\pi|_H\,\rho}$ 
for all $\rho$, $\sigma$ in $\wdh{H}$. 
\end{prop}

\begin{proof}
We use Corollary \ref{reduced space characterization}. 
The $G$-module $e_\rho C(H\setminus G)\otimes \mathbb{K}(L^2(G))e_\sigma$ 
is isomorphic to ${}^H (C(G)\otimes e_\rho \mathbb{K}(L^2(G))e_\sigma)$ 
by the previous lemma. 
Since quantum group $G$ acts on the first tensor component, 
the multiplicity space for $\pi$ is 
${}^H (H_\pi^{\ell}\otimes e_\rho \mathbb{K}(L^2(G))e_\sigma)$ which 
is linearly isomorphic to 
$\Hom_H(H_{\ovl{\pi|_H}}^{\ell}, K_\rho \otimes K_{\ovl{\sigma}})$. 
Hence its dimension is equal to 
$N_{\ovl{\pi|_H}}^{\rho\,\ovl{\sigma}}=N_\sigma^{\pi|_H\,\rho}$. 
\end{proof}

Define the vector $\boldsymbol{d}=(d_{\sigma})_{\sigma\in \wdh{H}}$ and 
we have $\mathbb{M}(\pi)\boldsymbol{d}=d_\pi \boldsymbol{d}$. 
Hence we have solved the eigenvector problem in the case of 
the quotient spaces. 

\begin{cor}
Let $H$ be a quantum subgroup of $G$ and consider the 
ergodic system $\{C(H\setminus G),\delta\}$. 
Then the multiplicity diagram for $\pi\in \wdh{G}$ coincides 
with the McKay diagram of $H\subset G$ with respect to $\pi$. 
\end{cor}

Next we study reduced ergodic systems associated to ergodic systems 
and show the heredity of information of 
multiplicity maps. 

\begin{defn}
Let $\{A,G,\alpha\}$ and $\{B,G,\alpha\}$ be 
compact quantum ergodic systems and 
take index sets $I,J$ with 
$A\rtimes_\alpha G \cong \bigoplus_{i\in I}\mathbb{K}(H_i)$ and 
$B\rtimes_\beta G \cong \bigoplus_{j\in J}\mathbb{K}({H_i}')$. 
We say two ergodic systems  are
of the same type if 
there exists a bijective map 
$\theta:I\longrightarrow J$ 
satisfying 
$\mathbb{M}^B(\pi)\circ \theta =\theta \circ \mathbb{M}^A(\pi)$ 
for all $\pi\in G$, 
where
$\mathbb{M}^A$ and $\mathbb{M}^B$ are multiplicity maps for ergodic systems 
$\{A,G,\alpha\}$ and $\{B,G,\alpha\}$ respectively. 
\end{defn}

\begin{thm}\label{reduced ergodic system}
Let $\{A,G,\alpha\}$ be a compact quantum ergodic system and 
take an index set $I$ with 
$A\rtimes_\alpha G\cong \bigoplus_{i\in I} \mathbb{K}(H_i)$. 
Take minimal projections $\{e_i\}_{i\in I}$ 
from $\{\mathbb{K}(H_i)\}_{i\in I}$. 
Then for all $i\in I$, two ergodic systems 
$\{A,G,\alpha\}$ and 
$\{e_i(A\otimes\mathbb{K}(L^2(G)))e_i,G,\tilde{\alpha}\}$ 
are of the same type. 
\end{thm}

\begin{proof}
Denote $\mathbb{K}(L^2(G))$ and $\mathbb{B}(L^2(G))$ 
simply by $\mathbb{K}$ and $\mathbb{B}$. 
We claim that 
the covariant system 
$\{A\otimes\mathbb{K}\otimes\mathbb{K},G,\tilde{\tilde{\alpha}}\}$ is 
strongly equivalent to 
$\{\mathbb{K}\otimes A\otimes\mathbb{K},G,\id\otimes\tilde{\alpha}\}$. 
Let us consider the map 
$\Ad(W_{23}^*):
A\otimes\mathbb{K}\otimes\mathbb{K}
\longrightarrow A\otimes\mathbb{K}\otimes\mathbb{K}$. 
For $x\in A, k\in \mathbb{K}\otimes\mathbb{K}$ we have: 
\begin{align*}
W_{23}^* \tilde{\tilde{\alpha}}(W_{23}(x\otimes k) W_{23}^*)W_{23}
=&
W_{23}^* W_{34}W_{24}W_{23} \alpha(x)_{14} (1\otimes k\otimes 1) 
W_{23}^* W_{24}^* W_{34}^* W_{23}\\
=&
W_{34} \alpha(x)_{14}(1\otimes k\otimes 1) W_{34}^*, 
\end{align*}
where we use the pentagonal identity for $W$, 
$W_{23}W_{13}W_{12}=W_{12}W_{23}$. 
Next flip the first and second tensor component of 
$A\otimes\mathbb{K}\otimes\mathbb{K}$ and the claim is proved. 
Hence we get an isomorphism 
$(A\otimes\mathbb{K})\rtimes_{\tilde{\alpha}} G 
\cong \mathbb{K}\otimes (A\rtimes_\alpha G)$. 
A projection $e_i\otimes 1$ is in the fixed point algebra of 
the multiplier algebra $M(A\otimes\mathbb{K}\otimes\mathbb{K})$ for 
$\tilde{\tilde{\alpha}}$. 
So it is mapped to a projection $f_i$ in the multiplier algebra 
$M(\mathbb{K}\otimes A\rtimes_\alpha G)$. 
Hence we see that 
$\{e_i(A\otimes\mathbb{K})e_i \otimes\mathbb{K},G,\tilde{\tilde{\alpha}}\}$ 
is equivalent to 
$\{f_i(\mathbb{K}\otimes A\otimes\mathbb{K})f_i,G,\id\otimes\tilde{\alpha}\}$. 
The projection $f_i$ is decomposed into a direct sum of 
projections $\{p_j\}_{j\in I}$ which are in 
$\{\mathbb{B}\otimes \mathbb{B}(H_j)\}_{j\in I}$. 
We claim they are non-zero projections. 
It suffices to 
prove the central support of $f_i$ in 
$\mathbb{B}\otimes M\rtimes_{\alpha}G$ is 
equal to $1$ 
by passing to von Neumann algebras, that is, 
isomorphic covariant systems 
$\{e_i(M\otimes\mathbb{B})e_i \otimes\mathbb{B},G,\tilde{\tilde{\alpha}}\}$
and 
$\{f_i(\mathbb{B}\otimes M\otimes\mathbb{B})f_i,G,\id\otimes\tilde{\alpha}\}$. 
Hence we show the central support of $e_i \otimes 1$ in 
$(M\otimes \mathbb{B}\otimes \mathbb{B})^{\tilde{\tilde{\alpha}}}
=(M\otimes\mathbb{B})\rtimes_{\tilde{\alpha}} G$
is equal to 1. 
By duality theorem, we get the isomorphism of inclusions,
$M\rtimes_{\alpha} G \otimes \C 
\subset (M\otimes\mathbb{B})\rtimes_{\tilde{\alpha}} G
\cong
\hat{\alpha}(M\rtimes_{\alpha} G)
\subset 
M\rtimes_{\alpha} G \otimes \mathbb{B}$. 
Hence we study the central support of 
$\hat{\alpha}(e_i)$ in 
$M\rtimes_{\alpha} G \otimes \mathbb{B}$. 
Fix $j\in I$ and we can take $\pi\in \widehat{G}$ with 
$\mathbb{M}(\pi)_{j,i}>0$ 
by Lemma \ref{not 0}. 
By the definition of $\mathbb{M}(\pi)$ 
we have 
$[\hat{\alpha}_{\pi}(e_i)]=\sum_{k\in I}\mathbb{M}(\pi)_{k,i}[e_k]$. 
Hence the $j$-th component of $\hat{\alpha}_{\pi}(e_i)$ is a 
non-zero projection, in particular the central support of 
$\hat{\alpha}(e_i)$ in $M\rtimes_{\alpha} G \otimes \mathbb{B}$ 
is equal to $1$. 
Let ${H_j}'$ be a Hilbert space $p_j( L^2(G)\otimes H_j)$ and 
we have an isomorphism 
$e_i(A\otimes\mathbb{K})e_i \rtimes_{\tilde{\alpha}} G 
\cong 
\bigoplus_{j\in I}\mathbb{K}({H_j}')$. 
Next we choose minimal projections 
$\{q_j\}_{j\in I}$ from 
$\mathbb{K}({H_j}')$. 
Let $p$ be a minimal projection in $\mathbb{K}$. 
For any $j\in I$, 
the projection $q_j$ is equivalent to 
$p\otimes e_j$ in 
$\mathbb{K}\otimes\mathbb{K}(H_j)\subset 
(\mathbb{K}\otimes A\otimes\mathbb{K})^{\id\otimes\tilde{\alpha}}$. 
Then we have a $G$-isomorphism between reduced spaces 
$q_j(\mathbb{K}\otimes A\otimes\mathbb{K})q_k$ 
and 
$(p\otimes e_j)(\mathbb{K}\otimes A\otimes\mathbb{K})(p\otimes e_k)
=
\C p \otimes e_j(A\otimes \mathbb{K})e_k$ 
for all $j,k\in I$. 
Therefore we have proved the equality 
$\dim\Hom_G(H_\pi,q_j(\mathbb{K}\otimes A\otimes\mathbb{K})q_k)
=\dim\Hom_G(H_\pi,e_j(A\otimes \mathbb{K})e_k)$ 
for all $\pi\in\widehat{G}$ and $j,k\in I$. 
\end{proof}

\begin{defn}
Let $\{A,G,\alpha\}$ be a compact quantum ergodic system. 
A unital $C^*$-subalgebra $B\subset A$ is called 
a $G$-invariant $C^*$subalgebra 
if 
it satisfies 
$\alpha(B)\subset B\otimes C(G)$. 
Then $\{B,G,\alpha\}$ is called a subsystem of 
$\{A,G,\alpha\}$ and 
denote this situation by 
$\{B,G,\alpha\}\subset\{A,G,\alpha\}$. 
\end{defn}

Let $B$ be a $G$-invariant $C^*$-subalgebra of $A$. 
The inclusion $G$-homomorphism $\iota:B\longrightarrow A$ induces 
the inclusion of crossed products, 
$B\rtimes_\alpha G \subset A\rtimes_\alpha G$. 
Taking index sets $I,J$ as 
$B\rtimes_\alpha G \cong \bigoplus_{i\in I}\mathbb{K}(H_i)$ 
and
$A\rtimes_\alpha G \cong \bigoplus_{j\in J}\mathbb{K}(K_i)$. 
Let $\Lambda$ be an inclusion matrix of 
$B\rtimes_\alpha G \subset A\rtimes_\alpha G$, 
that is, 
$\mathbb{K}(H_i)$ is amplified into 
$\mathbb{K}(K_j)$ by $\Lambda_{j,i}$ times. 
We can easily show the next proposition by the definition of 
multiplicity maps. 

\begin{prop}\label{intertwining property}
Let 
$\{B,G,\alpha\}\subset\{A,G,\alpha\}$ 
be an inclusion of compact quantum ergodic systems. 
Let 
$\mathbb{M}^B$ and $\mathbb{M}^A$ be the multiplicity maps. 
Then we have 
$\Lambda\,\mathbb{M}^B(\pi) =\mathbb{M}^A(\pi)\, \Lambda$ 
for all $\pi\in \widehat{G}$. 
\end{prop}

We recall an observation by A. Wassermann, which has 
an important role in his classification program in 
\cite{Wassermann3}. 
Let $p_0=(\id\otimes h)(V)$ be a minimal projection of $C_r^*(G)$. 
Note that it is central in not only $C_r^*(G)$ but also $C_{\ell}^*(G)$. 
It is also a minimal projection on $L^2(G)$. 
Then we have 
\[(1\otimes p_0)A\rtimes_\alpha G (1\otimes p_0)
=\big{(}A\otimes p_0\mathbb{K}(L^2(G))p_0\big{)}^{\tilde{\alpha}}
=A^\alpha\otimes \C p_0
=\C(1\otimes p_0).
\] 
This shows that $p_0$ is a minimal projection in $A\rtimes_\alpha G$. 
Hence if we take minimal projections in $\{e_i\}_{i\in I}$ as before, 
then there exists unique $i_0$ such that 
$e_{i_0}$ is equivalent to $p_0$ in $A\rtimes_\alpha G$. 
We often say this index corresponds to $p_0$. 
Then we obtain an isomorphism as covariant systems 
between $\{A,G,\alpha\}$ and 
$\{e_{i_0}A\otimes\mathbb{K}(L^2(G))e_{i_0},G,\tilde{\alpha}\}$. 
From this we get the following result. 

\begin{cor}\label{multiplicity vector for right coideal}
Let $A\subset C(G)$ be a right coideal. 
Then there exists an eigenvector $\boldsymbol{c}=(c_i)_{i\in I}$ 
which satisfies the following conditions. 
\begin{enumerate}
\item All entries of $\boldsymbol{c}$ are strictly positive integers and 
there exists an index $i_0\in I$ which satisfies $c_{i_0}=1$. 

\item It is a common eigenvector of the multiplicity map; 
$\mathbb{M}(\pi)\boldsymbol{c}=d_\pi \boldsymbol{c}$ 
for all $\pi\in\widehat{G}$. 

\item Two covariant systems $\{A,G,\alpha\}$ and 
$\{e_{i_0}A\otimes\mathbb{K}(L^2(G))e_{i_0},G,\tilde{\alpha}\}$ are 
isomorphic. 

\end{enumerate}
\end{cor}

\begin{proof}
Take a minimal 
projection $p_\pi$ of 
$\End(H_\pi)\subset C_r^*(G)$ for $\pi\in\widehat{G}\, .$ 
Let $\pi$ be an element of $\widehat{G}$. 
Consider the dual coaction 
$\hat{\delta}_{\pi}:C_r^*(G)\longrightarrow C_r^*(G)\otimes \End(H_\pi)$. 
We want to compute $\hat{\delta}_{\pi}(p_0)$ in $K_0(C_r^*(G))$. 
Since we have the isomorphism 
$K_0(C_r^*(G)) \cong R(G)=K_0^G(\C)$ by 
$[p_\pi]\rightarrow \pi$, 
we see $[\hat{\delta}_{\pi}(p_0)]=[p_\pi]$. 
Now let $\Lambda$ be an inclusion matrix for 
$A\rtimes_\alpha G\subset C(G)\rtimes_\delta G$. 
It is actually a row vector because $C(G)\rtimes_\alpha G$ is isomorphic 
to $\mathbb{K}(L^2(G))$. 
By the previous proposition we have 
$\mathbb{M}'(\pi)\Lambda=\Lambda\mathbb{M}(\pi)$ 
for all $\pi \in \widehat{G}$, 
where $\mathbb{M}'$ and $\mathbb{M}$ are the multiplicity maps for 
$\{A,G,\alpha\}$ and $\{C(G),G,\delta\}$ respectively. 
Let us take $1\otimes p_0$ for a minimal projection of 
$C(G)\rtimes_\delta G \cong \mathbb{K}(L^2(G))$. 
From the first part of this proof 
we have $\pi[p_0]=[p_\pi]=d_\pi[p_0]$ in 
$K_0(C(G)\rtimes_\delta G)$. 
Hence the action $\mathbb{M}'(\pi)$ is multiplication by $d_\pi$ and 
we get $\Lambda\mathbb{M}(\pi)=d_\pi \Lambda$. 
Transposing it, we have 
$\mathbb{M}(\overline{\pi}){}^T\Lambda=d_{\overline{\pi}} {}^T\Lambda$. 
Define a vector $\boldsymbol{c}={}^T\Lambda$. 
Let us take an index $i_0$ which corresponds to $p_0$ in $A\rtimes_\alpha G$. 
Since $1\otimes p_0$ is also minimal in $C(G)\rtimes_\delta G$, 
we get $c_{i_0}=\Lambda_{i_0}=1$. 
\end{proof}

We end this section with the following proposition. 
We use the same notations as before. 

\begin{prop}\label{criterion on eigenvector}
Let $A\subset C(G)$ be a right coideal and take 
an eigenvector $\boldsymbol{c}={}^T \Lambda$ 
for its multiplicity map. 
If there exists an index $j\in I$ with $c_j=1$, 
then the reduced ergodic system $e_j(A\otimes\mathbb{K}(L^2(G)))e_j$ 
is $G$-equivariantly embedded into $C(G)$. 
In particular it also becomes a right coideal of $C(G)$. 
\end{prop}
\begin{proof}
Let $p_0$ be a minimal projection in $\mathbb{K}(L^2(G))$ 
associated to the trivial representation. 
Since the ranks of $e_j$ and $p_0$ in $C(G)\rtimes_\delta G \cong 
\mathbb{K}(L^2(G))$ are equal, we have a partial isometry 
$v$ in $C(G)\rtimes_\delta G$ such that it satisfies 
$v^*v=e_j$ and $vv^*=p_0$. 
Note that $v$ belongs to 
$C(G)\rtimes_\delta G=(C(G)\otimes \mathbb{K}(L^2(G)))^G$. 
Then we have a desired $G$-equivariant $*$-homomorphism 
$e_j(A\otimes\mathbb{K}(L^2(G)))e_j \longrightarrow C(G)\otimes \C p_0$ 
defined by $a\mapsto vav^*$. 
\end{proof}

\begin{rem}\label{faithfulness}
Let $\{A,G,\alpha\}$ and $\{B,G,\beta\}$ be two ergodic 
systems and $\theta:A\longrightarrow B$ be a 
$G$-equivariant $*$-homomorphism. 
Then it must be faithful as we see below. 
Let $\varphi_A$ and $\varphi_B$ be the unique faithful 
invariant states. 
Since $\theta$ is $G$-equivariant, $\varphi_B\circ \theta$ is 
an invariant state. 
By its uniqueness we obtain $\varphi_A=\varphi_B\circ\theta$. 
Hence $\theta$ is faithful. 
\end{rem}

\section{Elementary results for $SU_q(2)$}

In this section we summarize basic facts for $\suq$ 
for readers' convenience. 
Readers are referred to, for example, 
\cite{Woronowicz1,Masuda1} for 
its basic theory and 
\cite{Podles1,Podles2,Koornwinder1,Noumi1} for 
the results on the quantum spheres or quantum subgroups in 
$\suq$. 
We adopt the same notation as \cite{Masuda1} in this paper. 
We treat a real number $q$ in $[-1,1]\setminus \{0\}$. 
The smooth function algebra $\asuq$ is the universal $*$-algebra 
generated by four elements $x,u,v$ and $y$ with the following relations
\begin{center}
$ux=q xu,\quad vx=q xv,\quad yu=q uy,\quad yv=qvy$,
\end{center}
\begin{center}
$uv=vu,\quad xy-q^{-1}uv=yx-quv=1,\quad
x^*=y,\quad u^*=-q^{-1}v$.
\end{center}
Its universal $C^*$-algebra is denoted by $\csuq$. 
We often use a positive operator $\zeta=-q^{-1}uv$. 
We make $\asuq$ a Hopf $*$-algebra by defining 
coproduct $\delta$, counit $\varepsilon$ and antipode $\kappa$ as 
follows
\[
\begin{pmatrix}
\delta(x) & \delta(u) \\
\delta(v) & \delta(y)\\
\end{pmatrix}
=
\begin{pmatrix}
x\otimes1 & u\otimes1 \\
v\otimes1 & y\otimes1\\
\end{pmatrix}
\cdot
\begin{pmatrix}
1\otimes x &1\otimes u \\
1\otimes v &1\otimes y \\
\end{pmatrix},
\]
\[
\begin{pmatrix}
\varepsilon(x) &\varepsilon(u) \\
\varepsilon(v) &\varepsilon(y) \\
\end{pmatrix}
=
\begin{pmatrix}
 1&0 \\
 0&1 \\
\end{pmatrix},
\]
\[
\begin{pmatrix}
\kappa(x) &\kappa(u) \\
\kappa(v) &\kappa(y) \\
\end{pmatrix}
=
\begin{pmatrix}
 y& -q u\\
 -q^{-1}v & x\\
\end{pmatrix}.
\]
Then the pair $SU_q(2)=(\csuq,\delta)$ becomes a compact quantum group 
and it is often called the \textit{twisted} $SU(2)$ \textit{group}. 
The maps $\delta$ and $\varepsilon$ are extended to $\csuq$ norm continuously 
and the map $\kappa$ is extended to $\csuq$ as a closed operator. 
Its Woronowicz characters $\{f_z\}_{z\in\C}$ are given by
\[
\begin{pmatrix}
f_z(x) & f_z(u) \\
f_z(v) & f_z(y)\\
\end{pmatrix}
=
\begin{pmatrix}
|q|^z & 0\\
 0& |q|^{-z}\\
\end{pmatrix}
\ 
\mbox{\textrm{for all}}
\ 
z\in \C. 
\] 
The equivalence classes of irreducible representations 
$\widehat{\suq}$ are indexed by spin numbers 
$\nu\in \frac{1}{2}\Z_{\geq0}$ 
and 
we fix a selection of irreducible representations 
$\{w(\pi_\nu)\}_{\nu}$ corresponding to spins $\nu$ as follows. 
The representation space $H_\nu=H_{w(\pi_\nu)}$ is 
$2\nu+1$-dimensional and 
fix the orthonormal basis $\{\xi_{r}^{\nu}\}_{r\in I_{\nu}}$ 
where the index set $I_{\nu}$ is 
$\{-\nu,-\nu+1,\ldots,\nu-1,\nu\}$. 
Then we define the matrix $w(\pi_\nu)\in \mathbb{B}(H_\nu)\otimes \asuq$ 
by setting $w(\pi_\nu)_{i,j}$ as follows. 
\begin{enumerate}

\item Case $i+j\leq0,\, i\geq j$: 
\[x^{-i-j}v^{i-j}q^{(\nu+j)(j-i)}{\nu+i \brack i-j}_{q^2}^{\frac{1}{2}}
{\nu-j \brack i-j}_{q^2}^{\frac{1}{2}} P_{\nu+j}^{(i-j,-i-j)}(\zeta;q^2),\]

\item Case $i+j\leq0,\, i\leq j$: 
\[x^{-i-j}u^{j-i}q^{(\nu+i)(i-j)}{\nu-i \brack j-i}_{q^2}^{\frac{1}{2}}
{\nu+j \brack j-i}_{q^2}^{\frac{1}{2}} P_{\nu+i}^{(j-i,-i-j)}(\zeta;q^2),\]

\item Case $i+j\geq0,\, i\leq j$: 
\[q^{(j-i)(j-\nu)}{\nu-i \brack j-i}_{q^2}^{\frac{1}{2}}
{\nu+j \brack j-i}_{q^2}^{\frac{1}{2}} P_{\nu-j}^{(j-i,i+j)}(\zeta;q^2)
u^{j-i}y^{i+j},\]

\item Case $i+j\geq0,\, i\geq j$: 
\[q^{(i-j)(i-\nu)}{\nu+i \brack i-j}_{q^2}^{\frac{1}{2}}
{\nu-j \brack i-j}_{q^2}^{\frac{1}{2}} P_{\nu-i}^{(i-j,i+j)}(\zeta;q^2)
v^{i-j}y^{i+j},\]

\end{enumerate}

where we have used the $q$-\textit{binomial coefficients} and 
the \textit{little} $q$-\textit{Jacobi polynomials}:
\[{m \brack n}_q=\frac{(q;q)_m}{(q;q)_n(q;q)_{m-n}}, \quad 
(t;q)_m=\prod_{s=0}^{m-1}(1-tq^s),
\]
\[
P_{n}^{(\alpha,\beta)}(z;q)=
\sum_{r\geq0}
\frac{(q^{-n};q)_r(q^{\alpha+\beta+n+1};q)_r}
     {(q;q)_r(q^{\alpha+1};q)_r}
 (qz)^r.
\]
This yields the following formula
\[
w(\pi_\nu)_{i,j}^*=(-q)^{i-j}w(\pi_\nu)_{-i,-j}
\]
for all $\pi_\nu\in \widehat{\suq}$ and $i,j\in I_\nu$. 
For an integer $n$, 
we define the \textit{q-integer} and its factorial by
\[
(n)_q=\frac{|q|^n-|q|^{-n}}{|q|-|q|^{-1}}\,,
\quad
(n)_q!=(n)_q(n-1)_q\cdots(1)_q\,.
\]
We summarize useful formulae as follows.
\begin{enumerate}

\item $f_z(w(\pi_{\nu})_{r,s})=\delta_{r,s}|q|^{-2rz},$
\vspace{0.2cm}
\item $\sigma_{t}^h(w(\pi_{\nu})_{r,s})=|q|^{-2(r+s)it}w(\pi_{\nu})_{r,s},$
\vspace{0.2cm}
\item $\tau_{t}(w(\pi_{\nu})_{r,s})=|q|^{-2(r-s)it}w(\pi_{\nu})_{r,s},$
\vspace{0.2cm}
\item $R(w(\pi_{\nu})_{r,s})=(-1)^{r-s}w(\pi_{\nu})_{-s,-r}, $

\end{enumerate}
for all $z\in\C$, $t\in \R$, $\nu\in \frac{1}{2}\Z_\geq 0$ and 
$r,s\in I_\nu$. 
The $F$-matrix $F_{\pi_\nu}=(\id\otimes f_1)(w(\pi_\nu))$ is equal 
to $\diag(q^{2\nu},q^{2\nu-2},\ldots,q^{-(2\nu-2)},q^{-2\nu})$.  
Hence we obtain $D_{\pi_\nu}=(2\nu+1)_q$. 
In the spin $\nu=\frac{1}{2}$ case we have the following matrices
\[
\begin{pmatrix}
\sigma_{t}^h(x )& \sigma_{t}^h(u )\\
\sigma_{t}^h(v )& \sigma_{t}^h(y) \\
\end{pmatrix}
=
\begin{pmatrix}
|q|^{2it}x & u \\
v & |q|^{-2it}y \\
\end{pmatrix},
\]\vspace{0.2cm}
\[
\begin{pmatrix}
\tau_t(x) & \tau_t(u )\\
\tau_t(v) & \tau_t(y) \\
\end{pmatrix}
=
\begin{pmatrix}
x & |q|^{2it}u \\
|q|^{-2it}v & y \\
\end{pmatrix},
\]\vspace{0.2cm}
\[
\begin{pmatrix}
R(x) & R(u) \\
R(v) & R(y) \\
\end{pmatrix}
=
\begin{pmatrix}
y & -v \\
-u & x \\
\end{pmatrix}.
\]
For the calculation about the Haar state, we have
\[
h(\zeta^n)=\frac{1-q^2}{1-q^{2(n+1)}}
\  \mbox{for all} \  n\in \Z_{\geq0}.
\]

We recall the embedding of $\T$ into 
$\suq$ for $-1\leq q <1$. 
The torus group $\T$ is identified with 
the set $\{z\in \C\mid |z|=1\}$ and 
its cyclic subgroup $\T_m\, (m\geq2)$ is generated by 
$\exp(\frac{2\pi\sqrt{-1}}{m})$. 
For $-1\leq q<1$ 
if $\theta:C(\suq)\longrightarrow C(\T)$ is a restriction map, 
then it must be as follows. 
\begin{equation*}\label{tor01}
\begin{pmatrix}
\theta(x) & \theta(u)\\
\theta(v) & \theta(y)\\
\end{pmatrix}
=
\begin{pmatrix}
z & 0\\
0 & \ovl{z}\\
\end{pmatrix}
\ \mbox{or}\ 
\begin{pmatrix}
\ovl{z} & 0\\
0 & z\\
\end{pmatrix}\, ,
\end{equation*}
where $z\in \T\subset \C$ is a usual coordinating map. 
Let $\pi_\T$ and $\pi_\T'$ be the first and second 
restriction maps 
in the above equalities, respectively. 
$\T$ is called the \textit{maximal torus subgroup} of $\suq$ 
and we always treat it with the restriction map $\pi_\T$. 
Note that quotient space $\T\setminus \suq$ or 
$\T_m\setminus\suq$ do not depend on the choice of 
$\pi_\T$ or $\pi_\T'$. 
The quotient space $C(\T\setminus \suq)$ is often 
called the \textit{canonical homogeneous sphere}, 
which is generated by $\{w(\pi_1)_{0,r}\}_{r\in I_1}$. 
The quotient space $C(\T_{2}\setminus \suq)$ is 
denoted by $\csoq$ which is also a compact quantum group 
by restricting the coproduct $\delta$. 
We can easily see $\csoq$ has the spectral pattern: 
$\oplus_{k\in\Z_{\geq0}}(2k+1)\pi_k$ and hence 
it does not depend on positivity or negativity of the 
parameter $q$. 
In fact we have an isomorphism 
$\Xi_q:\csoq\longrightarrow\csomq$ as compact quantum groups 
defined by the following equality, 
where $2$ by $2$ matrices 
$
\begin{pmatrix}
x & u \\
v & y
\end{pmatrix}
$
and 
$
\begin{pmatrix}
a & b \\
c & d
\end{pmatrix}
$ 
are fundamental representations of $\suq$ and $\sumq$ 
respectively: 
\begin{align*}
&
\begin{pmatrix}
\Xi_q(x^2) & \sqrt{1+q^2} \Xi_q(xu) & \Xi_q(u^2)\\
\sqrt{1+q^2}\Xi_q(xv) & \Xi_q(1+(q+q^{-1})uv) & \sqrt{1+q^2}\Xi_q(uy) \\
\Xi_q(v^2) & \sqrt{1+q^2}\Xi_q(vy) & \Xi_q(y^2)
\end{pmatrix}
\\
&\quad=
\begin{pmatrix}
a^2 & -i\sqrt{1+q^2}ab & b^2\\
i\sqrt{1+q^2}ac & 1-(q+q^{-1})bc & i\sqrt{1+q^2}bd\\
c^2 & -i\sqrt{1+q^2}cd & d^2
\end{pmatrix}
.
\end{align*}

With this isomorphism, we often identify $\soq$ with $\somq$. 

Now we recall classical results on closed subgroups of 
$SU(2)$. 
Let us use the Pauli matrices 
\[
\sigma_1=
\begin{pmatrix}
0 & 1 \\
1 & 0
\end{pmatrix}
,\quad 
\sigma_2=
\begin{pmatrix}
0 & -i \\
i & 0
\end{pmatrix}
,\quad 
\sigma_3=
\begin{pmatrix}
1 &0 \\
0 & -1
\end{pmatrix}
. 
\]
They give an orthonormal basis for a three dimensional 
real vector space $\R^3=\R \sigma_1 +\R \sigma_2 +\R \sigma_3$. 
We define a double covering $\pi_1:SU(2)\longrightarrow SO(3)$ by 
the adjoint action 
$\pi_1(g)(a)=gag^{-1}$ for all $g\in SU(2)$ and $a\in \R^3$. 
For a subgroup $H\subset SO(3)$ 
we write $H^*$ for its inverse image by $\pi_1$ 
and such a subgroup is called a \textit{binary subgroup}. 
A symmetric group and alternative groups 
$A_4$, $S_4$ and $A_5$ are 
embedded into $SO(3)$ as 
the \textit{tetrahedral group}, the \textit{octahedral group} 
and the \textit{icosahedral group}, respectively. 
The group $D_m\,(2\leq m\leq \infty)$ is a \textit{dihedral group}, 
and $\T$ and $\T_m\,(2\leq m)$ are the ordinary torus and 
the cyclic groups of order $m$. 
In $SU(2)$ all the closed subgroups are conjugate to one of 
the \textit{trivial group} $1$, the \textit{cyclic group} $\T_n$, 
$SU(2)$, the maximal torus $\T$, 
the \textit{binary dihedral group} $D_m^* \,(2\leq m\leq \infty)$, 
the \textit{binary tetrahedral group} $A_4^*$, 
the \textit{binary octahedral group} $S_4^*$ 
and 
the \textit{binary icosahedral group} $A_5^*$. 
Since we need explicit embedding of $\T_n$ and $D_n$ in the final 
section, 
we prepare a few notations. 
Let us define rotation matrices
\[
r^{12}(\theta)
=
\begin{pmatrix}
e^{-i\frac{\theta}{2}} & 0\\
0 & e^{i\frac{\theta}{2}}
\end{pmatrix}
,\,
r^{23}(\theta)
=
\begin{pmatrix}
\cos\frac{\theta}{2} & -i\sin\frac{\theta}{2}\\
-i\sin\frac{\theta}{2} & \cos\frac{\theta}{2}
\end{pmatrix}
,\,
r^{13}(\theta)
=
\begin{pmatrix}
\cos\frac{\theta}{2} & \sin\frac{\theta}{2}\\
-\sin\frac{\theta}{2} & \cos\frac{\theta}{2}
\end{pmatrix}
\hspace{-1.3mm},
\]
where $\pi_1(r^{ij}(\theta))$ gives the rotation of angle $\theta$ in 
$\sigma_i$-$\sigma_j$ plane. 
The cyclic group $\T_n$ is specified by two 
angles $\chi$ and $\psi$ and denoted by $\T_{n}^{\chi,\psi}$. 
It consists of rotations of the angle $\frac{2\pi}{n}$ around 
the axis; $-\cos\chi\sin\psi\sigma_1
-\sin\chi\sin\psi\sigma_2+\cos\psi\sigma_3$. 
The dihedral group $D_n$ is specified by three angles 
$\phi, \chi$ and $\psi$ and 
denoted by $D_{n}^{\phi,\chi,\psi}$. 
It is generated by $\T_{n}^{\chi,\psi}$ and the rotation of the angle $\pi$ 
around the axis; 
$\cos\phi\cos\chi\cos\psi\sigma_1+\sin\phi\sin\chi\sigma_2
+\cos\phi\sin\psi\sigma_3$. 
By definition we have 
$\T_{n}^{\chi,\psi}
=\Ad(\pi_1(r^{12}(\chi)r^{13}(\psi)))(\T_{n}^{0,0})$ 
and 
$D_{n}^{\phi,\chi,\psi}
=\Ad(\pi_1(r^{12}(\chi)r^{13}(\psi)r^{12}(\phi)))(D_{n}^{0,0,0})$. 

In \cite{Podles2} he has classified all the quantum 
subgroups of $\suq$ for $-1\leq q <1$. 
We summarize it for readers' convenience and analysis on 
$\summ$ later. 

(1) $0<|q|<1$ case. 
Its quantum subgroup is one of $1$, $\T_m\,(m\geq2)$, 
$\T$, and $\suq$. 
For $\T_m\,(m\geq2)$ and $\T$, their restriction map is $\pi_\T$ or 
$\pi_\T'$. 

(2) $q=-1$ case. 
For 
$g
=
\begin{pmatrix}
\alpha & -\ovl{\gamma} \\
\gamma & \ovl{\alpha}
\end{pmatrix}
\in SU(2)$, 
define the $*$-homomorphism 
$\upsilon_{g}:\csum\longrightarrow \mathbb{B}(\C^2)$ by 
\[
\begin{pmatrix}
\ups_g(x) & \ups_g(u) \\
\ups_g(v) & \ups_g(y)
\end{pmatrix}
=
\begin{pmatrix}
\alpha \sigma_1 & \ovl{\gamma} \sigma_2 \\
\gamma \sigma_2 & \ovl{\alpha} \sigma_1
\end{pmatrix}
.
\]
Now we define the irreducible representation $\tau_g$ 
of $\csum$ as a $C^*$-algebra. 
They give all irreducible representations of $\csum$. 

(i) $\alpha,\gamma\neq0$ case. $\tau_g=\ups_g$ 

(ii) $\alpha=0$ case. 
$\tau_g:\csum\longrightarrow \C$ is 
\[
\begin{pmatrix}
\tau_g(x) & \tau_g(u) \\
\tau_g(v) & \tau_g(y)
\end{pmatrix}
=
\begin{pmatrix}
0 & \ovl{\gamma} \\
\gamma & 0
\end{pmatrix}
.
\]

(iii) $\gamma=0$ case. 
$\tau_g:\csum\longrightarrow \C$ is 
\[
\begin{pmatrix}
\tau_g(x) & \tau_g(u) \\
\tau_g(v) & \tau_g(y)
\end{pmatrix}
=
\begin{pmatrix}
\alpha & 0 \\
0 & \ovl{\alpha}
\end{pmatrix}
.
\]
For a compact subset $Z\subset SU(2)$ 
we consider the direct sum representation $\pi_Z=\oplus_{g\in Z}\pi_g$. 
In \cite[Proposition 2.4.]{Podles2}, 
it is characterized when 
$\pi_Z:\csum \longrightarrow\pi_Z(\csum)$ gives a restriction map 
to a quantum subgroup. 
We write $C(G_Z)=\pi_Z(\csum)$. 
Before a summary we study the embedding $\T_n\subset\summ$ a little. 
For $a\in \csum$ and $g\in\T_n$ we write simply $a(g)$ for $r(a)(g)$ 
where $r$ is a restriction map. 
Let $g$ be a generator of $\T_n$. 
Since we have $xv=-vx$, $x(g)=0$ or $v(g)=0$. 
We consider $v(g)\neq0$. 
By the definition of the restriction map we have 
$
\Big{(}
\begin{smallmatrix}
x(g^2)& u(g^2) \\
v(g^2)& y(g^2) \\
\end{smallmatrix}
\Big{)}
=
\Big{(}
\begin{smallmatrix}
1     & 0 \\
0     & 1 \\
\end{smallmatrix}
\Big{)}
$. 
This shows $n$ must be $2$. We denote the restriction by $\pi_{D_1}$. 
Hence if $n\geq 3$, the embedding of $\T_n$ is unique and 
if $n=2$, two embedding arises. 
Since in this paper notations $\T$ and $\T_n$ 
are used for the maximal torus 
and its cyclic subgroups, we prepare the notation $D_1$ for 
a subgroup of order $2$ which is embedded by $\pi_{D_1}$. 
Of course subgroups $D_1$ and $\T_2$ are isomorphic, 
however, right coideals by them are not isomorphic. 

(a) $C(G_Z)$ is an abelian $C^*$-algebra case
 (that is, $G_Z$ is an ordinary group), then 
$G_Z$ is (isomorphic to) one of trivial group $1$, 
the maximal torus $\T$, 
the cyclic groups $\T_n(\subset\T)$$(n\geq2)$, 
$D_1$ and 
dihedral groups $D_n$$(2\leq n \leq\infty)$ containing $\T_n$. 
(This isomorphism gives the conjugation $\beta_z^L$ 
on right coideals by them.) 
About $D_n$$(1\leq n \leq \infty)$ its corresponding subset $Z_{D_n}$ is 
$\{r^{12}(\frac{2\pi k}{n}),r^{12}(\frac{2\pi k}{n})r^{23}(\pi)
\mid 0\leq k \leq n-1\}$. 
Note that $Z_{D_n}$ is not a subgroup in $SU(2)$ if odd $n$. 
If $G_Z$ is not $\T_n$ or $D_n$ (odd $n$), $Z$ is a binary subgroup 
and it shows the spectral pattern of the right coideal 
$C(G_Z\setminus \summ)$ consists of integer spins. 
If $G_Z$ is, then it also has half integer spin parts. 

(b) If $C(G_Z)$ is a non-abelian $C^*$-algebra, 
$Z$ is one of binary subgroups 
whose image by $\pi_1$ is listed in \cite[p. 11]{Podles2}. 
Note that it depends on the embedding of a closed subgroup 
into $SO(3)$. 
Then it is known that 
$C(G_Z\setminus\summ)\subset C(SO_{-1}(3))$ and 
$\Xi_{-1}(C(G_Z\setminus\summ))=C(H\setminus SO_{1}(3))$. 

The quantum universal enveloping algebra $\uqsu$ 
is generated by four elements $k,k^{-1},e$ and $f$ which 
satisfy the following relations:
       \[kk^{-1}=1=k^{-1}k,\]
       \[ kek^{-1}=qe,\  kfk^{-1}=q^{-1}f,\]
       \[ ef-fe=\frac{k^2-k^{-2}}{q-q^{-1}}.\]
$\uqsu$ is realized as a Hopf $*$-subalgebra of 
the algebraic functional space $\asuq^*$ as follows.

(1) $q>0$ case. 
\[
k^{\pm 1}(w(\pi_\nu)_{r,s})=q^{\mp r} \delta_{r,s},
\]
\[
e(w(\pi_\nu)_{r,s})=\delta_{r+1,s}\sqrt{(\nu+s)_q(\nu-s+1)_q},
\]
\[
f(w(\pi_\nu)_{r,s})=\delta_{r-1,s}\sqrt{(\nu-s)_q(\nu+s+1)_q}. 
\]

(2) $q<0$ case. 
\[
k^{\pm 1}(w(\pi_\nu)_{r,s})=q^{\mp r} \delta_{r,s},
\]
\[
e(w(\pi_\nu)_{r,s})=\delta_{r+1,s}\sqrt{-1}^{-2\nu+1}
\sqrt{(\nu+s)_q(\nu-s+1)_q},
\]
\[
f(w(\pi_\nu)_{r,s})=\delta_{r-1,s}\sqrt{-1}^{-2\nu+1}
\sqrt{(\nu-s)_q(\nu+s+1)_q},
\]
where $q^n=\sqrt{-1}^{2n}(-q)^{n}$ for a half integer $n$. 
The Hopf $*$-algebra structure of $\uqsu$ is given by 
\begin{enumerate}

\item $k^*=k, \  e^*=f, \  f^*=e$,

\item $\hat{\delta}(k)=k\otimes k,\ 
      \hat{\delta}(e)=e\otimes k + k^{-1}\otimes e , \ 
      \hat{\delta}(f)=f\otimes k + k^{-1}\otimes f , \ $

\item $\hat{\kappa}(k)=k^{-1},\ 
      \hat{\kappa}(e)=-q e ,\ 
      \hat{\kappa}(f)=-q^{-1} f .$ 
\end{enumerate}
Hence $\uqsu$ is a Hopf $*$-subalgebra of $\asuq$. 
For a smooth corepresentation of $\asuq$ 
$\alpha:K \longrightarrow K\otimes \asuq$, 
we prepare $\uqsu$-module structure on $K$ by 
$\theta\cdot \xi=(\id\otimes\theta)(w \xi)$ 
for all $\theta\in\uqsu$ and $\xi\in K$. 
This representation of $\uqsu$ is called a differential representation 
of a corepresentation $(\alpha,K)$. 
Consider irreducible representation $(w(\pi_\nu),H_\nu)$ and 
we have the following formulae about its $\uqsu$-module structure.

(1) $q>0$ case. 
\[
k^{\pm 1}\cdot \xi_{r}^{\nu}=q^{\mp r}\xi_{r}^{\nu}\ ,
\]
\[
e\cdot \xi_{r}^{\nu}=\sqrt{(\nu+r)_q (\nu-r+1)_q}\,\xi_{r-1}^{\nu}\ ,
\]
\[
f\cdot \xi_{r}^{\nu}=\sqrt{(\nu-r)_q (\nu+r+1)_q}\,\xi_{r+1}^{\nu}\ .
\]

(2) $q<0$ case. 
\[
k^{\pm 1}\cdot \xi_{r}^{\nu}=q^{\mp r}\xi_{r}^{\nu}\ ,
\]
\[
e\cdot \xi_{r}^{\nu}=\sqrt{-1}^{-2\nu+1}
\sqrt{(\nu+r)_q (\nu-r+1)_q}\,\xi_{r-1}^{\nu}\ ,
\]
\[
f\cdot \xi_{r}^{\nu}=\sqrt{-1}^{-2\nu+1}
\sqrt{(\nu-r)_q (\nu+r+1)_q}\,\xi_{r+1}^{\nu}\ .
\]
for all $r\in I_\nu$. 
For a half integer $n$ 
a vector $\xi\in K$ is called a \textit{highest weight vector} 
of \textit{weight} $n$
if it satisfies $k\cdot\xi=q^n \xi$ and $e\cdot \xi=0$. 
For example, the vector $\xi_{-\nu}^{\nu}\in H_\nu$ is a 
highest weight vector of weight $\nu$. 
It is well-known that a tensor product $\uqsu$-module 
$H_\mu \otimes H_\nu$ is isomorphic to the direct sum $\uqsu$-module 
$\bigoplus_{|\mu-\nu|\leq \ell \leq \mu+\nu}H_{\ell}$, 
where $\ell$ runs through half integers. 
If we want to make a highest weight vector of $H_{\ell}$ 
from those of $H_\mu$ and $H_\nu$, 
the following well-known lemmas are useful. 

\begin{lem}\label{Clebsh}
For positive $q$, 
consider the tensor product $\uqsu$-module $H_\mu \otimes H_\nu$. 
Let $\ell$ be an integer with $|\mu-\nu|\leq \ell \leq \mu+\nu$. 
Define the coefficients $(C_{\mu,\nu}^{\ell})_r$ for 
$0 \leq r \leq \mu+\nu-\ell$ by
\[(C_{\mu,\nu}^{\ell})_r=
q^{-\frac{1}{2}(\ell+1)(\mu+\nu-\ell)}
(-q^{\ell+1})^{r}
\prod_{t=1}^{r}
\sqrt{\frac{(\mu+\nu-\ell+1-t)_q (\mu-\nu+\ell+t)_q}{(t)_q (2\nu-t+1)_q}}
 .\]
Then a vector 
$\eta^{\ell}=\sum_{r=0}^{\mu+\nu-\ell} 
(C_{\mu,\nu}^{\ell})_r\, \xi_{-\ell+\nu-r}^\mu \otimes \xi_{-\nu+r}^{\nu}$ 
is a highest weight vector of weight $\ell$ in $H_\mu \otimes H_\nu$. 
\end{lem}

\begin{proof}
The action of $\uqsu$ on $H_\mu \otimes H_\nu$ is given via 
coproduct. 
With this, we can easily justify $k\cdot\eta^{\ell}=q^{-\ell} \eta^{\ell}$ and 
$e\cdot \eta^{\ell}=0$. 
\end{proof}

We state a negative 
$q$ version of Lemma \ref{Clebsh} as follows. 

\begin{lem}\label{Clebsh negative}
For negative $q$, 
consider the tensor product $\uqsu$-module $H_\mu \otimes H_\nu$. 
Let $\ell$ be an integer with $|\mu-\nu|\leq \ell \leq \mu+\nu$ 
and $\xi^\mu$ and $\xi^\nu$ be vectors of copy of $\pi_\mu$ and 
$\pi_\nu$ respectively. 
Define the coefficients $(C_{\mu,\nu}^{\ell})_r$ for 
$0 \leq r \leq \mu+\nu-\ell$ by
\[(C_{\mu,\nu}^{\ell})_r=
q_0^{-\frac{1}{2}(\ell+1)(\mu+\nu-\ell)}
(-1)^{(-\mu+\nu+\ell)r}
q_0^{r(\ell+1)}
\prod_{t=1}^{r}
\sqrt{\frac{(\mu\!+\!\nu\!-\!\ell\!+\!1\!-\!t)_q 
(\mu\!-\!\nu\!+\!\ell\!+\!t)_q}{(t)_q (2\nu\!+\!1\!-\!t)_q}}
 .\]
Then a vector 
$\eta^\ell=\sum_{r=0}^{\mu+\nu-\ell} 
(C_{\mu,\nu}^{\ell})_r\, \xi_{-\ell+\nu-r}^\mu \otimes \xi_{-\nu+r}^{\nu}$ 
is a highest weight vector of weight $\ell$ in $H_\mu \otimes H_\nu$.
\end{lem}
Let $Y_\mu$ be a linear space of 
$\pi_\mu$-eigenvectors of $\csuq$. 
We prepare the notation of eigenvectors 
$\bw^\mu_r:=(w(\pi_\mu)_{r,t})_{t\in I_{\mu}}$ for 
$r\in I_{\mu}$. 
They give an orthonormal basis of $Y_{\mu}$. 
Let us consider a covariant system $(A,\suq,\alpha)$. 
For its eigenvector spaces $\{X_\nu\}_{\nu\in \frac{1}{2}\Z_{\geq0}}$ 
we define the \textit{product of eigenvectors} 
$\Psi_{\ell}:X_\mu \times X_\nu \longrightarrow X_{\ell}$ 
by using Lemma \ref{Clebsh}:
\[
\Psi_{\ell}(\xi^\mu,\xi^\nu)_{-\ell}=\sum_{r=0}^{\mu+\nu-\ell} 
(C_{\mu,\nu}^{\ell})_r\, \xi_{-\ell+\nu-r}^\mu \,\xi_{-\nu+r}^{\nu}
\]
for all eigenvectors $\xi^\mu$ and $\xi^\nu$, 
where the coefficients $(C_{\mu,\nu}^{\ell})_r$ are 
given in Lemma \ref{Clebsh} and Lemma \ref{Clebsh negative}. 
In \cite{Podles1}, he has classified ergodic systems 
$\{A,\suq,\alpha\}$ with 
$\dim A_{\pi_1}=1$ and $A=C^*(A_{\pi_1})$ 
for $|q|<1$. 
It is uniquely determined up to $\suq$-isomorphism and 
is called the \textit{quantum sphere}. 
We summarize his classification. 
Let $X_{\pi_1}$ be the $\pi_1$-eigenvector space. 
Since this is one-dimensional, we can take 
the (unique) $\pi_1$-eigenvector 
$\xi$ with the following properties:
\[
T\xi_\lambda=\xi_\lambda,\quad 
(\xi_\lambda,\xi_\lambda)=1,\quad 
\Psi_1(\xi_\lambda, \xi_\lambda)= \lambda\, \xi_\lambda,
\]
where $\lambda$ is a non-negative real constant. 
From the second assumption the $C^*$-algebra $A$ is generated by 
$(\xi_\lambda)_{-1},$ $(\xi_\lambda)_0$ and $(\xi_\lambda)_1$. 
Write $\lambda_0$ for $(q^{-1}-q)^{-1}\lambda$. 
For $n\in \Z_{\geq1}$ we define a positive number 
$c(n)=\frac{q^{n+1}+q^{-n-1}}{\sqrt{(n)_q (n+2)_q}}$. 
Then the classification is done as follows.

Case 1: If $\lambda_0>1$, 
then there exists $n\in \Z_{\geq1}$ such that 
$\lambda=c(n)$ and we obtain a $G$-isomorphism 
$A\cong \End(H_{\pi_{\frac{n}{2}}})$. 
In this case the spectral pattern of $A$ is 
$\pi_0 \oplus \pi_{1} \oplus \cdots \oplus \pi_{n}$. 

Case 2: If $0\leq \lambda_0\leq 1$ and $q>0$, 
then the map 
$\xi_\lambda=((\xi_\lambda)_{-1},(\xi_\lambda)_0,(\xi_\lambda)_1) \mapsto 
\Big{(}q^{\frac{1}{2}}\sqrt{\frac{1-\lambda_0^2}{(2)_q}}
,\lambda_0
,-q^{-\frac{1}{2}}\sqrt{\frac{1-\lambda_0^2}{(2)_q}}
\Big{)}
\cdot w(\pi_1)$ 
gives a $G$-equivariant embedding $A\hookrightarrow \csuq$. 
In this case the spectral pattern of $A$ is 
$\oplus_{\ell\in \Z_{\geq0}} \pi_{\ell}$. 

Case 3: If $0\leq \lambda_0\leq 1$ and $q<0$, 
then the map 
$\xi_\lambda=((\xi_\lambda)_{-1},(\xi_\lambda)_0,(\xi_\lambda)_1) \mapsto 
\Big{(}(-q)^{\frac{1}{2}}\sqrt{\frac{1-\lambda_0^2}{(2)_q}}
,\lambda_0
,(-q)^{-\frac{1}{2}}\sqrt{\frac{1-\lambda_0^2}{(2)_q}}
\Big{)}
\cdot w(\pi_1)$ 
gives a $G$-equivariant embedding $A\hookrightarrow \csuq$. 
In this case the spectral pattern of $A$ is 
$\oplus_{\ell\in \Z_{\geq0}} \pi_{\ell}$. 

In this paper we use the notation $C(S_{q,\lambda}^2)$ for 
the right coideal which is defined in the above case $2$ and $3$. 
In the proof of Lemma \ref{pi1 of tn} we obtain a $\pi_1$-eigenvector 
of $C(\T^{0,\psi}\setminus SO(3))$. 
Comparing it with the above embedding, 
the quantum sphere $C(S_{q,\lambda}^2)$ is considered 
a $q$-deformation of $C(\T^{0,\psi}\setminus SO(3))$ 
with the parameter $\lambda_0=\cos \psi$. 
If $\lambda_0=1$, $C(S_{q,(q^{-1}-q)^{-1}})$ becomes 
the canonical homogeneous sphere $C(\T \setminus \suq)$. 
Let $\beta_z^{L}$ be a left action of the maximal torus 
on $\csuq$ defined by 
$(\ev_{z}\circ\pi_{\T}\otimes\id)\circ \delta$ for all $z\in \T$. 
It satisfies $(\beta_z^L\otimes\id)\circ\delta=\delta\circ\beta_z^L$, 
that is, $\beta_z^L$ gives an $\suq$-isomorphism. 
Note that 
all embeddings of quantum spheres into $\csuq$ are obtained by 
rotations of $\beta^{L}$ for the above given embedding. 
Similarly we define the right action $\beta^R$ of the maximal torus 
by $\beta^R=(\id\otimes\pi_\T)\circ\delta$. 

We show the list of all the connected graphs of norm $2$ in the Appendix. 
For its classification readers are referred to 
\cite[Lemma 1.4.1]{GHJ}. 
The labels except for $A_m'$ correspond to 
closed subgroups of $SU(2)$ or $SU_{-1}(2)$ by their 
McKay diagrams about 
fundamental representation $\pi_{\frac{1}{2}}$ 
or multiplicity diagrams of the 
right coideals obtained by quotients.  
We will see later that type $A_m'$ does not occur 
even if we investigate right coideals. 

\begin{rem}
We discuss the absence of $A_m'\,(3\leq m\leq\infty)$. 
Let $C(G)$ be a compact quantum group which has an 
irreducible unitary representation $w$ generating $R(G)$
 (that is, $C(G)$ is a compact matrix pseudogroup 
\cite{Woronowicz3}). 
We consider a compact quantum subgroup $C(H)$ whose McKay diagram 
about $w$ is of type $A_m'\,(m\geq3)$. 
We simply denote the restriction of $w$ onto $H$ by $w|_H$. 
Now vertices in $A_m'$ corresponds to the irreducible 
representations of $H$ and let 
$\rho_0,\rho_1\ldots$ be the irreducible representations from 
the left side in the Figure \ref{A_m'}, \ref{A_infty'}. 
Since the entries of the Perron-Frobenius eigenvector 
are equal to the dimensions of the corresponding irreducible modules, 
$\wdh{H}$ must be a group. 
By definition of McKay diagram, 
$ w|_H\cdot\rho_0=\rho_0+\rho_1$ holds. 
Hence we have $w|_H=0+\rho_1\rho_0^{-1}$ 
where $0$ is the trivial representation. 
Then it yields $w|_H\cdot\rho_1=\rho_1+\rho_1\rho_0^{-1}\rho_1$, 
and this shows there must exist a single loop at $\rho_1$. 
Hence the McKay diagram of type $A_m'$ $(m\geq3)$ 
does not appear. 
\end{rem}

\section{Classification of right coideals of $\csuq$: $0<q<1$ case}

In $q=1$ case the main theorem of \cite[p.\ 309]{Wassermann3} says that 
any ergodic system of $SU(2)$ is an induced system 
$\{\End(W),H,\beta\}$ where $H$ is a closed subgroup of $SU(2)$. 
A closed subgroup $H$ is conjugate to one of 
$1,\T_n\,(n\geq2),\T,SU(2),D_n^*\,(n\geq 2),D_\infty^*,A_4^*,S_4^*$ 
and $A_5^*$. 
Their corresponding multiplicity diagrams are exactly 
their McKay diagrams. 
In $0<q<1$ case there are a little quantum subgroups of $\suq$ 
\cite[Theorem 2.1.]{Podles2}. 
So we are interested in the absence of some multiplicity diagrams in the 
table. 
Finally we obtain the following results for right coideals. 

\begin{thm}\label{diagram of right coideal 0<q<1}
Let $A\subset \csuq$ be a right coideal. 
Then its multiplicity diagram is one of type 
$1, \T_n\,(n\geq2), \T, SU(2)$ and $D_\infty^*$. 
If it is of type $\T$, then it is one of the quantum spheres. 
Otherwise it is unique 
up to conjugation by $\beta^L$. 
\end{thm}

Now we start the proof of Theorem \ref{diagram of right coideal 0<q<1}. 
First we show that there do not exist right coideals 
whose multiplicity diagrams are of type 
$D_n^*\,(n\geq 2),D_\infty^*,A_4^*,S_4^*$ 
and $A_5^*$. 
Second we classify right coideals whose multiplicity diagrams are 
of type $1, \T_n\,(n\geq2), \T, SU(2), D_\infty^*$. 
In all cases we make use of products of eigenvectors. 

(I) $A_4^*$ case. $A$ has a spectral pattern 
$\pi_0\oplus \pi_3 \oplus \pi_4\oplus 2\pi_6\oplus \pi_7 \oplus \cdots$. 
Like the discussion of \cite[p. 321]{Wassermann3}, 
we focus on spectral gaps: $\pi_1,\pi_2$ and $\pi_5$. 
We will soon notice the importance of using the both of 
even and odd spin. 
Let $\eta=(\eta_r)_{r\in I_3}$ be a self-conjugate 
$\pi_3$-eigenvector of $A$. 
Take scalars $\{c_r\}_{r\in I_3}$ such that 
$\eta_r=\sum_{s\in I_3}c_s w(\pi_3)_{s,r}$ for all $r\in I_3$. 
Applying Lemma \ref{Clebsh}, we obtain the following 
highest weight vectors
\[
\Psi_5(\eta,\eta)_{-5}
=q^{-3} \eta_{-2}\eta_{-3}-q^3 \eta_{-3}\eta_{-2}\, ,
\]
\begin{align*}
\Psi_2(\eta,\eta)_{-2}
=&\,
q^{-6}\eta_{1}\eta_{-3}
-q^{-3}\sqrt{\frac{(4)_q(3)_q}{(6)_q}}\eta_{0}\eta_{-2}
+\frac{(4)_q(3)_q}{\sqrt{(6)_q(5)_q(2)_q}}\eta_{-1}\eta_{-1}\\
&\,
-
q^{3}\sqrt{\frac{(4)_q(3)_q}{(6)_q}}\eta_{-2}\eta_{0}
+
q^6\eta_{-3}\eta_{1}\, ,
\end{align*}

\begin{align*}
\Psi_1(\eta,\eta)_{-1}
=&\,
 q^{-5}\eta_{2}\eta_{-3} - q^{-3} 
\sqrt{\frac{(5)_q(2)_q}{(6)_q}}\eta_{1}\eta_{-2}\\
&\,
+q^{-1}\sqrt{\frac{(4)_q(3)_q}{(6)_q}} \eta_{0}\eta_{-1}
-q \sqrt{\frac{(4)_q(3)_q}{(6)_q}}\eta_{-1}\eta_{0}\\
&\,
+ 
q^3\sqrt{\frac{(5)_q(2)_q}{(6)_q}}\eta_{-2}\eta_{1}
-q^5\eta_{-3}\eta_{2} \, .
\end{align*}

They are actually $0$ vectors because of the absence of 
spectra. 
Applying the lowering operator $f\in \uqsu$ to these vectors, 
we obtain
\[
f\cdot\Psi_5(\eta,\eta)_{-5}=\sqrt{(5)_q(2)_q}\eta_{-1}\eta_{-3}
+\sqrt{(6)_q}(q^{-5}-q^5)\eta_{-2}^2
-\sqrt{(5)_q(2)_q}\eta_{-3}\eta_{-1}\, ,
\]
$
f^2\cdot\Psi_5(\eta,\eta)_{-5}=
q^3\sqrt{(4)_q(3)_q}\eta_{0}\eta_{-3}
+
(q^{-3}+q^{-1}-q^7)\sqrt{(6)_q(5)_q(2)_q}\eta_{-1}\eta_{-2}$
\begin{center}
$
{}\qquad\qquad+
(q^{-7}-q-q^3)\sqrt{(6)_q(5)_q(2)_q}\eta_{-2}\eta_{-1}+
q^{-3}\sqrt{(4)_q(3)_q}\eta_{-3}\eta_{0}\, 
$
\end{center}
and 

$
\sqrt{(6)_q}f\cdot\Psi_1(\eta,\eta)_{-1}=
q^{-2}(6)_q\eta_{3}\eta_{-3}
+\big{(}-q^{-1}(5)_q(2)_q +q^{-3}(6)_q \big{)}
\eta_{2}\eta_{-2}$
\begin{center}
$\qquad
+\big{(}(4)_q(3)_q-q^{-2}(5)_q(2)_q\big{)}
\eta_{1}\eta_{-1}
+
(q^{-1}-q)(4)_q(3)_q \eta_{0}^2
$
\end{center}
\begin{center}
$\qquad\quad\quad\!\!\,
+
\big{(}q^2(5)_q(2)_q-(4)_q(3)_q\big{)}
\eta_{-1}\eta_{1}
+
\big{(}q(5)_q(2)_q-q^3 (6)_q\big{)}
\eta_{-2}\eta_{2}
-
$
\end{center}

$
\qquad\qquad\qquad
-q^2 (6)_q\eta_{-3}\eta_{3}\, .
$

Recall a surjective $*$-homomorphism 
$\pi_{\T}:\csuq \longrightarrow C(\T)$. 
It sends $w(\pi_\nu)_{r,s}$ to $\delta_{r,s}z^{-2r}$ for 
$\nu\in \frac{1}{2}\Z_{\geq0}$ and $r,s \in I_\nu$. 
Then we have $\pi_\T(\eta_r)=z^{-2r}c_r$ for $r\in I_3$. 
From $0=\pi_\T(f\cdot\Psi_5(\eta,\eta)_{-5})$, we obtain $c_{-2}=0$. 
Then from $0=\pi_\T(\Psi_1(\eta,\eta)_{-1})
=\pi_\T(f^2\cdot\Psi_5(\eta,\eta)_{-5})$, 
we have $c_{-3}c_0=c_{-1}c_0=0$. 
If $c_0$ is not $0$, then $c_{-3}$ and $c_{-1}$ are equal to $0$. 
Hence $\{c_{r}\}_{r\in I_3}$ are all zero except for $r=0$. 
From $0=\pi_\T(f\cdot\Psi_1(\eta,\eta)_{-1})$ 
we have $c_0=0$, this is contradiction.
Therefore $c_0$ must be $0$. 
Then the last equality deduces to 
$(6)_q c_3c_{-3}=(5)_q(2)_q c_1c_{-1}$. 
Because of the self-conjugacy of $\eta$ we have 
$c_{-3}=-q^3\overline{c_3}$ and $c_{-1}=-q \overline{c_1}$. 
Then we have 
\begin{equation*}
(6)_q |c_3|^2=q^{-2}(5)_q(2)_q|c_1|^2\, . 
\end{equation*}
From the second one we get 
\begin{equation*}
(q^{-6}+q^6)c_1c_{-3}+\frac{(4)_q(3)_q}{\sqrt{(6)_q(5)_q(2)_q}}c_{-1}^2=0\, .
\end{equation*}
We can easily see that there does not exist a solution of 
the above two equalities except for the case $c_3=c_1=0$. 
Therefore we have proved no existence of right coideal of type $A_4^*$.

\begin{rem}
If we consider the $q=1$ or $q=-1$ case, then products of integer spin 
eigenvectors to odd integer spin eigenvectors become a trivial zero map. 
Therefore arguments in the above and below are non-trivial only in 
the case of $|q|\neq1$. 
\end{rem}

(II) $S_4^*$ case.  $A$ has a spectral pattern 
$\pi_0 \oplus \pi_4 \oplus \pi_6 \oplus \pi_8 \oplus \pi_9 
\oplus \pi_{10} \oplus \cdots$. 
There are spectral gaps of $\pi_\nu$ for $\nu=1,2,3,5,7$. 
Let $\eta=(\eta_s)_{s\in I_4}=
\sum_{s\in I_4}c_s \bw_s^4$ be a non zero self-conjugate 
$\pi_4$-eigenvector of $A$. 
We derive contradiction by showing the coefficients 
$c_s$ are all zero. 
Let us consider eigenvectors
\begin{align*}
\Psi_7(\eta,\eta)_{-7}
=&\,
q^{-4}\eta_{-3}\eta_{-4}-q^{4}\eta_{-4}\eta_{-3}\, ,\\
\Psi_5(\eta,\eta)_{-5}
=&\, 
q^{-9}\eta_{-1}\eta_{-4}
-q^{-3}\sqrt{\frac{(6)_q(3)_q}{(8)_q}}\eta_{-2}\eta_{-3}
+q^{3}\sqrt{\frac{(6)_q(3)_q}{(8)_q}}\eta_{-3}\eta_{-2}\\
&\quad
-q^{9}\eta_{-4}\eta_{-1} ,\\
\Psi_3(\eta,\eta)_{-3}
=&\,
q^{-10}\eta_{1}\eta_{-4}
-q^{-6}\sqrt{\frac{(5)_q(4)_q}{(8)_q}}\eta_{0}\eta_{-3} 
+q^{-2}\frac{(5)_q(4)_q}{\sqrt{(8)_q(7)_q(2)_q}}\eta_{-1}\eta_{-2}\\
&\quad 
-q^2 \frac{(5)_q(4)_q}{\sqrt{(8)_q(7)_q(2)_q}}\eta_{-2}\eta_{-1} 
+q^6\sqrt{\frac{(5)_q(4)_q}{(8)_q}}\eta_{-3}\eta_{0}
-q^{10}\eta_{-4}\eta_{1} , \\
\Psi_2(\eta,\eta)_{-2}
=&\,
q^{-9}\eta_{2}\eta_{-4}
-q^{-6}\sqrt{\frac{(6)_q(3)_q}{(8)_q}}\eta_1\eta_{-3}\\
&\quad
+q^{-3}\sqrt{\frac{(6)_q(5)_q(4)_q(3)_q}{(8)_q(7)_q(2)_q}}\eta_0\eta_{-2}
-\frac{(5)_q(4)_q}{(8)_q(7)_q(2)_q}\eta_{-1}^2\\
&\quad
+q^{3}\sqrt{\frac{(6)_q(5)_q(4)_q(3)_q}{(8)_q(7)_q(2)_q}}\eta_{-2}\eta_{0}
-q^{6}\sqrt{\frac{(6)_q(3)_q}{(8)_q}}\eta_1\eta_{-3}
+q^9 \eta_{-4}\eta_2,
\end{align*}
\begin{align*}
\Psi_1(\eta,\eta)_{-1}
=&\,
q^{-7}\eta_{3}\eta_{-4}
-q^{-5}\sqrt{\frac{(7)_q(2)_q}{(8)_q}}\eta_{2}\eta_{-3}
+q^{-3}\sqrt{\frac{(6)_q(3)_q}{(8)_q}}\eta_{1}\eta_{-2}\\
&\quad
-q^{-1}\sqrt{\frac{(5)_q(4)_q}{(8)_q}}\eta_{0}\eta_{-1}
+q\sqrt{\frac{(5)_q(4)_q}{(8)_q}}\eta_{-1}\eta_{0}
-q^{3}\sqrt{\frac{(6)_q(3)_q}{(8)_q}}\eta_{-2}\eta_{1}\\
&\quad
+q^{5}\sqrt{\frac{(7)_q(2)_q}{(8)_q}}\eta_{-3}\eta_{2}
-q^{7}\eta_{-4}\eta_{3}.
\end{align*}
They are in fact $0$ vectors by the gap at each spectrum. 
Applying lowering operator $f$ to these vectors, 
we get 
\begin{align*}
f\cdot\Psi_7(\eta,\eta)_{-7}
=&\,
\sqrt{(7)_q(2)_q}\eta_{-2}\eta_{-4}+(q^{-7}-q^7)\sqrt{(8)_q}\eta_{-3}^2
-\sqrt{(7)_q(2)_q}\eta_{-4}\eta_{-2},
\end{align*}
\begin{align*}
f\cdot\Psi_5(\eta,\eta)_{-5}
=&\,
q^{-5}\sqrt{(5)_q(4)_q}\eta_{0}\eta_{-4}
+\frac{q^{-10}(8)_q-(6)_q(3)_q}{\sqrt{(8)_q}}\eta_{-1}\eta_{-3}\\
&\! 
+(q^5-q^{-5})\sqrt{\frac{(7)_q(6)_q(3)_q(2)_q}{(8)_q}}\eta_{-2}\eta_{-2}
+\frac{-q^{10}(8)_q+(6)_q(3)_q}{\sqrt{(8)_q}}\eta_{-3}\eta_{-1}\\
&\!
-q^{5}\sqrt{(5)_q(4)_q}\eta_{-4}\eta_{0},
\end{align*}
\begin{align*}
f\cdot\Psi_3(\eta,\eta)_{-3}
=&\,
q^{-6}\sqrt{(6)_q(3)_q}\eta_{2}\eta_{-4}
+
\frac{q^{-9}\sqrt{(8)_q}-q^{-3}(5)_q(4)_q}{\sqrt{(8)_q}}\eta_{1}\eta_{-3}\\
&\quad
+
((5)_q(4)_q-q^{-6}(7)_q(2)_q)
\frac{(5)_q(4)_q}{(8)_q(7)_q(2)_q}\eta_{0}\eta_{-2}\\
&\quad
+
(q^{-3}-q^3)(5)_q(4)_q\frac{(6)_q(3)_q}{(8)_q(7)_q(2)_q}\eta_{-1}\eta_{-1} \\
&\quad
+
(q^6(7)_q(2)_q-(5)_q(4)_q)\frac{(5)_q(4)_q}{(8)_q(7)_q(2)_q}\eta_{-2}\eta_{0}\\
&\quad
+
\frac{q^{3}(5)_q(4)_q-q^{9}\sqrt{(8)_q}}{\sqrt{(8)_q}}\eta_{-3}\eta_{1}
-
q^6\sqrt{(6)_q(3)_q}\eta_{-4}\eta_{2} ,
\end{align*}
\begin{align*}
f\cdot\Psi_1(\eta,\eta)_{-1}
=&\,
q^{-3}\sqrt{(8)_q}\eta_{4}\eta_{-4}
+
\frac{q^{-4}(8)_q-q^{-2}(7)_q(2)_q}{\sqrt{(8)_q}}\eta_{3}\eta_{-3}\\
&\quad
+
\frac{q^{-1}(6)_q(3)_q-q^{-3}(7)_q(2)_q}{\sqrt{(8)_q}}\eta_{2}\eta_{-2} 
+
\frac{q^{-2}(6)_q(3)_q-(5)_q(4)_q}{\sqrt{(8)_q}}\eta_{1}\eta_{-1}\\
&\quad
+
\frac{(q-q^{-1})(5)_q(4)_q}{\sqrt{(8)_q}}\eta_{0}\eta_{0}
+
\frac{(5)_q(4)_q-q^2(6)_q(3)_q}{\sqrt{(8)_q}}\eta_{-1}\eta_{1} \\
&\quad
+
\frac{q^3(7)_q(2)_q-q(6)_q(3)_q}{\sqrt{(8)_q}}\eta_{-2}\eta_{2}
+
\frac{q^2(7)_q(2)_q-q^4(8)_q}{\sqrt{(8)_q}}\eta_{-3}\eta_{3} \\
&\quad
-
q^3\sqrt{(8)_q}\eta_{-4}\eta_{4}\, .
\end{align*}
Applying $\pi_\T$ to these vectors, we obtain the following 
equations: 
\begin{align}
&
c_{-3}
=0\, , \label{seqn01}\\
&(q^{-5}-q^5)\sqrt{(5)_q(4)_q}c_0 c_{-4}
+(q^5-q^{-5})\sqrt{\frac{(7)_q(6)_q(3)_q(2)_q}{(8)_q}}c_{-2}^2
=0 \label{seqn02}\, ,\\ 
&
(q^{-6}-q^6)\sqrt{(6)_q(3)_q}c_{2}c_{-4}
+
(q^6-q^{-6})\frac{(5)_q(4)_q}{(8)_q} c_{0}c_{-2} \label{seqn03}\\
&\hspace{30mm}
+
(q^{-3}-q^3)(5)_q(4)_q \frac{(6)_q(3)_q}{(8)_q(7)_q(2)_q} c_{-1}^2
=0\, , \notag\\
&(q^{-3}-q^3)\sqrt{(8)_q}c_{4}c_{-4}
+
\frac{(q^{-1}-q)(6)_q(3)_q-(q^3-q^{-3})(7)_q(2)_q}{\sqrt{(8)_q}}
c_{2}c_{-2} \label{seqn04} \\
&\hspace{30mm}
+
\frac{(q^{-2}-q^2)(6)_q(3)_q}{\sqrt{(8)_q}}
c_{1}c_{-1}
+
\frac{(q-q^{-1})(5)_q(4)_q}{\sqrt{(8)_q}}c_{0}^2 
=0 \, . 
\notag
\end{align}
Also noticing $c_{-3}=0$, we get the following equations 
via original eigenvectors $\Psi_\nu(\eta,\eta)$ 
for $\nu=5,3,2$:
\begin{align}
&\hspace{-1.2mm}
(q^{-9}-q^9)c_{-1}c_{-4}
= 0 \, ,\label{seqn05}\\
&\hspace{-1.2mm}
(q^{-10}-q^{10})c_{1}c_{-4}+(q^{-2}-q^2)
\frac{(5)_q(4)_q}{(8)_q(7)_q(2)_q}c_{-1}c_{-2}
= 0 \, ,\label{seqn06}\\
&\hspace{-1.2mm}
(q^{-9}+q^9)c_{2}c_{-4}
+
(q^3+q^{-3})
\sqrt{\frac{(6)_q(5)_q(4)_q(3)_q}{(8)_q(7)_q(2)_q}}c_0 c_{-2}
-
\frac{(5)_q(4)_q}{(8)_q(7)_q(2)_q}c_{-1}^2
= 0 .\label{seqn07}
\end{align}
From (\ref{seqn05}) $c_{-1}$ or $c_{-4}$ are $0$. 
If $c_{-4}$ is $0$, then we have $c_{-2}=c_{-1}=c_{0}=0$ by 
(\ref{seqn02}), (\ref{seqn03}) and (\ref{seqn04}). 
This is contradiction. 
Hence $c_{-4}$ is not $0$ and $c_{-1}$ 
is equal to $0$. 
Next we can derive $c_2c_{-4}=c_0 c_{-2}=0$ by (\ref{seqn03}) 
and (\ref{seqn07}). 
Then we have $c_{-2}=0$. 
This yields $c_0=0$ through (\ref{seqn02}). 
Finally we get $c_{-4}=0$ by (\ref{seqn04}), however, this is a contradiction.

(III) $A_5^*$ case. A $C^*$-algebra $A$ has a spectral pattern 
$\pi_0 \oplus \pi_6 \oplus \pi_{10} \oplus \pi_{12} \oplus \cdots$. 
Let $\eta=(\eta_s)_{s\in I_6}=\sum_{s\in I_6}c_s \bw_s^6$ be a 
non-zero self-conjugate $\pi_6$-eigenvector. 
Equality $c_{-s}=(-q)^s\overline{c_s}$ follows from the self-conjugacy. 
We show coefficients $c_s$ are all zero and derive contradiction. 
In order to do it we make use of spectral gaps at $\pi_\nu$ for 
$\nu=11,9,8,7,4,2$ and $1$. 

First we see the gap at $11$. 
We construct a $\pi_{11}$-eigenvector $\Psi_{11}(\eta,\eta)$ 
and apply the lowering operator $f$ to $\Psi_{11}(\eta,\eta)_{-11}$. 
Then we have the following equalities
\begin{align*}
\Psi_{11}(\eta,\eta)_{-11}
=&
q^{-6}\eta_{-5}\eta_{-6}-q^{6}\eta_{-6}\eta_{-5}\, ,\\
f\cdot\Psi_{11}(\eta,\eta)_{-11}
=&
\sqrt{\!(11)_q(2)_q}\hspace{0.2mm}\eta_{-4}\eta_{-6}
+(q^{-11}\!-q^{11})
\sqrt{\!(12)_q}\hspace{0.2mm}\eta_{-5}^2
-\sqrt{\!(11)_q(2)_q}\hspace{0.2mm}\eta_{-6}\eta_{-4}.
\end{align*}
Applying $\pi_\T$ to the second one, we get $c_{-5}=0\,$. 

Next we see the gap at $9$. 
We make $\Psi_{9}(\eta,\eta)_{-9}$ and 
$f\cdot\Psi_{9}(\eta,\eta)_{-9}$ as follows. 
\begin{align*}
\sqrt{(12)_q}\,\Psi_{9}(\eta,\eta)_{-9}
=&\hspace{0.3mm}
q^{-15}\sqrt{(12)_q}\eta_{-3}\eta_{-6}
-q^{-5}\sqrt{(10)_q(3)_q}\eta_{-4}\eta_{-5} \\
&\quad 
+q^5\sqrt{(10)_q(3)_q}\eta_{-5}\eta_{-4} 
-q^{15}\sqrt{(12)_q}\eta_{-6}\eta_{-3} \\
\sqrt{(12)_q}\,f\cdot\Psi_{9}(\eta,\eta)_{-9}
=&\hspace{0.3mm}
q^{-9}\sqrt{(12)_q(9)_q(4)_q}\eta_{-2}\eta_{-6}
+\big{(}q^{-18}(12)_q\!-\!(10)_q(3)_q\big{)}\eta_{-3}\eta_{-5} \\
&
+(q^9\!-\!q^{-9})\sqrt{(11)_q(10)_q(3)_q(2)_q}\eta_{-4}\eta_{-4} \\
&
+\big{(}(10)_q(3)_q\!-\!q^{18}(12)_q\big{)}\eta_{-5}\eta_{-3} 
-q^{9}\sqrt{(12)_q(9)_q(4)_q}\eta_{-6}\eta_{-2}.
\end{align*}
Applying $\pi_\T$ to the first one, we get $c_{-3}c_{-6}=0$. 
Similarly from the second one we get the following 
equality
\begin{equation}\label{a501}
\sqrt{(12)_q(9)_q(4)_q}c_{-2}c_{-6}
-
\sqrt{(11)_q(10)_q(3)_q(2)_q}c_{-4}^2
=0\, .
\end{equation}

Next we look at the gap at $8$. 
In this case we need two operators
\begin{align*}
&\sqrt{(12)_q(11)_q(2)_q}\,\Psi_8(\eta,\eta)_{-8}\\
=&\,
q^{-18}\sqrt{(12)_q(11)_q(2)_q}\eta_{-2}\eta_{-6}
-q^{-9}\sqrt{(11)_q(9)_q(4)_q(2)_q}\eta_{-3}\eta_{-5}\\
&\quad
+\sqrt{(10)_q(9)_q(4)_q(3)_q}\eta_{-4}\eta_{-4}
-q^{9}\sqrt{(11)_q(9)_q(4)_q(2)_q}\eta_{-5}\eta_{-3}\\
&\quad
+q^{18}\sqrt{(12)_q(11)_q(2)_q}\eta_{-6}\eta_{-2}\, ,\\
&
\sqrt{(12)_q(11)_q(2)_q}\,f^2\cdot\Psi_8(\eta,\eta)_{-8}\\
=&\hspace{0.3mm}
q^{-6}\sqrt{(12)_q(11)_q(8)_q(7)_q(6)_q(5)_q(2)_q}\eta_{0}\eta_{-6}\\
&\hspace{-0.4mm}
+
\!\big{(}(q^{-13}+q^{-15})(12)_q-q(9)_q(4)_q\big{)}
\sqrt{(11)_q(8)_q(5)_q(2)_q}\eta_{-1}\eta_{-5}\\
&\hspace{-0.4mm}
+\!
\big{\{}q^8(10)_q(9)_q(4)_q(3)_q
\!+\!q^{-22}(12)_q(11)_q(2)_q \!-\!
(q^{-6}\!+\!q^{-8})(11)_q(8)_q(4)_q(2)_q\big{\}}\eta_{-2}\eta_{-4}\\
&\hspace{-0.4mm}
+\!
\big{(}(10)_q(3)_q-(q^{-15}+q^{15})(11)_q\big{)}(2)_q
\sqrt{(10)_q(9)_q(4)_q(3)_q}
\eta_{-3}\eta_{-3}\\
&\hspace{-0.4mm}
+\!
\big{\{}q^{-8}(10)_q(9)_q(4)_q(3)_q
+q^{22}(12)_q(11)_q(2)_q
-(q^{6}+q^{8})(11)_q(8)_q(4)_q(2)_q\big{\}}
\eta_{-4}\eta_{-2}\\
&\hspace{-0.4mm}
+\!
\big{(}(q^{13}+q^{15})(12)_q-q^{-1}(9)_q(4)_q\big{)}
\sqrt{(11)_q(8)_q(5)_q(2)_q}\eta_{-5}\eta_{-1} \\
&\hspace{-0.4mm}
+\!
q^{6}\sqrt{(12)_q(11)_q(8)_q(7)_q(6)_q(5)_q(2)_q}\eta_{-6}\eta_{0}
.
\end{align*}
From the first equality we get the following one
\begin{equation}\label{a502}
(q^{-18}+q^{18})\sqrt{(12)_q(11)_q(2)_q}c_{-2}c_{-6}
+
\sqrt{(10)_q(9)_q(4)_q(3)_q}c_{-4}^2
=0\, .
\end{equation}
Equalities (\ref{a501}) and (\ref{a502}) shows 
$c_{-2}c_{-6}=0$ and $c_{-4}=0$. 
Hence by the above second equality we get
\begin{align}
&(q^{-6}+q^6)\sqrt{(12)_q(11)_q(8)_q(7)_q(6)_q(5)_q(2)_q}c_{0}c_{-6}
\label{a503}\\
&\hspace{12mm}
+
\big{(}(10)_q(3)_q-(q^{-15}+q^{15})(11)_q\big{)}(2)_q
\sqrt{(10)_q(9)_q(4)_q(3)_q}
c_{-3}^2
=0 .\notag
\end{align}

Next we see the gap at $7$. 
The operators 
$\sqrt{(12)_q(11)_q(2)_q}\,\Psi_7(\eta,\eta)_{-7}$ and 
$\sqrt{(12)_q(11)_q(2)_q}\,f\cdot\Psi_7(\eta,\eta)_{-7}$ 
are as follows.
\begin{align*}
\sqrt{(12)_q(11)_q(2)_q}&\,\Psi_7(\eta,\eta)_{-7}\\
&\ =
q^{-20}\sqrt{(12)_q(11)_q(2)_q}\eta_{-1}\eta_{-6}
-q^{-12}\sqrt{(11)_q(8)_q(5)_q(2)_q}\eta_{-2}\eta_{-5}\\
&\hspace{7mm}
+q^{-4}\sqrt{(9)_q(8)_q(5)_q(4)_q}\eta_{-3}\eta_{-4}
-q^{4}\sqrt{(9)_q(8)_q(5)_q(4)_q}\eta_{-4}\eta_{-3}\\
&\hspace{7mm}
+q^{12}\sqrt{(11)_q(8)_q(5)_q(2)_q}\eta_{-5}\eta_{-2}
-q^{20}\sqrt{(12)_q(11)_q(2)_q}\eta_{-6}\eta_{-1},
\end{align*}
\begin{align*}
\sqrt{(12)_q(11)_q(2)_q}\,f\cdot\Psi_7(\eta,\eta)_{-7}
=&\hspace{0.3mm}
q^{-14}\sqrt{(12)_q(11)_q(7)_q(6)_q(2)_q}\eta_{0}\eta_{-6}\\
&
+
\big{(}q^{-21}(12)_q-q^{-7}(8)_q(5)_q\big{)}
\sqrt{(11)_q(2)_q}\eta_{-1}\eta_{-5}\\
&
+
\big{(}(9)_q(4)_q-q^{-14}(11)_q(2)_q\big{)}\sqrt{(8)_q(5)_q}
\eta_{-2}\eta_{-4}\\
&
+
(q^{-7}-q^7)\sqrt{(10)_q(9)_q(8)_q(5)_q(4)_q(3)_q}
\eta_{-3}\eta_{-3} \\
&
+
\big{(}(9)_q(4)_q-q^{-14}(11)_q(2)_q\big{)}\sqrt{(8)_q(5)_q}
\eta_{-4}\eta_{-2}\\
&
+
\big{(}q^{7}(8)_q(5)_q-q^{21}(12)_q\big{)}
\sqrt{(11)_q(2)_q}
\eta_{-5}\eta_{-1}\\
&
-q^{14}\sqrt{(12)_q(11)_q(7)_q(6)_q(2)_q}\eta_{-6}\eta_{0}.
\end{align*}
From the second one we get the following equation
\begin{align}
(q^{-14}-q^{14})&\sqrt{(12)_q(11)_q(7)_q(6)_q(2)_q}c_0c_{-6}\label{a504}\\
&\quad
+
(q^{-7}-q^7)\sqrt{(10)_q(9)_q(8)_q(5)_q(4)_q(3)_q}
c_{-3}^2
=0\, . \notag
\end{align}
On equations (\ref{a503}) and (\ref{a504}) 
we can easily show that the determinant of 
the following matrix is not $0$ for $0<q<1$
\[
\begin{pmatrix}
(q^{-6}+q^6)& 
\big{(}(10)_q(3)_q-(q^{-15}+q^{15})(11)_q\big{)}(2)_q\\
(q^{-14}-q^{14})& (q^{-7}-q^7)(8)_q(5)_q\\
\end{pmatrix}\,.
\]
Hence we obtain $c_0c_{-6}=0$ and $c_{-3}=0$. 

Next we see the gap at $4$. 
We use an eigenvector $\Psi_4(\eta,\eta)$
\begin{align*}
&\sqrt{(12)_q(11)_q(10)_q(9)_q(4)_q(3)_q(2)_q}\,\Psi_4(\eta,\eta)_{-4}\\
&\quad\hspace{3cm}
=
q^{-20}\sqrt{(12)_q(11)_q(10)_q(9)_q(4)_q(3)_q(2)_q}\eta_{2}\eta_{-6}\\
&\hspace{38mm}
-q^{-15}\sqrt{(11)_q(10)_q(9)_q(8)_q(5)_q(4)_q(3)_q(2)_q}\eta_{1}\eta_{-5}\\
&\hspace{38mm}
+q^{-10}\sqrt{(10)_q(9)_q(8)_q(7)_q(6)_q(5)_q(4)_q(3)_q}\eta_{0}\eta_{-4}\\
&\hspace{38mm}
-q^{-5}(7)_q(6)_q\sqrt{(9)_q(8)_q(5)_q(4)_q}\eta_{-1}\eta_{-3}\\
&\hspace{38mm}
+(8)_q(7)_q(6)_q(5)_q\eta_{-2}^2\\
&\hspace{38mm}
-q^{5}(7)_q(6)_q\sqrt{(9)_q(8)_q(5)_q(4)_q}\eta_{-3}\eta_{-1}\\
&\hspace{38mm}
+q^{10}\sqrt{(10)_q(9)_q(8)_q(7)_q(6)_q(5)_q(4)_q(3)_q}\eta_{-4}\eta_{0}\\
&\hspace{38mm}
-q^{15}\sqrt{(11)_q(10)_q(9)_q(8)_q(5)_q(4)_q(3)_q(2)_q}\eta_{-5}\eta_{1}\\
&\hspace{38mm}
+q^{-20}\sqrt{(12)_q(11)_q(10)_q(9)_q(4)_q(3)_q(2)_q}\eta_{-6}\eta_{2}\, .
\end{align*}
We know $c_2c_{-6}=0$ because of $c_{-2}c_{-6}=0$. 
Hence the above equality derives $c_{-2}=0$. 

Next we see the gap at $2$. 
The operator $\Psi_2(\eta,\eta)_{-2}$ is
\begin{align*}
\sqrt{(12)_q(11)_q(2)_q}&\,\Psi_2(\eta,\eta)_{-2}\\
&\ 
=
q^{-15}\sqrt{(12)_q(11)_q(2)_q}\eta_{4}\eta_{-6}
-q^{-12}\sqrt{(11)_q(10)_q(3)_q(2)_q}\eta_{3}\eta_{-5}\\
&\quad
+q^{-9}\sqrt{(10)_q(9)_q(4)_q(3)_q}\eta_{2}\eta_{-4}
-q^{-6}\sqrt{(9)_q(8)_q(5)_q(4)_q}\eta_{1}\eta_{-3}\\
&\quad
+q^{-3}\sqrt{(8)_q(7)_q(6)_q(5)_q}\eta_{0}\eta_{-2}
-(7)_q(6)_q\eta_{-1}^2\\
&\quad
+q^{3}\sqrt{(8)_q(7)_q(6)_q(5)_q}\eta_{-2}\eta_{0}
-q^{6}\sqrt{(9)_q(8)_q(5)_q(4)_q}\eta_{-3}\eta_{1}\\
&\quad
+q^{9}\sqrt{(10)_q(9)_q(4)_q(3)_q}\eta_{-4}\eta_{2}
-q^{12}\sqrt{(11)_q(10)_q(3)_q(2)_q}\eta_{-5}\eta_{3}\\
&\quad
+q^{15}\sqrt{(12)_q(11)_q(2)_q}\eta_{-6}\eta_{4}\, .
\end{align*}
This shows $c_{-1}=0$. 

Finally we see the gap at $1$. 
The operators $\Psi_1(\eta,\eta)_{-1}$ and $f\cdot\Psi_1(\eta,\eta)_{-1}$ 
are
\begin{align*}
\sqrt{(12)_q}\,\Psi_1(\eta,\eta)_{-1}
\!
=&\hspace{0.3mm}
q^{-11}\sqrt{\!(12)_q}\eta_{5}\eta_{-6}
\!-\!q^{-9}\sqrt{\!(11)_q(2)_q}\eta_{4}\eta_{-5}
\!+\!q^{-7}\sqrt{\!(10)_q(3)_q}\eta_{3}\eta_{-4}\\
&\hspace{-0.5mm}
\!-\!q^{-5}\sqrt{\!(9)_q(4)_q}\eta_{2}\eta_{-3}
\!+\!q^{-3}\sqrt{\!(8)_q(5)_q}\eta_{1}\eta_{-2}
\!-\!q^{-1}\!\sqrt{\!(7)_q(6)_q}\eta_{0}\eta_{-1}\\
&
+q\sqrt{(6)_q(7)_q}\eta_{-1}\eta_{0}
-q^3\sqrt{(5)_q(8)_q}\eta_{-2}\eta_{1}
+q^5\sqrt{(4)_q(9)_q}\eta_{-3}\eta_{2}\\
&
+q^{7}\sqrt{(10)_q(3)_q}\eta_{-4}\eta_{3}
-q^{9}\sqrt{(11)_q(2)_q}\eta_{-5}\eta_{4}
-q^{11}\sqrt{(12)_q}\eta_{-6}\eta_{5}\, ,
\end{align*}
\begin{align*}
\sqrt{(12)_q}&\,f\cdot\Psi_1(\eta,\eta)_{-1}\\
=
&\,
q^{-5}(12)_q\eta_{6}\eta_{-6}
+\big{(}q^{-6}(12)_q-q^{-4}(11)_q(2)_q\big{)}\eta_{5}\eta_{-5}\\
&\hspace{0.5mm}+
\big{(}q^{-3}(10)_q(3)_q-q^{-5}(11)_q(2)_q\big{)}\eta_{4}\eta_{-4}
+\big{(}q^{-4}(10)_q(3)_q-q^{-2}(9)_q(4)_q\big{)}\eta_{3}\eta_{-3}\\
&\hspace{0.5mm}+
\big{(}q^{-1}(8)_q(5)_q-q^{-3}(9)_q(4)_q\big{)}\eta_{2}\eta_{-2}
+\big{(}q^{-2}(8)_q(5)_q-(7)_q(6)_q\big{)}\eta_{1}\eta_{-1}\\
&\hspace{0.5mm}+
(q-q^{-1})(7)_q(6)_q\eta_{0}\eta_{0}
+\big{(}(7)_q(6)_q-q^{2}(8)_q(5)_q\big{)}\eta_{-1}\eta_{1}\\
&\hspace{0.5mm}+
\big{(}q^{3}(9)_q(4)_q-q (8)_q(5)_q\big{)}\eta_{-2}\eta_{2}
+\big{(}q^{2}(9)_q(4)_q-q^{4}(10)_q(3)_q\big{)}\eta_{-3}\eta_{3}\\
&\hspace{0.5mm}+
\big{(}q^{5}(11)_q(2)_q-q^{3}(10)_q(3)_q\big{)}\eta_{-4}\eta_{4}
+\big{(}q^{4}(11)_q(2)_q-q^{6}(12)_q\big{)}\eta_{-5}\eta_{5}\\
&\hspace{0.5mm}-q^{5}(12)_q\eta_{-6}\eta_{6}\, .
\end{align*}
Since we have $c_s=0$ for $s=5,4,3,2$ and $1$, 
the above equality shows
\begin{equation*}
(q^{-5}-q^5)(12)_q c_6c_{-6}
+(q-q^{-1})(7)_q(6)_q c_0^2
=0\,.
\end{equation*}
We also know $c_0c_6=0$ because of $c_0c_{-6}=0$. 
Hence we obtain $c_{-6}=c_0=0$. 

(IV)$D_\infty^*$ case. 
We conclude that this case actually occurs and a right 
coideal of this type is unique up to conjugation by $\beta_z^L$. 

\begin{lem}\label{pi2-eigenvector}
Consider the quantum sphere $C(S_{q,0}^2)$. 
Then a vector 
\[\Big{(}q\sqrt{(3)_q!},0,-\sqrt{(4)_q},0,q^{-1}\sqrt{(3)_q!}\Big{)}
\,w(\pi_2)\] 
is a $\pi_2$-eigenvector of $C(S_{q,0}^2)$.
\end{lem}

\begin{proof}
Recall a highest weight vector of weight $1$ $(\xi_0^1)_{-1}$ in 
$C(S_{q,0}^2)$, 
$(\xi_0^1)_{-1}
=
q^{\frac{1}{2}}\sqrt{(2)_q}^{-1}x^2
-
q^{-\frac{1}{2}}\sqrt{(2)_q}^{-1}v^2\, .$ 
Then we can easily compute a highest weight vector of weight 
$2$ in $C(S_{q,0}^2)$ as follows, 
$(\xi_0^1)_{-1}^2=
q (2)_q^{-1}w(\pi_2)_{-2 -2} 
-
(q^{-2}+q^2)\sqrt{(4)_q!}^{-1}w(\pi_2)_{0,-2}
+
q^{-1}(2)_q^{-1} w(\pi_2)_{2,-2}\, .$ 
Multiply $(2)_q\sqrt{(3)_q!}$ to both sides we get 
the desired highest weight vector. 
\end{proof}

Let $\xi=(\xi_{-1},\xi_0,\xi_1)$ be the canonical $\pi_1$-eigenvector 
of $C(S_{q,0}^2)$. 
By its definition, $\{\xi_r\}_{r\in I_1}$ satisfies the following conditions
\begin{equation}\label{Dinf01}
\xi_{-1}^*=-q\,\xi_1,\ \xi_0^*=\xi_0,\ 
\xi_{\pm1}\xi_0=q^{\pm2}\xi_0\xi_{\pm1}, 
\end{equation}
\begin{align}
\xi_1\xi_{-1}
=&\,
(q+q^{-1})^{-1}(q^2\xi_0^2-1),\label{Dinf02}\\
\xi_{-1}\xi_{1}
=&\,
(q+q^{-1})^{-1}(q^{-2}\xi_0^2-1).\label{Dinf03}
\end{align}
An embedding of $C(S_{q,0}^2)$ into $\csuq$ is given by
\[
\xi=(\sqrt{1+q^{-2}}^{-1},0,-\sqrt{1+q^2}^{-1})\,w(\pi_1)\,,
\]
or more precisely 
\[
\xi_{-1}=\sqrt{1+q^{-2}}^{-1} x^2 -\sqrt{1+q^{2}}^{-1} v^2\, ,\quad
\xi_0=q xu -vy\, ,\quad
\]
\[
\xi_{-1}=\sqrt{1+q^{-2}}^{-1} u^2 -\sqrt{1+q^{2}}^{-1} y^2\, .
\]

Consider the smooth part of $C(S_{q,0}^2)$, that is, 
the $*$-algebra generated by $\{\xi_r\}_{r\in I_1}$. 
Denote it by $\mathcal{B}$. 
Let $\mathcal{B}_{\mbox{e}}$ be a $\suq$-invariant 
$*$-subalgebra of $\mathcal{B}$ 
generated by $\{\xi_{r}\xi_{s}\}_{r,s\in I_1}$. 
In \cite[Proposition 2.9.]{D-Koornwinder}, it is shown that 
$\{\xi_0^m \xi_{-1}^n\}_{m,n\in \Z_{\geq0}}$ and 
$\{\xi_0^m \xi_{1}^n\}_{m,n\in \Z_{\geq0}}$ are 
basis for a vector space $\mathcal{B}$. 
Since the $*$-subalgebra $\mathcal{B}_{\mbox{e}}$ is 
generated by words of even length 
(this is because of the previous equalities), 
we see 
$\xi_r$ is not in $\mathcal{B}_{\mbox{e}}$ for $r\in I_1$. 
In particular, we have $\mathcal{B}_{\mbox{e}}\subsetneq \mathcal{B}$. 
Let $B_{\mbox{e}}$ be the norm closure of $\mathcal{B}_{\mbox{e}}$. 
By $\suq$-invariance of $\mathcal{B}_{\mbox{e}}$ 
we see that $B_{\mbox{e}}$ is also $\suq$-invariant and 
$B_{\mbox{e}}\subsetneq C(S_{q,0}^2)$. 
Next we prove that $\pi_3$-spectral subspace of $B_{\mbox{e}}$ is $0$. 
Since $\xi_{-1}^3\neq0$ is 
the highest weight vector in $C(S_{q,0}^2)_{\pi_3}$, 
it is not contained in $\mathcal{B}_{\mbox{e}}$. 
Again $\suq$-invariance yields that $\xi_{-1}^3$ is not contained in 
$B_{\mbox{e}}$. 
In a similar way we can prove that all the odd spin spectral subspaces 
of $B_{\mbox{e}}$ are $0$. 
Of course, even spin spectral subspaces 
of $B_{\mbox{e}}$ are not $0$ because of $\xi_{-1}^{2n}\neq0$ for all 
$n\in \Z_{\geq0}$. 
The spectral pattern $\oplus_{k\in \Z_{\geq0}}\pi_{2k}$ 
follows from the spectral pattern of $C(S_{q,0}^2)$ 
$\oplus_{k\in \Z_{\geq0}}\pi_{k}$. 
Therefore, we have shown that $B_{\mbox{e}}$ is of type $D_{\infty}^*$. 
Next lemma says that a right coideal of $D_{\infty}^*$-type 
is unique up to conjugation of the automorphism $\beta_{z}^{L}$. 

\begin{lem}\label{uniqueness of infinite dihedral}
Let $A\subset \csuq$ be a right coideal of type $D_{\infty}^*$. 
Then 
there exists $z\in\C$ such that 
$A$ is $\beta_{z}^L\big{(}B_{\mbox{e}}\big{)}$. 
\end{lem}

\begin{proof}
By assumption $A$ has the spectral pattern 
$\pi_0\oplus \pi_2 \oplus \pi_4\oplus \cdots$. 
In the following discussion 
we do not need to use the gap at $\pi_3$. 
Let $\eta=(\eta_s)_{s\in I_2}=\sum_{s\in I_2}c_s \bw_s^2$ be 
a self conjugate $\pi_2$-eigenvector of $A$. 
We consider zero vectors $\Psi_1(\eta,\eta)_{-1}$ and 
$f\cdot \Psi_1(\eta,\eta)_{-1}$. 
They become as follows.
\begin{align*}
\Psi_1(\eta,\eta)_{-1}
=&\,
q^{-3}\,\eta_1 \eta_{-2} 
-q^{-1}\sqrt{\frac{(3)_q(2)_q}{(4)_q}}\,\eta_0\eta_{-1}
+q \sqrt{\frac{(3)_q(2)_q}{(4)_q}}\,\eta_{-1}\eta_{0}
-q^3\,\eta_{-2} \eta_{1}\, ,\\
f\cdot\Psi_1(\eta,\eta)_{-1}
=&\,q^{-1}\sqrt{(4)_q}\,\eta_2\eta_{-2}
+\sqrt{(4)_q}^{-1}(q^{-2}(4)_q-(3)_q!)\,\eta_{1}\eta_{-1}\\
&\,
+\sqrt{(4)_q}^{-1}(q-q^{-1})(3)_q!\,\eta_0^2 
+\sqrt{(4)_q}^{-1}((3)_q!-q^{2}(4)_q)\,\eta_{-1}\eta_{1}\\
&\,
-q\sqrt{(4)_q}\,\eta_{-2}\eta_{2}\, .
\end{align*}
Apply the map $\pi_\T$ to the above equalities and 
we obtain
\begin{align*}
&\,(q^{-3}-q^3)c_1 c_{-2}+(q-q^{-1})\sqrt{\frac{(3)_q!}{(4)_q}}c_{-1}c_0
=0\ ,
\\
&\,
(q^{-1}-q)\sqrt{(4)_q}c_2c_{-2}+(q^{-2}-q^{2})\sqrt{(4)_q}c_1c_{-1} \, ,
\\
&\qquad+(q-q^{-1})\sqrt{(4)_q}^{-1}(3)_q!c_0^2=0\, .
\end{align*}
For simplicity, we may assume that $c_1$ is a real number by 
applying a conjugation $\beta_z^L$. 
Then we have $c_{-1}=-q c_1$ by the self-conjugacy of $\eta$. 
Assume $c_1$ is not $0$. 
Then from the first equality we get 
$c_{-2}=-q\sqrt{\frac{(2)_q}{(4)_q(3)_q}}c_0$. 
If we put this one to the second equality, we get
\[
\frac{(2)_q-(3)_q^2(2)_q}{(3)_q}c_0^2
+(2)_q (4)_qc_1c_{-1}=0\, .
\]
Since we know $(3)_q>1$ and $c_1c_{-1}=-q c_1^2 <0$, 
the left hand side is strictly negative unless $c_0=c_1=0$. 
This also shows $c_2=0$ and this is a contradiction. 
Hence we get $c_1=0$. 
Then from the second equality we get 
$|c_2|=q^{-1}\sqrt{\frac{(3)_q!}{(4)_q}}|c_0|$. 
We may assume $c_0=-\sqrt{(4)_q}$ by scalar multiplication 
and $c_{2}=q^{-1}\sqrt{(3)_q!}$ by a conjugation of $\beta_z^L$.

\end{proof}

From now on, the right coideal of type $D_\infty^*$ 
generated by 
$q\sqrt{(3)_q!}w(\pi_2)_{-2,s}-\sqrt{(4)_q}w(\pi_2)_{0,s}
+q^{-1}\sqrt{(3)_q!}w(\pi_2)_{2,s}$ for $s\in I_2$ is 
denoted by $A_{D_\infty^*}$. 
We remark that $\bw_0^2$ does not generates a right coideal 
of type $D_\infty^*$. 
In fact, it generates the canonical homogeneous sphere 
$C(\T\setminus SU_q(2))$. 
This is essentially due to the effect of $q\neq1$ 
(see also Lemma \ref{product of pi2 negative}). 
The $C^*$-algebra $A_{D_\infty^*}$ is actually 
the Toeplitz algebra as we see below. 
By \cite[Lemma 3.2]{Masuda2} there exists the matrix units 
$\{e_{i,j}^{+}\}_{i,j\in \Z_{\geq0}}$ and 
$\{e_{i,j}^{-}\}_{i,j\in \Z_{\geq0}}$ in $C(S_{q,0}^2)$ 
which satisfy the following equalities
\begin{equation}\label{Dinf04}
\xi_{-1}e_{k,k}^{\pm}
=
\pm\gamma_k e_{k+1,k}^{\pm},\ 
\xi_{1}e_{k,k}^{\pm}
=
\mp q^{-1}\gamma_{k-1} e_{k-1,k}^{\pm},
\end{equation}
\begin{equation}
\xi_0
=
\sum_{k=0}^{\infty}-q^{2k+1}e_{k,k}^{-}+
\sum_{k=0}^{\infty}q^{2k+1}e_{k,k}^{+}
\end{equation}
for all $k\in \Z_{\geq0}$, where 
$\gamma_k=q^{\frac{1}{2}}
(q+q^{-1})^{-\frac{1}{2}}(1-q^{4k+4})^{\frac{1}{2}}$. 
Note the equality 
$1=\sum_{k=0}^\infty e_{k,k}^{-}+\sum_{k=0}^\infty e_{k,k}^{+}$ 
in $C(\suq)$ where the summation converges in the strong operator 
topology in $\mathbb{B}(L^2(\suq))$. 
In order to verify it, it suffices to show 
$1=\sum_{k=0}^\infty h(e_{k,k}^{-})+\sum_{k=0}^\infty h(e_{k,k}^{+})$, 
because the Haar state $h$ is normal and faithful. 
This equality holds by the following lemma. 

\begin{lem}
We have $h(e_{k,k}^{\pm})=2^{-1}q^{2k}(1-q^2)$ for all $k\in \Z_{\geq0}$. 
\end{lem}

\begin{proof}
Let $X_{\pi_1}=\C\xi$ be the $\pi_1$-eigenvector space of 
$C(S_{q,0}^2)$. 
Recall the conjugation map $T$ and the one-parameter unitary group 
$U_t$ on $X_{\pi_1}$ which are defined by 
$(T\xi)_r=(-q)^{-r}\xi_{-r}^*$ and 
$(U_t\xi)_r=(F_{\pi_1})_{r,r}^{-it}\sigma_t^h(\xi_r)$ for all 
$r\in I_{\pi_1}$ and $t\in \R$, 
where $\sigma_t^h$ is the modular automorphism group with respect to 
the invariant state $h$ on $C(S_{q,0}^2)$. 
As is proved in Lemma \ref{T*T=U}, we have $IT^*TI^*=U_i$ where 
the unitary map $I$ is defined by $(I\xi)_r=\xi_{-r}$ for $r\in I_{\pi_1}$. 
It shows that the spectrum of $U_i$ is inverse-closed and hence 
$U_i=\id$. 
It yields the formula $\sigma_t^h(\xi_r)=q^{-2rit}\xi_r$. 
We have the equality 
$h(\xi_1 e_{k,k}^{\pm}\xi_{-1})=-q^{-1}\gamma_{k-1}^2h(e_{k-1,k-1}^{\pm})$ 
by (\ref{Dinf04}). 
The left hand side is equal to 
\begin{align*}
h(\xi_1 e_{k,k}^{\pm}\xi_{-1})
=&\,
h(\sigma_i^h(\xi_{-1})\xi_1 e_{k,k}^{\pm})\\
=&\,
q^{-2}h(\xi_{-1}\xi_1 e_{k,k}^{\pm})\\
=&\,
q^{-2}(q+q^{-1})^{-1}(q^{4k}-1)h(e_{k,k}^{\pm}), 
\end{align*}
where we used (\ref{Dinf02}) for the third equality. 
Hence we have 
$h(e_{k,k}^{\pm})=q^{2k}h(e_{0,0}^{\pm})$ for $k\geq 1$. 
Since the equality $h(\xi_0)=0$ holds, we get 
$h(e_{0,0}^{+})=h(e_{0,0}^{-})$. 
Using Lemma \ref{innerproduct}, we see $h(\xi_0^2)=(3)^{-1}$ and obtain 
$h(e_{0,0}^{\pm})=2^{-1}(1-q^2)$. 
\end{proof}
Now we consider the invariant subalgebra 
$A_{D_\infty^*}\subset C(S_{q,0}^2)$. 
Let us define the matrix units $e_{i,j}=e_{i,j}^{-}+e_{i,j}^{+}$ for 
$i,j\in \Z_{\geq0}$. 
The $C^*$-algebra $A_{D_\infty^*}$ is generated by them and 
the unilateral shift. 
It also concludes that the right coideal von Neumann algebra 
${A_{D_\infty^*}}''$ is isomorphic to $\mathbb{B}(\ell^2(\Z_{\geq0}))$. 
Especially the quantum group 
$\suq$ can act ergodically on every separable type I factor. 

(V) $D_m^*\, (m\geq 2)$ case. 
We will derive the non-existence of this case. 
We study a $C^*$-algebra $A$ whose spectral pattern is 
$\oplus_{k\in\Z_{\geq0}}\big{(}\frac{1+(-1)^k}{2}
+\big{[}\frac{k}{m}\big{]}\big{)}\pi_k$. 
Making use of a gap at $\pi_1$ as Lemma \ref{uniqueness of infinite dihedral}, 
we can similarly prove that 
elements of 
one $\pi_2$-eigenvector generates a right coideal of type 
$D_\infty^*$. 
(the $\pi_2$-eigenvector space is two-dimensional when $m$ is equal to $2$.) 
So we may assume there is a $\pi_2$-eigenvector $\xi_0^2$ in $A$. 
Recall that $\xi_0^2$ is defined by 
$\xi_0^2=\Big{(}q\sqrt{(3)_q(2)_q},0,
-\sqrt{(4)_q},0,q^{-1}\sqrt{(3)_q(2)_q}\Big{)}
\,w(\pi_2)$. 
Note that the index $0$ of $\xi_0^2$ means $\lambda=0$ of $\CSq$. 

\begin{lem}\label{product of pi2}
Let $n$ be a half integer with $n \geq 2$. 
Then we have the following equalities for all $t\in I_n$
\begin{align*}
&
\Psi_{n-2}(\bw_{-2}^2,\bw_t^n)
=q^{-2t+6}
\sqrt{
\frac{(n\!+\!t)_q(n\!+\!t\!-\!1)_q(n\!+\!t\!-\!2)_q(n\!+\!t\!-\!3)_q}
{(2n)_q(2n\!-\!1)_q(2n\!-\!2)_q(2n\!-\!3)_q}
}\,\bw_{t-2}^{n-2}
\, ,\\
&
\Psi_{n-2}(\bw_{0}^2,\bw_t^n)
=
q^{-2n-2t+4}
\sqrt{\frac{(4)_q(3)_q}{(2)_q}}
\sqrt{
\frac{(n\!+\!t)_q(n\!+\!t\!-\!1)_q(n\!-\!t)_q(n\!-\!t\!-\!1)_q}
{(2n)_q(2n\!-\!1)_q(2n\!-\!2)_q(2n\!-\!3)_q}
}\,\bw_{t}^{n-2}
,\\
&
\Psi_{n-2}(\bw_{2}^2,\bw_t^n)
=q^{-4n-2t+2}
\sqrt{
\frac{(n\!-\!t)_q(n\!-\!t\!-\!1)_q(n\!-\!t\!-\!2)_q(n\!-\!t\!-\!3)_q}
{(2n)_q(2n\!-\!1)_q(2n\!-\!2)_q(2n\!-\!3)_q}
}\,\bw_{t+2}^{n-2}
\, .
\end{align*}

\end{lem}

\begin{proof}
For the first equality we use 
$w(\pi_2)_{-2,2-r}={4 \brack r}_{q^2}^{\frac{1}{2}}x^r u^{4-r}$ 
for $r\in I_{2}$. 
From this we obtain
\[
\Psi_{n-2}(\bw_{-2}^2,\bw_{t}^n)_{-(n-2)}=\sum_{r=0}^4
(C_{n-2}^{2,n})_{r}{4 \brack r}_{q^2}^{\frac{1}{2}}x^r u^{4-r}
w(\pi_n)_{t,-n+r}
\]
for all $t\in I_n$ where the coefficients $(C_{n-2}^{2,n})_{r}$ are 
given 
as follows.
\begin{center}
$(C_{n-2}^{2,n})_{0}=q^{-2(n-1)}\, ,$ 
$(C_{n-2}^{2,n})_{1}=-q^{-(n-1)}\sqrt{\frac{(4)_q}{(2n)_q}}\, ,$ \\
\vspace{0.3cm}
$(C_{n-2}^{2,n})_{2}=\sqrt{\frac{(4)_q(3)_q}{(2n)_q(2n-1)_q}}\, ,$
$(C_{n-2}^{2,n})_{3}=-q^{n-1}
\sqrt{\frac{(4)_q(3)_q(2)_q}{(2n)_q(2n-1)_q(2n-2)_q}}\, ,$\\
\vspace{0.3cm}
$(C_{n-2}^{2,n})_{4}=q^{2(n-1)}
\sqrt{\frac{(4)_q(3)_q(2)_q}{(2n)_q(2n-1)_q(2n-2)_q(2n-3)_q}}\, .$
\end{center}
This must be a scalar multiple of $x^{n-t}v^{n+t-4}$. 
So we may compute only the $r=4$ term
\[
(C_{n-2}^{2,n})_{4}\,x^4 
w(\pi_n)_{t,-n+4}
.\]
If $w(\pi_n)_{t,-n+4}$ contains the word $u$ such that 
$t\lneq -n+4$, we can disregard the effect of this term. 
In order to simplify calculations we use the symbol $\sim$ which 
ignores terms with a power of $u$. 
Assume $t\geq -n+4$ then we have
\begin{center}
$w(\pi_n)_{t,-n+4}\sim q^{4(-n-t+4)}{n+t \brack 4}_{q^2}^{\frac{1}{2}}
{2n-4 \brack n-t}_{q^2}^{\frac{1}{2}} x^{n-t-4}v^{n+t-4}$ if $t-n+4\leq 0$\\  
$w(\pi_n)_{t,-n+4}\sim q^{(n+t-4)(t-n)}{n+t \brack 4}_{q^2}^{\frac{1}{2}}
{2n-4 \brack n-t}_{q^2}^{\frac{1}{2}}v^{n+t-4}y^{-n+t+4}$ if $t-n+4\geq 0$. 
\end{center}
In both cases 
we have the same result 
\begin{align*}
x^4 w(\pi_n)_{t,-n+4}\sim&\,
q^{4(-n-t+4)}{n+t \brack 4}_{q^2}^{\frac{1}{2}}
{2n-4 \brack n-t}_{q^2}^{\frac{1}{2}} x^{n-t}v^{n+t-4}\\
=&\,
q^{4(-n-t+4)}{n+t \brack 4}_{q^2}^{\frac{1}{2}}
w(\pi_{n-2})_{t-2,-(n-2)}
\end{align*} 
by simple calculation. 
Hence we obtain the first equality
\begin{align*}
\Psi_{n-2}(\bw_{-2}^2,\bw_{t}^n)_{-(n-2)}
=&\,
q^{2(n-1)}q^{4(-n-t+4)}
\sqrt{\frac{(4)_q !}{(2n)_q\cdots(2n-3)_q}}\\
&\qquad
\cdot{n+t \brack 4}_{q^2}^{\frac{1}{2}}\,w(\pi_{n-2})_{t-2,-(n-2)}\\
=&\,
q^{-2t+6}
\sqrt{\frac{(n+t)_q\cdots(n+t-3)_q}{(2n)_q\cdots(2n-3)_q}}
\,w(\pi_{n-2})_{t-2,-(n-2)}\, .
\end{align*}
Next we show the second equality. 
Then the matrix element $w(\pi_2)_{0,2-r}$ is as follows.
\begin{center}
$w(\pi_2)_{0,2-r}\sim q^{(4-r)(2-r)}{2 \brack r-2}_{q^2}^{\frac{1}{2}}
{r \brack 2}_{q^2}^{\frac{1}{2}} x^{r-2}v^{r-2}$ \quad if $r\geq 2\, ,$\\
\vspace{0.3cm}
$w(\pi)_{0,2-r}\sim 0$ \quad if $r\leq 1$. 
\end{center}
Hence we may only compute the coefficient of $x^{n-t-2}v^{n+t-2}$ of 
the next one
\[
\sum_{r\geq2}(C_{n-2}^{2,n})_{r}
q^{(4-r)(2-r)}{2 \brack r-2}_{q^2}^{\frac{1}{2}}
{r \brack 2}_{q^2}^{\frac{1}{2}} x^{r-2}v^{r-2}
w(\pi_n)_{t,-n+r}
.\]
If $t\geq -n+4$, then we have 
\begin{align*}
w(\pi_2)_{0,2-r}w(\pi_n)_{t,-n+r}
\sim& \,
q^{r^2-(2t+4)r-2n+2t+8}
{n+t \brack r}_{q^2}^{\frac{1}{2}}
{2n-r \brack n-t}_{q^2}^{\frac{1}{2}}\\
&\,{2 \brack r-2}_{q^2}^{\frac{1}{2}}
{r \brack 2}_{q^2}^{\frac{1}{2}}
{2n-4 \brack n+t-2}_{q^2}^{-\frac{1}{2}}
\, 
w(\pi_{n-2})_{t,-(n-2)}\, .
\end{align*}
When $t< -n+4$, then 
$w(\pi_n)_{t,-n+r}$ contains the word $u$ if and only if 
$2n-r < n-t$. 
In this case the binomial ${2n-r \brack n-t}_{q^2}^{\frac{1}{2}}$ is 
$0$ by definition. 
From this observation we may continue the computation 
under $t\geq -n+4$ in order to get the result about general $t$. 
Then we can carry out the following calculation 
\begin{align*}
&\sum_{r=2}^{4}
(C_{n-2}^{2,n})_{r}
q^{r^2-(2t+4)r-2n+2t+8}
{n+t \brack r}_{q^2}^{\frac{1}{2}}
{2n-r \brack n-t}_{q^2}^{\frac{1}{2}}
{2 \brack r-2}_{q^2}^{\frac{1}{2}}
{r \brack 2}_{q^2}^{\frac{1}{2}}
{2n-4 \brack n+t-2}_{q^2}^{-\frac{1}{2}}\\
=&\,
(C_{n-2}^{2,n})_{2} q^{-2n-2t+4}
{n+t \brack 2}_{q^2}^{\frac{1}{2}}
{2n-2 \brack n-t}_{q^2}^{\frac{1}{2}}
{2n-4 \brack n+t-2}_{q^2}^{-\frac{1}{2}}
\\
&\quad+
(C_{n-2}^{2,n})_{3}
q^{-2n-4t+5}
{n+t \brack 3}_{q^2}^{\frac{1}{2}}
{2n-3 \brack n-t}_{q^2}^{\frac{1}{2}}
{2 \brack 1}_{q^2}^{\frac{1}{2}}
{3 \brack 2}_{q^2}^{\frac{1}{2}}
{2n-4 \brack n+t-2}_{q^2}^{-\frac{1}{2}}\\
&\quad+
(C_{n-2}^{2,n})_{4}
q^{-2n-4t+5}
{n+t \brack 4}_{q^2}^{\frac{1}{2}}
{2n-4 \brack n-t}_{q^2}^{\frac{1}{2}}
{2 \brack 1}_{q^2}^{\frac{1}{2}}
{4 \brack 2}_{q^2}^{\frac{1}{2}}
{2n-4 \brack n+t-2}_{q^2}^{-\frac{1}{2}}\\
=&\,
{2n-4 \brack n+t-2}_{q^2}^{-\frac{1}{2}}
\sqrt{\frac{(4)_q(3)_q}{(2)_q}}
\sqrt{\frac{(2n-4)_q!}{(2n)_q\cdots(2n-3)_q}}
\sqrt{\frac{(n+t)_q!}{(n-t)_q!}}
\,q^{\frac{1}{2}n^2-\frac{1}{2}t^2-n-t+1}\\
&\quad\cdot
\Big{\{}
\frac{(2n-2)_q(2n-3)_q}{(n+t-2)_q!}q^{-n+t+1}
-\frac{(2n-3)_q(2)_q}{(n+t-3)_q!}
+\frac{q^{n-t-1}}{(n+t-4)_q!}
\Big{\}}\\
=&\,
q^{-2n-2t+4}
\sqrt{\frac{(4)_q(3)_q}{(2)_q}}
\sqrt{\frac{(n+t)_q(n+t-1)_q(n-t)_q(n-t-1)_q}{(2n)_q(2n-1)_q(2n-2)_q(2n-3)_q}}
\,.
\end{align*}
Hence we have proved the second equality. 

Finally we prove the third equality. 
Then the matrix element $w(\pi_2)_{2,2-r}$ is as follows.
\[
w(\pi_2)_{2,2-r}
=
{4 \brack r}_{q^2}^{\frac{1}{2}}v^r y^{4-r}
\]
for all $0\leq r \leq 4$. 
If $t\geq -n+4$, we have 
\begin{align*}
v^r y^{4-r} w(\pi_n)_{t,-n+r}
\sim
&\,
q^{r(r-2t-4)}{n+t \brack r}_{q^2}^{\frac{1}{2}}
{2n-r \brack n-t}_{q^2}^{\frac{1}{2}}
x^{n-t-4}v^{n+t}\\
=&\,
q^{r(r-2t-4)}{n+t \brack r}_{q^2}^{\frac{1}{2}}
{2n-r \brack n-t}_{q^2}^{\frac{1}{2}}
{2n-4 \brack t+2}_{q^2}^{-\frac{1}{2}}
\,
w(\pi_{n-2})_{t+2,-(n-2)}\,.
\end{align*}

For $t< -n+4$, $w(\pi_n)_{t,-n+r}$ contains the word $u$ 
if and only if $n+t < r$. 
Then the binomial ${n+t \brack r}_{q^2}$ is equal to $0$. 
As in the proof of the second equality, we may treat only $t\geq -n+4$. 
We compute the coefficient of the desired equality as follows. 
\begin{align*}
&\,\sum_{r=0}^{4}
(C_{n-2}^{2,n})_{r}
{4 \brack r}_{q^2}^{\frac{1}{2}}
q^{r(r-2t-4)}{n+t \brack r}_{q^2}^{\frac{1}{2}}
{2n-r \brack n-t}_{q^2}^{\frac{1}{2}}
{2n-4 \brack t+2}_{q^2}^{-\frac{1}{2}}\\
=&\,
q^{-2(n-1)}{2n \brack n-t}_{q^2}^{\frac{1}{2}}
{2n-4 \brack t+2}_{q^2}^{-\frac{1}{2}}\\
&\quad 
-
q^{-(n-1)}q^{-2t-3}\sqrt{\frac{(4)_q}{(2)_q}}
{4 \brack 1}_{q^2}^{\frac{1}{2}}
{n+t \brack 1}_{q^2}^{\frac{1}{2}}
{2n-1 \brack n-t}_{q^2}^{\frac{1}{2}}
{2n-4 \brack t+2}_{q^2}^{-\frac{1}{2}}\\
&\quad
+
q^{2(-2t-2)}\sqrt{\frac{(4)_q(3)_q}{(2n)_q(2n-1)_q}}
{4 \brack 2}_{q^2}^{\frac{1}{2}}
{n+t \brack 2}_{q^2}^{\frac{1}{2}}
{2n-2 \brack n-t}_{q^2}^{\frac{1}{2}}
{2n-4 \brack t+2}_{q^2}^{-\frac{1}{2}}\\
&\quad
-
q^{n-1}q^{3(-2t-1)}
\sqrt{\frac{(4)_q!}{(2n)_q\cdots(2n-2)_q}}
{4 \brack 3}_{q^2}^{\frac{1}{2}}
{n+t \brack 3}_{q^2}^{\frac{1}{2}}
{2n-3 \brack n-t}_{q^2}^{\frac{1}{2}}
{2n-4 \brack t+2}_{q^2}^{-\frac{1}{2}}\\
&\quad
+
q^{2(n-1)}q^{-8t}
\sqrt{\frac{(4)_q!}{(2n)_q\cdots(2n-3)_q}}
{n+t \brack 4}_{q^2}^{\frac{1}{2}}
{2n-4 \brack n-t}_{q^2}^{\frac{1}{2}}
{2n-4 \brack t+2}_{q^2}^{-\frac{1}{2}}
\end{align*}
\begin{align*}
=&\,
q^{\frac{1}{2}n^2-\frac{1}{2}t^2-2t-4}
\sqrt{\frac{(n+t)_q!}{(n-t)_q!}}(n+t)_q!^{-1}
\frac{\sqrt{(2n)_q!}}{(2n)_q\cdots(2n-3)_q}
{2n-4 \brack t+2}_{q^2}^{-\frac{1}{2}}\\
&\,\cdot
\Big{\{}
q^{-2n+2t+6}(2n)_q\cdots(2n-3)_q
-
q^{-n+t+3}(4)_q(2n-1)_q\cdots(2n-3)_q(n+t)_q\\
&\quad\, 
+
\frac{(4)_q(3)_q}{(2)_q}(2n-2)_q(2n-3)_q(n+t)_q(n+t-1)_q\\
&\quad\, 
-
q^{n-t-3}(4)_q(2n-3)_q(n+t)_q\cdots(n+t-2)_q\\
&\quad\, 
+
q^{2n-2t-6}(n+t)_q\cdots(n+t-3)_q
\Big{\}}\\
=&\,
q^{\frac{1}{2}n^2-\frac{1}{2}t^2-2t-4}
\sqrt{\frac{(n+t)_q!}{(n-t)_q!}}\frac{1}{(n+t)_q!}
\frac{\sqrt{(2n)_q!}}{(2n)_q\cdots(2n-3)_q}
{2n-4 \brack t+2}_{q^2}^{-\frac{1}{2}}\\
&\,\cdot
q^{-6n-2t+6}(n-t)_q\cdots(n-t-3)_q\\
=&\,
q^{-4n-2t+2}
\sqrt{\frac{(n-t)_q\cdots(n-t-3)_q}{(2n)_q\cdots(2n-3)_q}}. 
\end{align*}

\end{proof}

Let $n$ be the smallest odd integer which satisfies 
$\big{[}\frac{n}{m}\big{]}=1$. 
Then $\pi_n$ appears in the spectral pattern of $A$ once and 
$\pi_{n-2}$ does not. 
Therefore $\Psi_{n-2}: X_{2}\times X_{n}\longrightarrow X_{n-2}$ 
is a $0$-map. 
Let $\eta=\sum_{t\in I_n}d_t \bw_t^n$ be a self-conjugate 
$\pi_n$-eigenvector. 
From the self-conjugacy 
we obtain $d_{-t}=(-q)^{t}\overline{d_t}$ for all $t \in I_n$. 
We shall show the complex numbers $\{d_t\}_{t\in I_n}$ are all $0$. 
It derives non-existence of type $D_m^*$ for $m\geq2$. 

\begin{lem}\label{kernel of pi2-product}
The complex numbers $\{d_t\}_{t\in I_n}$ satisfy the following 
recurrence formula for $t\in I_n=\{-n,-n+1,\ldots,n-1,n\}$.
\begin{align*}
&q^{2n-1}(2)_q \sqrt{(n+t+2)_q(n+t+1)_q(n+t)_q(n+t-1)_q} \,d_{t+2}
\\
-
&
(4)_q \sqrt{(n+t)_q(n+t-1)_q(n-t)_q(n-t-1)_q} \, d_t
\\
+&
q^{-2n+1}(2)_q \sqrt{(n-t+2)_q(n-t+1)_q(n-t)_q(n-t-1)_q}\, d_{t-2}
=0,
\end{align*}
where $d_t=0$ if $|t|\geq n+1$. 
\end{lem}

\begin{proof}
Take the $\pi_2$-eigenvector $\xi_0^2=\sum_{s\in I_2} c_s \bw_{s}$ 
where $c_{\pm2}=q^{\mp1}\sqrt{(3)_q!}$, $c_{\pm1}=0$ and 
$c_0=-\sqrt{(4)_q}$. 
Then we have $\Psi_{n-2}(\xi_0^2,\zeta)=0$. 
Multiplying $q^{2n-4}\sqrt{(2n)_q\cdots(2n-3)_q}$ to 
the left hand side, we can derive
\begin{align*}
&\,\sum_{t\in I_n}
\Big{(}c_{-2}q^{2n-2t-2}\sqrt{(n+t+2)_q\cdots(n+t-1)_q}\,d_{t+2}\\
&\quad 
+c_0 q^{-2t}\sqrt{\frac{(4)_q(3)_q}{(2)_q}}
\sqrt{(n+t)_q(n+t-1)_q(n-t)_q(n-t-1)_q}\, d_t\\
&\quad 
+
c_2 q^{-2n-2t+2}\sqrt{(n-t+2)_q\cdots(n-t-1)_q}\,d_{t-2}
\Big{)}\bw_t^{n-2}=0.
\end{align*}
Then we obtain the desired formula. 
\end{proof}

Let us write 
\begin{align*}
\alpha_t=&\,q^{2n-1}(2)_q \sqrt{(n+t+2)_q(n+t+1)_q(n+t)_q(n+t-1)_q}\, ,\\
\beta_t=&\,-(4)_q \sqrt{(n+t)_q(n+t-1)_q(n-t)_q(n-t-1)_q}\, ,\\ 
\gamma_t=&q^{-2n+1}(2)_q \sqrt{(n-t+2)_q(n-t+1)_q(n-t)_q(n-t-1)_q}\, .
\end{align*} 
Then $\{d_t\}_{t\in I_n}$ satisfies
\begin{equation}\label{recu01}
\alpha_t d_{t+2}+\beta_t d_{t}+\gamma_t d_{t-2}=0\, .
\end{equation} 
By definition we obtain $\alpha_{-t}=q^{4n-2}\gamma_t$ and 
$\beta_{-t}=\beta_t$. 
From the above recurrence formula at $-t$, 
we obtain 
\[
q^{4n-2}\gamma_t d_{-t+2}+\beta_{t}d_{-t}+q^{-4n+2}\alpha_{t} d_{-t-2}=0\, .
\]
Since we have assumed the self-conjugacy 
$d_{-t}=(-q)^t\overline{d_{t}}$, 
we obtain the following formula.
\begin{equation}\label{recu02}
q^{-4n+4}\alpha_t d_{t+2}+\beta_{t}d_{t}+q^{4n-4}\gamma_t d_{t-2}=0\, .
\end{equation}
Through formulae (\ref{recu01}) and (\ref{recu02}), 
we obtain
\begin{equation}\label{recu03}
d_{t+2}=q^{-2}
\sqrt{\frac{(n-t+2)_q\cdots(n-t-1)_q}{(n+t+2)_q\cdots(n+t-1)_q}} d_{t-2}
\end{equation} 
for $t\geq-n+2$. 
Again from (\ref{recu01}), we have 
$\beta_t d_t=(-1-q^{4n-4})\gamma_t d_{t-2}$. 
Hence for $|t|\leq n-2$, we obtain the following formula.
\begin{equation}\label{recu04}
d_t=q^{-1}(q^{-2n+2}+q^{2n-2})\frac{(2)_q}{(4)_q}
\sqrt{\frac{(n-t+2)_q(n-t+1)_q}{(n+t)_q(n+t-1)_q}}\, d_{t-2}\, .
\end{equation}
Then we consider the above equalities (\ref{recu03}) 
and (\ref{recu04}). 
For $-n+2\leq t \leq n-4$ we obtain by using (\ref{recu04}) twice
\begin{equation}\label{recu05}
d_{t+2}=q^{-2}(q^{-2n+2}+q^{2n-2})^2\frac{(2)_q^2}{(4)_q^2}
\sqrt{\frac{(n-t+2)_q\cdots(n-t-1)_q}{(n+t+2)_q\cdots(n+t-1)_q}}\,d_{t-2}\, .
\end{equation}
By (\ref{recu03}) and (\ref{recu05}) we have a equality
\begin{equation}\label{recu06}
d_{t-2}=(q^{-2n+2}+q^{2n-2})^2\frac{(2)_q^2}{(4)_q^2}\, d_{t-2}
\end{equation}
for $-n+2\leq t \leq n-4$. 
We easily see that $1=(q^{-2n+2}+q^{2n-2})^2\frac{(2)_q^2}{(4)_q^2}$ holds 
if and only if $n=0$ or $2$, 
however, this does not occur 
because $n$ is an odd number. 
Thus we have $d_{t-2}=0$ for $-n+2\leq t \leq n-4$. 
It follows $d_{t}=0$ for $-n \leq t \leq n-6$. 
For $n\geq7$ we can derive $d_{t}=0$ immediately. 
So we have to consider the cases $n=3$ and $n=5$. 

(1) $n=5$ case. We have already known $d_{t}=0$ for $-5\leq t \leq -1$. 
Using (\ref{recu04}), we have $d_0=0$. 
Hence we have $d_t=0$ for all $t\in I_5$. 

(2) $n=3$ case. We have already known $d_{-3}=0$. 
Using (\ref{recu04}) for $t=-1,1$, we have $d_t=0$ for odd $t$. 
From (\ref{recu03}) and (\ref{recu04}) for $t=0$ we obtain 
\begin{align*}
d_{2}
=&\,
q^{-2}d_{-2}\, ,\\
d_0
=&\, 
q^{-1}(q^{-4}+q^4)\frac{(2)_q}{(3)_q}
\sqrt{\frac{(5)_q(4)_q}{(3)_q(2)_q}}\,d_{-2} \, .
\end{align*}
This shows $d_{-2}$ and $d_{2}$ are real numbers. 
Now we consider a $\pi_1$-eigenvector $\Psi_1(\eta,\eta)=0$. 
As in the case $A_4^*$, we obtain 
\begin{align*}
\sqrt{(6)_q}f\cdot \Psi_1(\eta,\eta)_{-1}
=&\,
q^{-2}(6)_q\eta_{3}\eta_{-3}
+\big{(}-q^{-1}(5)_q(2)_q +q^{-3}(6)_q \big{)}
\eta_{2}\eta_{-2}\\
&\quad\,
+\big{(}(4)_q(3)_q-q^{-2}(5)_q(2)_q\big{)}
\eta_{1}\eta_{-1}
+
(q^{-1}-q)(4)_q(3)_q \eta_{0}^2\\
&\quad\,
+
(q^2(5)_q(2)_q-(4)_q(3)_q)
\eta_{-1}\eta_{1}
+
(q(5)_q(2)_q-q^3 (6)_q)
\eta_{-2}\eta_{2}\\
&\quad\,
-q^2 (6)_q\eta_{-3}\eta_{3}
\, .
\end{align*}
Applying $\pi_\T$ to the above both sides, we have 
\[
0=
\big{(}(q-q^{-1})(5)_q(2)_q+(q^{-3}-q^3)(6)_q\big{)}d_{-2}d_2
+
(q^{-1}-q)(4)_q(3)_q d_0^2\, .
\]
Since $d_{-2}$ and $d_0$ are real numbers, we obtain
\[
q^{-2}((6)_q(3)_q-(5)_q(2)_q)d_{-2}^2+(4)_q(3)_qd_0^2=0\, ,
\]
however, the left hand side is positive because the above $d_{t}$ are real. 
Hence we can get $d_{-2}=d_0=d_2=0$. 

(VI) $\T_m\, (m\geq2)$ case. 
We treat a $C^*$-algebra $A$ whose spectral pattern is 
one of the following

$\T_{2\ell-1}$\, $(\ell\geq2)$ 
$\oplus_{k\in\Z_{\geq0}}
\big{(}1
+2\big{[}\frac{k}{2\ell-1}\big{]}\big{)}\pi_k
\oplus 
\oplus_{k\in \Z_{\geq0}} 
2\big{[}\frac{k+\ell}{2\ell-1}\big{]}\pi_{k+\frac{1}{2}}$
,

$\T_{2\ell}$\, $(\ell\geq1)$ 
$\oplus_{k\in \Z_{\geq0}}(1+2\big{[}\frac{k}{\ell}\big{]})\pi_k$
.

If $A$ is of type $\T_2$, 
then we can easily derive $A=C(\T_{2}\setminus\suq)$. 
Hence we study the case $\T_m$ for $m\geq3$. 
In the cases, the $\pi_1$-multiplicity are one, 
so the linear subspace $A_{\pi_1}$ generates a quantum sphere 
$\CSq$. 
We shall prove this is the canonical homogeneous sphere 
$C(\T\setminus\suq)$, 
that is, the parameter $\lambda_0=(q^{-1}-q)^{-1}\lambda$ 
is equal to $1$. 
In order to do this, we take the same strategy 
as in the case of type $D_m^*$. 
We have to prepare the following lemma 
which is proved in similar way to Lemma \ref{product of pi2} 
and we omit the proof. 
Recall a $\pi_1$-eigenvector of $\CSq$ 
\[\xi_\lambda^1
=q^{\frac{1}{2}}\sqrt{\frac{1-\lambda_0^2}{(2)_q}} \bw_{-1}^1
+
\lambda_0 \bw_0^1
-
q^{-\frac{1}{2}}\sqrt{\frac{1-\lambda_0^2}{(2)_q}} \bw_{1}^1
.
\]

\begin{lem}\label{product of pi1}
Let $n$ be a half integer with $n\geq1$. 
Then we have the following equalities for all $t\in I_n$.

$\Psi_{n-1}(\bw_{-1}^1,\bw_t^n)
=q^{-t+2}\sqrt{\frac{(n+t)_q(n+t-1)_q}{(2n)_q(2n-1)_q}}\bw_{t-1}^{n-1}\, ,$

$\Psi_{n-1}(\bw_{0}^1,\bw_t^n)
=-q^{-n-t+1}\sqrt{\frac{(2)_q(n-t)_q(n+t)_q}{(2n)_q(2n-1)_q}}\bw_t^{n-1}\, ,$

$\Psi_{n-1}(\bw_{1}^1,\bw_t^n)
=q^{-2n-t}\sqrt{\frac{(n-t)_q(n-t-1)_q}{(2n)_q(2n-1)_q}}\bw_{t+1}^{n-1}\, .$

\end{lem}

We immediately derive the following lemma. 

\begin{lem}\label{products of 1-eigenvector of quantum sphere}
Let $n$ be a half integer with $n\geq1$. 
Then we have the following equalities for all $t\in I_n$.
\begin{align*}
q^{n-1}\sqrt{(2)_q(2n)_q(2n-1)_q} &\, \Psi_{n-1}(\xi_{\lambda}^1,\bw_t^n)\\
&=
\ q^{n-t+\frac{3}{2}}
\sqrt{1-\lambda_0^2}
\,\sqrt{(n+t)_q(n+t-1)_q}\,\bw_{t-1}^{n-1}\\
&\quad
-q^{-t}(2)_q\lambda_0 
\sqrt{(n-t)_q(n+t)_q} \,\bw_t^{n-1}\\
&\quad
-q^{-n-t-\frac{3}{2}}
\sqrt{1-\lambda_0^2}
\,\sqrt{(n-t)_q(n-t-1)_q}\, \bw_{t+1}^{n-1}
\, .
\end{align*}

\end{lem}

We need another lemma whose proof is also same as for 
Lemma \ref{kernel of pi2-product}. 

\begin{lem}\label{kernel of pi1-product}
Let $n$ be a half integer with $n\geq1$ and 
$\eta=\sum_{t\in I_n}d_t \bw_t^n$ be a $\pi_n$-eigenvector of 
$\csuq$ such that $\Psi_{n-1}(\xi_\lambda^1,\eta)=0$. 
Then $\{d_t\}_{t\in I_n}$ satisfies the following recurrence 
equation for all $t\in I_n=\{-n,-n+1,\ldots,n-1,n\}$.
\begin{align*}
&
q^{n+\frac{1}{2}}\sqrt{1-\lambda_0^2}
\sqrt{(n+t+1)_q(n+t)_q}\,d_{t+1}\\
&\,
-(2)_q\lambda_0 \sqrt{(n+t)_q(n-t)_q}\,d_t\\
&\,
-q^{-n-\frac{1}{2}}\sqrt{1-\lambda_0^2}
\sqrt{(n-t+1)_q(n-t)_q}\,d_{t-1}=0
\end{align*}
where we define $d_t=0$ if $|t|\geq n+1$. 
\end{lem}

We show $A=C(\T_m\setminus\suq)$ for 
in each case of odd and even order cyclic groups. 
 
(1) $\T_{2\ell-1}\, (\ell\geq 2)$ case. 
We focus on $\pi_{\ell-\frac{1}{2}}$-spectral subspace. 
Write $n$ for $\ell-\frac{1}{2}$. 
Since $\pi_n$-eigenvector space is two-dimensional 
and $\pi_{n-1}$-eigenvector space is $0$, 
the map 
$\Psi_{n-1}(\xi_\lambda^1,\cdot):X_{\pi_n}\longrightarrow X_{\pi_{n-1}}$ 
is $0$-map. 
Let $\eta=\sum_{t\in I_n}d_t \bw_t^n$ be 
a non-zero $\pi_n$-eigenvector for $A$. 
Recall its conjugate eigenvector 
$T\eta=\sum_{t\in I_n}(-q)^{-t}\overline{d_{-t}}\bw_t^n$, 
where we have defined $(-q)^s=\sqrt{-1}^{2s}|q|^s$ for 
all real number $s$. 
Then $\eta$ and $T\eta$ 
are in the kernel of $\Psi_{n-1}(\xi_\lambda^1,\cdot)$. 
By Lemma \ref{kernel of pi1-product} 
we obtain the following recurrence equation.
\begin{align}
\alpha_t d_{t+1}+\beta_t d_t+ \gamma_t d_{t-1}=&\,0\, 
\label{Todd01}\\
q^{-2n}\alpha_t d_{t+1}+\beta_t d_t+ q^{2n}\gamma_t d_{t-1}=&\,0\, ,
\label{Todd02}
\end{align}
where we put 
\begin{align*}
\alpha_t=&\,q^{n+\frac{1}{2}}\sqrt{1-\lambda_0^2}
\sqrt{(n+t+1)_q(n+t)_q}\, ,\\
\beta_t=&\,-(2)_q\lambda_0 \sqrt{(n+t)_q(n-t)_q}\, ,\\
\gamma_t=&\,-q^{-n-\frac{1}{2}}\sqrt{1-\lambda_0^2}
\sqrt{(n-t+1)_q(n-t)_q}\, .
\end{align*}
Assume that $\lambda_0$ is not equal to $1$. 
Then we know $\alpha_t=0$ if and only if $t=-n$ and 
$\gamma_t=0$ if and only if $t=n$. 
By (\ref{Todd01}) and (\ref{Todd02}) we obtain 
\begin{equation}\label{Todd03}
\alpha_t d_{t+1}=q^{2n}\gamma_{t}d_{t-1}\, .
\end{equation}
Hence again by (\ref{Todd01}) we have 
\[
\beta_t d_t=-(1+q^{2n})\gamma_t d_{t-1}\, .
\]
Assume $\lambda_0$ is equal to $0$. 
Then the above equality shows $d_{t-1}=0$ for $-n\leq t\leq n-1$. 
By (\ref{Todd03}) $d_{n-1}$ and $d_n$ are also equal to $0$, 
however, this is a contradiction to non triviality of $\eta$. 
Next we consider the case $0<\lambda_0<1$. 
From the above equality 
$d_t$ for $-n+1 \leq t \leq n-1$ are uniquely determined by $d_{-n}$. 
The number $d_n$ is determined by the previous equation as 
$d_n=q^{2n}\alpha_{n-1}^{-1}\gamma_{n-1}d_{n-2}$ 
where we use the fact that $\alpha_{n-1}$ is never $0$ 
for $n\geq \frac{3}{2}$. 
Therefore $(d_t)_{t\in I_n}$ is a scalar multiple of a vector. 
In particular the $\pi_n$-eigenvector space $X_{\pi_n}$ becomes 
one-dimensional. 
This is a contradiction. 
It concludes $\lambda_0$ is equal to $1$ in this case. 
From simple computation 
we obtain $X_n=\C\bw_{-n}^n+\C\bw_{n}^n$. 
We show $A=C(\T_{2\ell-1}\setminus \suq)$. 
The $C^*$-algebra $C(\T_{2\ell-1}\setminus \suq)$ is 
generated by the matrix elements $w(\pi_\nu)_{ns,t}$ 
for all $\nu\in\frac{1}{2}\Z_{\geq0}$ and $t\in I_\nu$, 
where 
$s$ is an integer such that $ns\in I_\nu$. 
If $A$ contains $C(\T_{2\ell-1}\setminus \suq)$, then 
they coincide because of its common spectral pattern. 
By $\suq$-invariance of $A$, it suffices to show 
$w(\pi_\nu)_{ns,-\nu} 
={2\nu \brack \nu+ns}_{q^2}^{\frac{1}{2}}x^{\nu+ns}v^{\nu-ns}$ 
are contained in $A$ for non-negative $s$. 
This is equal to 
${2\nu \brack \nu+ns}_{q^2}^{\frac{1}{2}}x^{2ns}x^{\nu-ns}v^{\nu-ns}$. 
Since $A$ contains $xv\in C(\T\setminus\suq)$ and 
$x^{2n}={2\nu \brack \nu+ns}_{q^2}^{-\frac{1}{2}}w(\pi_n)_{-n,-n}$, 
the element $x^{\nu+ns}v^{\nu-ns}$ is in $A$. 
This shows $A$ is in fact $C(\T_{2\ell-1}\setminus \suq)$. 

(2) $\T_{2\ell}\, (\ell\geq 2)$ case. 
Write $n$ for $\ell$ in this case. 
We analyze the recurrence equation 
in Lemma \ref{kernel of pi1-product} 
under the condition that a $\pi_n$-eigenvector 
$\eta=\sum_{t\in I_n}d_t \bw_t^n$ is self-conjugate. 
Assume that $\lambda_0$ is not equal to $1$. 
When $\lambda_0$ is equal to $0$, $\eta$ is the zero vector 
as in the previous case. 
In the case $0<\lambda_0<1$ 
we know that the space of its solution is one dimensional 
and we can assume that $d_t$ are all real number. 
Then we have 
\[
d_t=-(1+q^{2n})\beta_t^{-1}\gamma_t d_{t-1}
\]
for all $|t|\leq n-1$. 
For $-n+1 \leq t \leq n-2$ we obtain the following equality.
\begin{align*}
d_{-t}
=&\, 
-(1+q^{2n})\frac{\gamma_{-t}}{\beta_{-t}}d_{-t-1}\\
=&\,
-(-q)^{t+1}(1+q^{2n})\frac{\gamma_{-t}}{\beta_{-t}}d_{t+1}\\
=&\,
(-q)^{t+1}(1+q^{2n})^2 
\frac{\gamma_{-t}\gamma_{t+1}}{\beta_{-t}\beta_{t+1}}d_{t}\\
=&\,
-q(1+q^{2n})^2 
\frac{\gamma_{-t}\gamma_{t+1}}{\beta_{-t}\beta_{t+1}}d_{-t}\, .
\end{align*}
Since its coefficient 
$\frac{\gamma_{-t}\gamma_{t+1}}{\beta_{-t}\beta_{t+1}}$ 
is strictly positive, 
we have $d_{t}=0$ for $-n+2 \leq t \leq n-1$. 
We now treat 
a positive integer $n\geq 2$. 
It deduces $-n+2\leq n-2$ and 
we obtain $d_n=0$ by (\ref{Todd03}). 
Hence we have $d_t=0$ for all $t\in I_n$ 
because of $-n+2\leq 0$ and the self-conjugacy of $\eta$. 
Under the condition $\lambda_0<1$ 
we have shown $\eta=0$ 
if 
$\eta$ is self-conjugate and $\Psi_{n-1}(\xi_\lambda^1,\eta)=0$. 
Now we consider the map 
$\Psi_{n-1}(\xi_\lambda^1,\cdot):X_{n}\longrightarrow X_{n-1}$. 
We have known $X_{n}$ is three-dimensional and $X_{n-1}$ is 
one-dimensional. 
Take two linearly independent $\pi_n$-eigenvectors 
$\zeta_1$ and $\zeta_2$ from $X_n$ 
which satisfy 
$\Psi_{n-1}(\xi_\lambda^1,\zeta_1)=\Psi_{n-1}(\xi_\lambda^1,\zeta_2)=0$. 
For conjugate eigenvectors $T\zeta_1$ and $T\zeta_2$, 
we can take complex numbers $\mu_1$ and $\mu_2$ which are defined by 
$\Psi_{n-1}(\xi_\lambda^1,T\zeta_1)=\mu_1 
\Psi_{n-1}(\xi_\lambda^1,\xi_\lambda^n)$ 
and 
$\Psi_{n-1}(\xi_\lambda^1,T\zeta_2)=\mu_2 
\Psi_{n-1}(\xi_\lambda^1,\xi_\lambda^n)$, 
where $\xi_{\lambda}^\nu$ is a self-conjugate $\pi_\nu$-eigenvector 
for $C(S_{q,\lambda}^2)$. 
Note $\Psi_{n-1}(\xi_\lambda^1,\xi_\lambda^n)$ is not the zero vector. 
Multiplying complex numbers to $\zeta_1$ and $\zeta_2$ respectively, 
we may assume $\mu_1$ and $\mu_2$ are real numbers. 
If $\mu_1$ is equal to $0$, 
then the self-conjugate vectors $\zeta_1+T\zeta_1$ and 
$\sqrt{-1}(\zeta_1-T\zeta_1)$ are in the kernel of 
$\Psi_{n-1}(\xi_\lambda^1,\cdot)$. 
From the above discussion they are $0$. 
This shows immediately $\zeta_1=0$. 
This is a contradiction. 
Hence we have a non-zero real number $\mu_1$. 
Similarly we show $\mu_2$ is non-zero. 
We may assume they are equal to $1$. 
For $i=1,2$ we have 
$\Psi_{n-1}(\xi_\lambda^1,\zeta_i+T\zeta_i-\xi_\lambda^n)=0$. 
Then it yields $\zeta_i+T\zeta_i=\xi_\lambda^n$ 
because of the self-conjugacy of $\zeta_i+T\zeta_i-\xi_\lambda^n$ 
for each $i$. 
It follows $\zeta_1-\zeta_2=-T(\zeta_1-\zeta_2)$. 
In particular $T(\zeta_1-\zeta_2)$ is 
in the kernel of $\Psi_{n-1}(\xi_\lambda^1,\cdot)$. 
This shows $\zeta_1-\zeta_2=0$ 
but this is a contradiction to linear independence of 
$\zeta_1$ and $\zeta_2$. 
Therefore we can derive $\lambda_0=1$. 
Then we obtain $X_n=\C\bw_{-n}^n+\C\bw_0^n+\C\bw_{n}^n$. 
As in the proof of the case $\T_{2\ell-1}$, we can prove 
$A=C(\T_{2\ell}\setminus\suq)$.

\section{Classification of right coideals 
of $\csuq$: $-1<q<0$ case}

We use the same strategy as in the previous section 
in order to classify right coideals of $\csuq$ for 
negative $q$. 
We prepare a positive parameter $q_0$ defined by $q_0=-q$. 
Recall we have defined $q$-integer 
$(n)_q$ by $q_0$-integer, more precisely 
$(n)_q=(n)_{q_0}$ and 
$q^n=\sqrt{-1}^{2n}q_0^{n}$ for a half integer $n$. 

Now we investigate each case of multiplicity diagrams. 
Contrary to the case of positive $q$, 
we have to treat other graphs 
which have a single loop at a vertex as 
listed in Appendix. 
We recall the notion of the quantum spheres $\CSq$. 
We use an embedding it into $\csuq$ by defining 
its $\pi_1$-eigenvector as 
\[
\xi_\lambda^1
=\bigg{(}q_0^{\frac{1}{2}}\sqrt{\frac{1-\lambda_0^2}{(2)_q}},\lambda_0,
q_0^{-\frac{1}{2}}\sqrt{\frac{1-\lambda_0^2}{(2)_q}}\bigg{)}\, w(\pi_1) \, .
\]

If a right coideal $A$ is of type $\T$, 
then it is in fact $\beta^L_{z}(\CSq)$ for some $z\in \T$ 
by Podle\'{s}' classification results 
on the quantum spheres \cite[Theorem 1]{Podles1}. 
We know that $\csoq$ and $\csomq$ are isomorphic 
as compact quantum groups. 
Hence we have the complete collection of right coideals 
which have only spectral subspaces with integer spins. 
In this way we can deny the existence of type 
$A_4^*$, $S_4^*$, $A_5^*$ and $D_m^*\, (m\geq2)$. 
For right coideals of $D_\infty^*$ and $\T_{2n}\,(n\geq1)$
we have shown their uniqueness up to conjugation by $\beta^L$. 
Hence we have to investigate right coideals of 
type $\T_{n}$ (odd $n\geq3$), $D_n$ (odd $n\geq1$), 
$A_m'\, (m\geq3)$. 
The main classification result is as follows. 

\begin{thm}\label{main negative}
Let $A\subset \csuq$ be a right coideal. 
Then its multiplicity diagram is one of 
type $1,\T_n\,(n\geq2),\T,SU(2),D_\infty^*$ 
and $D_1$. 
If it is of type $\T$, it is one of the quantum spheres. 
Otherwise it is unique up to conjugation by $\beta^L$. 
\end{thm}

(I) $\T_{2\ell-1}\, (\ell\geq2)$ case. 
We treat a $C^*$-algebra $A$ whose spectral pattern is 
$\oplus_{k\in\Z_{\geq0}}\big{(}1+2\big{[}\frac{k}{2\ell-1}\big{]}\big{)}\pi_k
\oplus 
\oplus_{k\in \Z_{\geq0}} 
2\big{[}\frac{k+\ell}{2\ell-1}\big{]}\pi_{k+\frac{1}{2}}$. 
As in the previous section we study the $C^*$-subalgebra 
generated by $\pi_1$-spectral subspace and assert 
that is the canonical homogeneous sphere. 
We begin classifying this case 
by stating a negative $q$ version of Lemma \ref{product of pi1}. 

\begin{lem}\label{product of pi1 negative}
Let $n$ be a half integer with $n\geq1$. 
Then we have the following equalities for all $t\in I_n$.

$\Psi_{n-1}(\bw_{-1}^1,\bw_t^n)
=q_0^{-t+2}\sqrt{\frac{(n+t)_q(n+t-1)_q}{(2n)_q(2n-1)_q}}\,\bw_{t-1}^{n-1}\, ,$

$\Psi_{n-1}(\bw_{0}^1,\bw_t^n)
=
-(-1)^{n-t}
q_0^{-n-t+1}\sqrt{\frac{(2)_q(n-t)_q(n+t)_q}{(2n)_q(2n-1)_q}}\,\bw_t^{n-1}\, ,$

$\Psi_{n-1}(\bw_{1}^1,\bw_t^n)
=q_0^{-2n-t}\sqrt{\frac{(n-t)_q(n-t-1)_q}{(2n)_q(2n-1)_q}}\,\bw_{t+1}^{n-1}\, .$

\end{lem}

If necessary, we consider the conjugation by $\beta^L$ and 
may assume $\CSq\subset A$. 

\begin{lem}
Let $n$ be a half integer with $n\geq1$. 
Then we get the following equality for $t\in I_n$. 
\begin{align*}
\Psi_{n-1}(\xi_\lambda^1,\bw_t^n)
=&\,
q_0^{-t+\frac{5}{2}}\sqrt{1-\lambda_0^2}
\sqrt{\frac{(n+t)_q(n+t-1)_q}{(2)_q(2n)_q(2n-1)_q}}\,\bw_{t-1}^{n-1}\\
&\,
-
(-1)^{n-t}
q_0^{-n-t+1}\lambda_0
\sqrt{\frac{(2)_q(n-t)_q(n+t)_q}{(2n)_q(2n-1)_q}}\,\bw_t^{n-1}\\
&\,
+
q_0^{-2n-t-\frac{1}{2}}\sqrt{1-\lambda_0^2}
\sqrt{\frac{(n-t)_q(n-t-1)_q}{(2)_q(2n)_q(2n-1)_q}}
\,\bw_{t+1}^{n-1}\, .
\end{align*}

\end{lem}

\begin{lem}
Let $n$ be a half integer with $n\geq1$ 
and $\eta=\sum_{t\in I_n}d_t \bw_t^n$ be a $\pi_n$-eigenvector 
with $\Psi_{n-1}(\xi_\lambda^1,\eta)=0$. 
Then we get the following recurrence equation for $t\in I_n$.
\begin{align*}
&
q_0^{n+\frac{1}{2}}\sqrt{1-\lambda_0^2}
\sqrt{(n+t)_q(n+t+1)_q}\, d_{t+1}\\
&\,
-
(-1)^{n-t}\lambda_0 (2)_q \sqrt{(n-t)_q(n+t)_q}\, d_t\\
&\,
+
q_0^{-n-\frac{1}{2}}\sqrt{1-\lambda_0^2}
\sqrt{(n-t)_q(n-t+1)_q}\, d_{t-1}=0\, ,
\end{align*}
where we define $d_t=0$ for $|t|\geq n+1$. 
\end{lem}

We write $n$ for $\ell-\frac{1}{2}$. 
Then $\pi_n$-eigenvector space of $A$ is two-dimensional. 
We may assume $\pi_1$-part generates a quantum sphere $C(S_{q,\lambda}^2)$ 
as usual. 
We derive $\lambda_0=1$ and it immediately derives 
$A=C(\T_{2\ell-1}\setminus\csuq)$. 
Now Let $\eta=\sum_{t\in I_n}d_t \bw_t^n$ be a 
$\pi_n$-eigenvector of $A$. 
By the absence of $\pi_{n-1}$-part we have 
$\Psi_{n-1}(\xi_\lambda^1,\eta)=0=\Psi_{n-1}(\xi_\lambda^1,T\eta)$. 
From the previous lemmas we get 
\begin{align*}
\alpha_t d_{t+1}+\beta_t d_t+\gamma_t d_{t-1}=&\,0\,, \\
q_0^{-2n}\alpha_t d_{t+1}-\beta_t d_t
+q_0^{2n}\gamma_t d_{t-1}=&\,0\,, 
\end{align*}
where we put 
\begin{align*}
\alpha_t
=&\,
q_0^{n+\frac{1}{2}}\sqrt{1-\lambda_0^2}
\sqrt{(n+t)_q(n+t+1)_q}\, , \\
\beta_t
=&\,
-
(-1)^{n-t}\lambda_0 (2)_q \sqrt{(n-t)_q(n+t)_q}\, , \\
\gamma_t
=&\, 
q_0^{-n-\frac{1}{2}}\sqrt{1-\lambda_0^2}
\sqrt{(n-t)_q(n-t+1)_q}\, .
\end{align*}

Then we can prove $(d_t)_{t\in I_n}$ is a scalar multiple of 
a vector as in the $\T_{2\ell-1}$ case in the previous section. 
This is a contradiction. 
Hence we have proved $\lambda_0=1$ and 
$A=C(\T_{2\ell-1}\setminus \suq)$. 

(II) $D_n$ (odd $n\geq3$) case. 
Before proof of the non-existence of this case, 
we shall state some basic lemmas without proofs 
which are proved similarly by direct calculations 
as before. 

We state a negative $q$ and $\Psi_{n-1}$ version of Lemma 
\ref{product of pi2}. 

\begin{lem}\label{product of pi2 negative}
Let $n$ be a half integer with $n \geq \frac{3}{2}$. 
Then we have the following equalities for all $t\in I_n$. 

\begin{align*}
&
\Psi_{n-1}(\bw_{-2}^2,\bw_t^n)
=
(-1)^{n-t}q_0^{-\frac{1}{2}n-2t+6}\frac{1}{(2n-2)_q}\\
&\hspace{40mm}
\sqrt{
\frac{(4)_q(n-t+1)_q(n+t)_q(n+t-1)_q(n+t-2)_q}
{(2n)_q(2n-1)_q}
}\,\bw_{t-2}^{n-1}
\, ,\\
&
\Psi_{n-1}(\bw_{0}^2,\bw_t^n)
=(-1)^{n-t}q_0^{-\frac{3}{2}n-2t+3} 
\frac{q_0^{-n-1}(n-t-1)_q-q_0^{n+1}(n+t-1)_q}{(2n-2)_q}\\
&\hspace{40mm}
\sqrt{
\frac{(3)_q!(n+t)_q(n-t)_q}
{(2n)_q(2n-1)_q}
}\,\bw_{t}^{n-1}
\, ,\\
&
\Psi_{n-1}(\bw_{2}^2,\bw_t^n)
=
-(-1)^{n-t}
q_0^{-\frac{5}{2}n-2t}
\frac{1}{(2n-2)_q}\\
&\hspace{40mm}
\sqrt{
\frac{(4)_q(n+t+1)_q(n-t)_q(n-t-1)_q(n-t-2)_q}
{(2n)_q(2n-1)_q}
}\,\bw_{t+2}^{n-1}
\, .
\end{align*}
\end{lem}

We take a vector 
$\xi_0^2=(q_0\sqrt{(3)_q!},0,-\sqrt{(4)_q},0,q_0^{-1}\sqrt{(3)_q!}\,)\,
w(\pi_2)$ 
as a $\pi_2$-eigenvector of $C(S_{q,0}^2)$. 
Recall its entries generate the right coideal of type $D_\infty^*$. 
We denote it by $A_{D_\infty^*}$. 

\begin{lem}

Let $n$ be a half integer with $n \geq \frac{3}{2}$. 
Then we have the following equalities for all $t\in I_n$. 
\begin{align*}
&-q_0^{\frac{3}{2}n-3}\sqrt{\frac{(2n)_q(2n-1)_q}{(4)_q!}}
\Psi_{n-1}(\xi_0^2,\bw_t^n)\\
&\hspace{3mm}
=
-(-1)^{n-t}q_0^{n-2t+4}\sqrt{(n-t+1)_q(n+t)_q(n+t-1)_q(n+t-2)_q}
\,\bw_{t-2}^{n-1}\\
&\hspace{8mm} 
+
(-1)^{n-t}q_0^{-2t}
(q_0^{-n-1}(n-t-1)_q-q_0^{n+1}(n+t-1)_q)
\sqrt{(n+t)_q(n-t)_q}\bw_t^{n-1}\\
&\hspace{8mm}
+
(-1)^{n-t}q_0^{-n-2t-4}
\sqrt{(n+t+1)_q(n-t)_q(n-t-1)_q(n-t-2)_q}
\,\bw_{t+2}^{n-1}\, .
\end{align*}

\end{lem}

\begin{lem}
Let $n$ be a half integer with $n \geq \frac{3}{2}$ 
and $\eta=\sum_{t\in I_n}d_t \bw_t^n$ 
be a $\pi_n$-eigenvector of $A$ such that 
it satisfies $\Psi_{n-1}(\xi_0^2,\eta)=0$. 
Then we have the following recurrence equation for $t\in I_n$. 
\begin{align*}
&-q_0^n \sqrt{(n-t-1)_q(n+t+2)_q(n+t+1)_q(n+t)_q}\, d_{t+2}\\
&\,
+
(q_0^{-n-1}(n-t-1)_q-q_0^{n+1}(n+t-1)_q)\sqrt{(n+t)_q(n-t)_q}\, d_t\\
&\,
+
q_0^{-n}\sqrt{(n+t-1)_q(n-t+2)_q(n-t+1)_q(n-t)_q}\, d_{t-2}
=0, 
\end{align*}
where we define $d_t=0$ for $|t|\geq n+1$. 
\end{lem}

\begin{lem}\label{pi2 product gap criterion negative}
Let $n$ be a half integer in $\frac{3}{2}+\Z_{\geq0}$ and $A$ be 
a right coideal in $\csuq$ such that $A$ contains $A_{D_\infty^*}$. 
If a $\pi_n$-eigenvector $\eta$ 
of $A$ satisfies $\Psi_{n-1}(\xi_0^2,\eta)=0$ and 
$\Psi_{n-1}(\xi_0^2,T\eta)=0$, then $\eta$ is the zero vector. 
\end{lem}

\begin{proof}
Take complex numbers $\{d_t\}_{t\in I_n}$ with 
$\eta=\sum_{t\in I_n}d_t \bw_t^n$. 
Note that $T\eta=\sum_{t\in I_n} (-q)^{-t}\overline{d_{-t}}\bw_t^n$. 
Let us prepare the following notations 
$\alpha_t$, $\beta_t$ and $\gamma_t$.
\begin{align*}
\alpha_t
=&\, 
-q_0^n \sqrt{(n-t-1)_q(n+t+2)_q(n+t+1)_q(n+t)_q}\\
\beta_t
=&\,
(q_0^{-n-1}(n-t-1)_q-q_0^{n+1}(n+t-1)_q)\sqrt{(n+t)_q(n-t)_q}\\
\gamma_t
=&\,
q_0^{-n}\sqrt{(n+t-1)_q(n-t+2)_q(n-t+1)_q(n-t)_q}\, .
\end{align*}
From our assumption on $\eta$ and $T\eta$, we get the 
following recurrence equations for $t\in I_n$ 
by the previous lemma. 
\begin{align}
\alpha_t d_{t+2}+\beta_t d_t +\gamma_t d_{t-2}=&\,0 \label{D"01}\,,\\
q_0^{-2n+2}\alpha_t d_{t+2}-\beta_{-t}d_t+q_0^{2n-2}\gamma_t d_{t-2}=&\,0
\notag\, .
\end{align}
In the above equations put $t=-n+1$ and we get 
\begin{align*}
-q_0^{n}\sqrt{(2n-2)_q(3)_q!}\,d_{-n+3}
+q_0^{-n-1}(2n-2)_q\sqrt{(2n-1)_q}\,d_{-n+1}
=&\,0\,,\\
-q_0^{-n+2}\sqrt{(2n-2)_q(3)_q!}\,d_{-n+3}
+q_0^{n+1}(2n-2)_q\sqrt{(2n-1)_q}\,d_{-n+1}
=&\,0\, .
\end{align*}
This derives $d_{-n+1}=d_{-n+3}=0$ because $0<q_0<1$. 
We know $\alpha_t=0$ if and only if $t=-n,n-1$. 
By using (\ref{D"01}) inductively, we get $d_{-n+2k-1}=0$ for 
$k=1,\ldots n+\frac{1}{2}$. 
Similarly we get $d_{n-2k+1}=0$ for $k=1,\ldots n+\frac{1}{2}$. 
Hence we have proved $\eta=0$. 
\end{proof}

Let $A$ be a right coideal of type $D_n$ (odd $n\geq3$). 
Then its spectral pattern is 
\[\bigoplus_{k\in \Z_{\geq0}}
\bigg{(}\bigg{[}\frac{k}{n}\bigg{]}+\frac{1+(-1)^k}{2}\bigg{)}\pi_k
\oplus
\bigoplus_{k\in \Z_{\geq0}}
\bigg{(}\bigg{[}\frac{2k+1}{n}\bigg{]}
-\bigg{[}\frac{k}{n}\bigg{]}\bigg{)}
\pi_{k+\frac{1}{2}}\, . 
\]
If we look at the integer spin spectral pattern, 
then we have $\pi_0\oplus\pi_2\oplus\cdots$. 
Therefore the $\pi_2$-spectral subspace $A_{\pi_2}$ generates 
a right coideal of type $D_\infty^*$. 
Considering $\beta_z^L(A)$ for some $z\in \T$, 
we may assume $A$ contains $A_{D_\infty^*}$. 
Next look at the half integer part and we see that 
$\pi_{\frac{n}{2}}$ appears once and $\pi_{\frac{n}{2}-1}$ 
does not there. 
Then we can make use of Lemma \ref{pi2 product gap criterion negative} 
and conclude such a right coideal does not exist. 

(III) $D_1$ case. 
We show the existence of this case and its uniqueness 
up to conjugation. 
In order to do we need some elementary lemmas. 

\begin{lem}
Take a positive integer $n$. 
Then we have the following equalities. 
\[
f\cdot x^n =\sqrt{-1}^{\,-n+1}q_0^{\frac{n-1}{2}}(n)_q x^{n-1}u \, ,
\quad
f\cdot v^n =\sqrt{-1}^{\,-n+1}q_0^{\frac{n-1}{2}}(n)_q v^{n-1}y\, .
\] 
\end{lem}
\begin{proof}
We know $x^n=w(\pi_\frac{n}{2})_{-\frac{n}{2},-\frac{n}{2}}$. 
Hence we have:
\begin{align*}
f\cdot x^n
=&\,
\sqrt{-1}^{\,-n+1}\sqrt{(n)_q}
w(\pi_\frac{n}{2})_{-\frac{n}{2},-\frac{n}{2}+1}\\
=&\,
\sqrt{-1}^{\,-n+1}\sqrt{(n)_q}
 { n \brack 1}^{\frac{1}{2}}x^{n-1}u\\
=&\,
\sqrt{-1}^{\,-n+1}
 q_0^{\frac{n-1}{2}}(n)_qx^{n-1}u\, .
\end{align*}
Similarly we get the assertion on $f\cdot v^n$. 
\end{proof}

\begin{lem}
For $r,s\geq0$ we have the following equality.
\[
h(x^r y^r u^s v^s)= q_0^{r(s+1)}\frac{(r)_q!(s)_q!}{(r+s+1)_q!}      
\, .
\]
\end{lem}
\begin{proof}
Recall $f\cdot y^r=f\cdot u^{s-1}=0$ and we get 
\begin{align*}
f\cdot x^{r+1}y^r u^{s-1} v^s
=&\,
(f\cdot x^{r+1})\cdot(k\cdot y^r u^{s-1} v^s)
+(k^{-1}\cdot x^{r+1})\cdot(f\cdot y^r u^{s-1} v^s)\\
=&\,
\sqrt{-1}^{\,-r}q_0^{\frac{r}{2}}(r+1)_q
\cdot
q^{-\frac{r}{2}-\frac{s-1}{2}+\frac{s}{2}} x^{r}u y^r u^{s-1} v^s\\
&\,
+
q^{\frac{-(r+1)}{2}}\cdot 
q^{\frac{r}{2}+\frac{s-1}{2}}
\sqrt{-1}^{\,-s+1}q_0^{\frac{s}{2}}(s)_q
x^{r+1}y^r u^{s-1} v^{s-1}y\\
=&\,
\sqrt{-1}^{\,-r}q_0^{\frac{r}{2}}(r+1)_q
\cdot
q^{-\frac{3r}{2}+\frac{1}{2}} x^{r}y^r u^{s} v^s\\
&\,
+
q^{\frac{s}{2}-1}
\sqrt{-1}^{\,-s+1}q_0^{\frac{s}{2}}(s)_q
q^{-2(s-1)}
x^{r+1}y^{r+1} u^{s-1} v^{s-1}\\
=&\,
\sqrt{-1} \, q_0^{-r+\frac{1}{2}}(r+1)_q
\cdot
x^{r}y^r u^{s} v^s\\
&\,
-
\sqrt{-1} \, q_0^{-s+\frac{1}{2}}(s)_q
x^{r+1}y^{r+1} u^{s-1} v^{s-1} \, .
\end{align*}
The Haar state $h$ has the property that 
$h(f\cdot a)=0$ for $a\in \asuq$. 
Apply $h$ to the above both sides and we have
\[
h(x^{r+1}y^{r+1}u^{s-1}v^{s-1})
=
q_0^{s-r}
\frac{(r+1)_q}{(s)_q}
h(x^{r}y^{r}u^{s}v^{s})\, .
\]
Inductively we can calculate as desired. 
\end{proof}

Recall the $\pi_{\frac{1}{2}}$-spectral subspace 
$Y_{\frac{1}{2}}=\C x+\C v +\C u +\C y$. 
Let $P_{\frac{1}{2}}:\csuq\longrightarrow Y_{\frac{1}{2}}$ be 
the projection introduced in the second section. 
It is an orthogonal projection with respect to 
the Haar state. 

\begin{lem}\label{proj half}

For $r,s\geq0$ we have the following equalities. 
\begin{align*}
P_{\frac{1}{2}}(x^r y^r u^s v^s u)
=
&\,
q_0^{r(s+2)}\frac{(2)_q(r)_q!(s+1)_q!}{(r+s+2)_q!} u\, ,\\
P_{\frac{1}{2}}(x^r y^r u^s v^s y)
=
&\,
q_0^{r(s+1)-s}\frac{(2)_q(r+1)_q!(s)_q!}{(r+s+2)_q!}u\, .
\end{align*}

\end{lem}

\begin{proof}
Recall two maps $\beta^L =(\pi_\T\otimes \id)\circ \delta$ 
and $\beta^R=(\id\otimes \pi_\T)\circ \delta$ 
which act on $x,v,u$ and $y$ as follows.
\begin{align*}
\begin{pmatrix}
  \beta_z^L(x) &  \beta_z^L(u)   \\
  \beta_z^L(v) &  \beta_z^L(y)
\end{pmatrix}
=&\,
\begin{pmatrix}
  z x& z u\\
  \bar{z} v& \bar{z} y\\
\end{pmatrix}\, ,\\
\begin{pmatrix}
  \beta_z^R(x) &  \beta_z^R(u)   \\
  \beta_z^R(v) &  \beta_z^R(y)
\end{pmatrix}
=&\,
\begin{pmatrix}
  z x& \bar{z} u\\
  z v& \bar{z} y\\
\end{pmatrix}\, ,
\end{align*}
where $z$ runs on the torus $\T\subset \C$. 
Hence the above four elements have the different spectrums 
with respect to $\beta^L$ and $\beta^R$. 
It is easy to see that the element 
$x^r y^r u^s v^s u$ has the same spectrum as $u$. 
Hence it has the expansion in $L^2(\suq)$ 
$x^r y^r u^s v^s u=\lambda u +\cdots$ where $\lambda$ is 
a complex number. 
By orthogonality of $u$ and $w(\pi_\nu)_{i,j}$ 
for $\nu\geq1$ and $i,j\in I_\nu$, we get
\begin{align*}
\lambda
=&\,
h(uu^*)^{-1}h(x^r y^r u^s v^s uu^*)\\
=&\,
q_0(2)_q 
(-q^{-1})
h(x^r y^r u^{s+1} v^{s+1})\\ 
=&\
q_0^{r(s+2)}
\frac{(2)_q(r)_q!(s+1)_q!}{(r+s+2)_q!}\, .
\end{align*}
We also get the desired result on $y$ in the same way. 

\end{proof}

We propose the following lemma. 
Later we see this equality guarantees the existence of the 
right coideal of type $D_1$. 
Let us write $p_r=\prod_{t=1}^r (1+q^{-t})$. 

\begin{lem}\label{binomial}

For $k\geq0$ we have the following equality.
\[
\sum_{r=0}^k {k \brack r}_q (-1)^{k-r}q^{r(r-k+1)}p_r p_{k-r+1}
=
q^{-k}\sum_{r=0}^k {k \brack r}_q 
(-1)^{k-r}q^{r(r-k+2)}p_{r+1} p_{k-r}\, .
\]
Moreover, this is precisely equal to 
$(1+q^{-1})(q^{-2};q^{-1})_k$. 
\end{lem}

\begin{proof}
Denote the above left and right hand side by $a_k$ and $b_k$, 
respectively. 
We want to lead the recurrence formula of them. 
In computations below, we use the following equalities.
\[
{ n+1 \brack r}_q={ n \brack r-1}_q+q^r { n \brack r}_q \, ,\quad
p_{n+1}=p_n + q^{-n-1} p_n\, .
\]
Change variables $r$ to $k-r$ in the formula of
 $a_k$ and $b_k$ and we get 
\begin{align*}
a_k
=&\,
q^k
\sum_{r=0}^k {k \brack r}_q 
(-1)^{r}q^{r(r-k-1)}p_{r+1} p_{k-r}\, ,\\
b_k
=&\,
q^k
\sum_{r=0}^k {k \brack r}_q 
(-1)^{r}q^{r(r-k-2)}p_{r} p_{k-r+1}\, .
\end{align*}
For $a_k$ we compute as follows. 
\begin{align*}
q^{-k-1}a_{k+1}
=&
\sum_{r=0}^{k+1} 
{k+1 \brack r}_q (-1)^{r}q^{r(r-k-2)}p_{r+1} p_{k-r+1}\\
=&
\sum_{r=0}^{k+1} 
\bigg{(}{ k \brack r-1}_q+q^r { k \brack r}_q \bigg{)}
(-1)^{r}q^{r(r-k-2)}p_{r+1} p_{k-r+1} \\
=&
\sum_{r=0}^{k} \!
{ k\brack r}_q \!(-1)^{r+1}q^{(r+1)(r-\!k-\!1)}p_{r+2} p_{k-r} 
\!+\!
\sum_{r=0}^{k} \!
{k \brack r}_q \!(-1)^{r}q^{r(r-\!k-\!1)}p_{r+1} p_{k-r+1}\\
=&
-q^{-k-1}
\bigg{(}
\sum_{r=0}^{k} 
{ k\brack r}_q (-1)^{r+1}q^{r(r-k)}(1+q^{-r-2})p_{r+1} p_{k-r}
\bigg{)}\\
&\quad+
\sum_{r=0}^{k} 
{k \brack r}_q (-1)^{r}q^{r(r-k-1)}(1+q^{-k+r-1})p_{r+1} p_{k-r}\\
=&
(1-q^{-k-3})
\sum_{r=0}^{k} 
{k \brack r}_q (-1)^{r}q^{r(r-k-1)}p_{r+1} p_{k-r}\\
=&
q^{-k}(1-q^{-k-3})a_k\, .
\end{align*}
Hence we get $a_{k+1}=(q-q^{-k-2})a_k$. 
Similarly we compute $b_{k+1}$ as follows.
\begin{align*}
q^{-k-1}b_{k+1}
=&
\sum_{r=0}^{k+1} {k+1 \brack r}_q 
(-1)^{r}q^{r(r-k-3)}p_{r} p_{k-r+2}\\
=&
\sum_{r=0}^{k+1} 
\bigg{(}{ k \brack r-1}_q+q^r { k \brack r}_q \bigg{)}
(-1)^{r}q^{r(r-k-3)}p_{r} p_{k-r+2}\\
=&
\sum_{r=0}^{k}\hspace{-1mm}{k \brack r}_q \!
(-1)^{r+1}q^{(r+1)(r-\!k-\!2)}p_{r+1} p_{k-r+1}
\!+\!
\sum_{r=0}^{k}\hspace{-1mm}{k \brack r}_q \!
(-1)^{r}q^{r(r-\!k-\!2)}p_{r} p_{k-r+2}\\
=&
-q^{-k-2}
\bigg{(}
\sum_{r=0}^{k} {k \brack r}_q 
(-1)^{r}q^{r(r-k-1)}(p_{r}+q^{-r-1}p_r) p_{k-r+1}
\bigg{)}\\
&\quad
+
\sum_{r=0}^{k} {k \brack r}_q 
(-1)^{r}q^{r(r-k-2)}(1+q^{-k+r-2})p_{r} p_{k-r+1}\\
=&
(-q^{-k-3}+1)
\sum_{r=0}^{k} {k \brack r}_q 
(-1)^{r}q^{r(r-k-2)}p_{r} p_{k-r+1}\\
=&
q^{-k}(1-q^{-k-3})b_k\, .
\end{align*}
Hence we get $b_{k+1}=(q-q^{-k-2})b_k$, which is the same form 
as $a_{k+1}$. 
We can easily check $a_0=p_1=b_0$ and it deduces the desired 
equality. 
\end{proof}

At last we prove the existence of a right coideal 
$A$ of type $D_1$. 
It has a spectral pattern, 
$\oplus_{k\in\Z_{\geq0}}\big{(}k+\frac{1+(-1)^k}{2}\big{)}\pi_k
\oplus \oplus_{k\in\Z_{\geq0}}(k+1)\pi_{k+\frac{1}{2}}$. 
If $A$ really exists, 
its spectral subspace $A_{\pi_{\frac{1}{2}}}$ generates 
the whole $C^*$-algebra $A$. 
In fact the generated $C^*$-algebra contains $\pi_{\frac{1}{2}}$ 
and the only $D_1$ case admits such an algebra. 
We have to clarify a self-conjugate 
$\pi_{\frac{1}{2}}$-eigenvector 
$\eta=\sum_{s\in I_\frac{1}{2}} d_s \bw_s^{\frac{1}{2}}$. 
By the self-conjugacy we get 
$d_{-\frac{1}{2}}=(-q)^{\frac{1}{2}}\overline{d_{\frac{1}{2}}}$. 
We may assume $d_{\frac{1}{2}}$ is real number by applying $\beta_z^L$ 
for some $z\in\T$. 
Uniqueness up to conjugation follows from this observation. 
Hence we conclude $\eta=(\sqrt{q_0}x+v, \sqrt{q_0}u+y)$ is 
a desired $\pi_\frac{1}{2}$-vector. 
Let us write $a=\sqrt{q_0}x+v$ and $b=\sqrt{q_0}u+y$. 
We study the $C^*$-algebra $A$, which is 
generated by $a$ and $b$, 
and derive the result that it is really of type $D_1$. 
By direct calculation we obtain
\[
a^*=\sqrt{q_0}\,b\, ,\quad
\sqrt{q_0}\,ab+\sqrt{q_0}^{-1}ba=1+q_0\, .
\]
Because of this equality the smooth part of $A$ is linearly 
spanned by $a^kb^\ell$ for $k,\ell\geq0$. 
We shall show $P_{\frac{1}{2}}(a^k b^\ell)\in \C a+ \C b$. 
If it is done, then it shows the $\pi_{\frac{1}{2}}$-multiplicity 
is exactly one. 
Since other types which have a single loop at a vertex 
do not occur, we can conclude $A$ is of type $D_1$. 
Now we start a proof. 
Applying $\beta^R$ to $a$ and $b$, we get 
$\beta_z^R(a)=z a$ and $\beta_z^L(b)=\bar{z}b$ for $z\in \T$. 
Hence the element $P_{\frac{1}{2}}(a^k b^\ell)$ must have 
the following form.
\[
\left\{
\renewcommand{\arraycolsep}{1.5pt}
\begin{array}{rcrcrcc}
P_{\frac{1}{2}}(a^{\ell+1} b^\ell) &=&\lambda a \, ,&  \\
P_{\frac{1}{2}}(a^k b^{k+1}) &=& \mu b \, ,&\\
P_{\frac{1}{2}}(a^k b^\ell)     &=&0\ & \mbox{if}\  |k-\ell|\neq1 \, ,\\
\end{array}
\right.
\]
where $\lambda$ and $\mu$ are complex number. 
If we show $a^k b^{k+1}=\mu b+\cdots$, 
then we have $a^{k+1}b^k=\bar{\mu}a+\cdots$. 
Hence it suffices to show 
$P_{\frac{1}{2}}(a^k b^{k+1})=\mu b$. 
Recall the following well-known binomial equality.
\[
(c+d)^n=\sum_{r=0}^n {n \brack r}_q c^r d^{n-r}\quad
\mbox{if}\quad dc=qcd\, .
\]
Then we have
\[
a^k b^{k+1}=
\sum_{
\begin{subarray}{1}
0\leq r \leq k\\
0\leq s \leq k+1
\end{subarray}
}
{k\brack r}_q{k+1\brack s}_q
\sqrt{q_0}^{s+r}
x^rv^{k-r} u^s y^{k+1-s}\, .
\]
We want to take out the coefficient of $u$ and $y$. 
For the sake of this, we make use of $\beta^L$. 
Since we have 
$\beta_z^L(x^rv^{k-r}u^sy^{k+1-s})=z^{2r+2s-2k-1}x^rv^{k-r}u^sy^{k+1-s}$, 
we obtain 
\begin{align}
P_{\frac{1}{2}}(a^k b^{k+1})=&\,
\sum_{r=0}^k
{k\brack r}_q{k+1\brack k+1-r}_q
\sqrt{q_0}^{k+1}
P_{\frac{1}{2}}(x^rv^{k-r} u^{k+1-r} y^{r})\notag\\
&\,
+
\sum_{r=0}^k
{k\brack r}_q{k+1\brack k-r}_q
\sqrt{q_0}^{k}
P_{\frac{1}{2}}(x^rv^{k-r} u^{k-r} y^{r+1})\notag\\
=&\,
\sum_{r=0}^k
{k\brack r}_q{k+1\brack r}_q
\sqrt{q_0}^{k+1}q^{-r(-2r+2k+1)}
P_{\frac{1}{2}}(x^r y^{r}u^{k-r}v^{k-r} u)\notag\\
&\,
+
\sum_{r=0}^k
{k\brack r}_q{k+1\brack k-r}_q
\sqrt{q_0}^{k}q^{-r(-2r+2k)}
P_{\frac{1}{2}}(x^r y^{r}u^{k-r} v^{k-r} y)\, \notag\\
=&\,
\frac{\sqrt{q_0}^{k\!+\!1}(2)_q}{(k\!+\!2)_q}
\!\sum_{r=0}^k\hspace{-0.2mm}
{k\brack r}_q{k\!+\!1\brack r}_q
\!(-1)^{r(k-r)} 
q^{r(r-k+1)}
\frac{(r)_q!(k\!-\!r\!+\!1)_q!}{(k\!+\!1)_q!}\,u
\label{A201}\\
&
\hspace{-3mm}
\!+\!
\frac{\sqrt{q_0}^{k}q_0^{-k}(2)_q}{(k\!+\!2)_q}
\!\sum_{r=0}^k\hspace{-0.2mm}
{k\brack r}_q{k\!+\!1\brack k\!-\!r}_q
\!(-1)^{r(k-r)+k}
q^{r(r-k+2)}
\frac{(r\!+\!1)_q!(k\!-\!r)_q!}{(k\!+\!1)_q!}
y,\notag
\end{align}
where we have used Lemma \ref{proj half} in the last equality. 
Moreover we use the following formula.
\[
\frac{(n)_q!}{(r)_q!(n-r)_q!}=
(-1)^{r(n-r)}
\frac{p_n}{p_r p_{n-r}}
{n \brack r}_q\, .
\]
Then (\ref{A201}) is equal to
\begin{align*}
(\ref{A201})
=&\,
\frac{(-1)^k\sqrt{q_0}^{k+1}
(2)_q}{(k+2)_q p_{k+1}}
\sum_{r=0}^k
{k\brack r}_q
(-1)^{k-r} q^{r(r-k+1)}
p_r p_{k-r+1}\,u
\\
&\,
+
\frac{(-1)^k\sqrt{q_0}^{k}q_0^{-k}(2)_q}{(k+2)_q!p_{k+1}}
\sum_{r=0}^k
{k\brack r}_q
(-1)^{k-r}
q^{r(r-k+2)}
p_{r+1}p_{k-r}
\,y\\
=&\,
\frac{(-1)^k\sqrt{q_0}^{k}
(2)_q}{(k+2)_q p_{k+1}}(1+q^{-1})(q^{-2};q^{-1})_k
(\sqrt{q_0}u+y)\, ,
\end{align*}
where we have used Lemma \ref{binomial}. 
Hence we have obtained 
$P_{\frac{1}{2}}(a^k b^{k+1})\in\C
(\sqrt{q_0}u+y)$ 
and this completes the proof of the existence of a right coideal 
of type $D_1$. 

(IV) $A_m'\, (3\leq m\leq \infty)$ case. 
We assume that there exists a right coideal of type $A_m'$. 
Although we only treat the finite $m$ case here, 
we can similarly derive a contradiction in the case of $m=\infty$. 
Look at the first vertex from the left of Figure \ref{A_m'}. 
Since the entry of the Perron-Frobenius eigenvector of this vertex is $1$,
the corresponding ergodic system becomes a right coideal 
of $\csuq$ by Lemma \ref{criterion on eigenvector}. 
Let us denote the right coideal by $A$. 
It is also of type $A_m'$ by Theorem \ref{reduced ergodic system}. 
The spectral subspace $A_{\pi_{\frac{1}{2}}}$ 
generates a right coideal of type $D_1$, which is denoted by $B$. 
Now we consider the subsystem $B\subset A$. 
Let $\Lambda$ be the inclusion matrix of 
$B\rtimes_\delta \suq \subset A\rtimes_\delta \suq$. 
Apply 
Proposition \ref{intertwining property} and 
we get the equality 
$\Lambda\,\mathbb{M}^B(\pi_{\frac{1}{2}}) 
=\mathbb{M}^A(\pi_{\frac{1}{2}})\, \Lambda$, 
where $\mathbb{M}^B$ and $\mathbb{M}^A$ are the multiplicity maps 
of $B$ and $A$, respectively. 
They have the following form.
\begin{center}
$\mathbb{M}^B(\pi_{\frac{1}{2}})
=
\begin{pmatrix}
1 & 1\\
1 & 1\\
\end{pmatrix}
,$\quad
$\mathbb{M}^A(\pi_{\frac{1}{2}})
=
\begin{pmatrix}
1      & 1       & 0        & 0        &  \dots  & 0      & 0 \\
1      & 0       & 1        & \ddots   &  \ddots & 0      & 0 \\

0      & 1       & 0        & \ddots   & \ddots  & \vdots & \vdots \\

0      & 0       & \ddots   & \ddots   & \ddots  & 0      & 0 \\
\vdots & \ddots  & \ddots   & \ddots   & \ddots  & 1      & 0 \\
0      & 0       & \dots    & 0        & 1       & 0      & 1 \\
0      & 0       & \dots    & 0        & 0       & 1      & 1 
\end{pmatrix}
.$

\end{center}
We know that the minimal projection $p_0\in C_r^*(\suq)$ 
is also minimal in $B\rtimes_\delta \suq $ and 
$A\rtimes_\delta \suq$. 
We discuss the corresponding vertex of $p_0$. 
About $B\rtimes_\delta \suq$ we may assume $p_0$ sits at 
the left vertex of Figure \ref{D_1}. 
About $A\rtimes_\delta \suq$ we may assume $p_0$ sits at 
the first vertex from the left in Figure \ref{A_m'} 
because the reduced system $p_0(A\otimes\mathbb{K}(L^2(\suq)))p_0$ 
is canonically isomorphic to $A$. 
Hence the inclusion matrix $\Lambda$ must be as the following form.
\[
\Lambda
=
\begin{pmatrix}
1 & \lambda_1 \\
0 & \lambda_2 \\
\vdots & \vdots \\
0 & \lambda_{m-1}\\
0 & \lambda_{m}
\end{pmatrix}
.
\]
Then we want to solve the equation 
$\Lambda\,\mathbb{M}^B(\pi_{\frac{1}{2}}) 
=\mathbb{M}^A(\pi_{\frac{1}{2}})\, \Lambda$ 
about the non-negative integers $\lambda_1\dots\lambda_m$, 
however, we immediately see it has no solutions. 
This is a contradiction and hence there is not 
a right coideal of type $A_m'$.

\section{Classification of right coideals of $\csum$}
In this final section we complete the classification program 
of right coideals associated to the quantum $SU(2)$ group. 
Its representation theory such as actions of $\umsu$ on $\csum$ or 
the Clebsh-Gordan coefficients is obtained by the limit of 
$q_0\rightarrow 1$. 
The continuous function algebra $\csum$ is 
generated by $x,u,v$ and $y$ which satisfy
\[
ux=-xu\,,\quad vx=-xv\, , \quad uy=-yu\, ,\quad vy=-yv\, ,\quad 
uv=vu\, ,
\]
\[
xy+uv=yx+uv=1\, ,\quad x^*=y\, ,\quad u^*=v\, .
\]
By simple calculation, we see $\kappa^2=\id$. 
Hence $\csum$ is a compact Kac algebra. 
Refer the theory of general Kac algebras 
to \cite{Enock Schwarz}. 
One of differences between cases of 
$q=-1$ and $q^2\neq1$ is 
amount of the quantum subgroups (or right coideals). 
In fact as we have seen 
there are not a lot of right coideals 
in $\csuq$, on the contrary, 
in \cite{Podles2} he completely collects plentiful quantum subgroups 
of $\summ$. 
The main result in this section is as follows. 
The definition of $\eta^{\frac{n}{2}}$ and 
$\wdh{\eta}^{\frac{n}{2}}$ is given in the case of type $\T_n$ with 
odd $n\geq3$. 

\begin{thm}
If a right coideal $A$ is not of type $\T_n$ (odd $n\geq3$) or 
$D_n$ (odd $n\geq1$), 
there exists a closed subgroup $H$ in $SO_1(3)$ such that 
$A$ is $C(H\setminus \som)$. 
If a right coideal $A$ is of type $\T_n$ (odd $n\geq3$), 
$A$ is conjugated to 
$C(\T_n\setminus \summ)$ or 
$C^*(\eta^{\frac{n}{2}},\wdh{\eta}^{\frac{n}{2}})$. 
If a right coideal $A$ is of type $D_1$, then 
$A$ is conjugated to $C(D_1\setminus \summ)$. 
If a right coideal $A$ is of type $D_n$ (odd $n\geq3$), 
$A$ is conjugated to 
$C(D_n\setminus \summ)$ or $C^*(\eta^{\frac{n}{2}})$. 
Here conjugation is given by $\beta_z^L$ for some $z\in\T$. 
\end{thm}
Although 
$C(\T_n\setminus \summ)$ is not isomorphic to 
$C^*(\eta^{\frac{n}{2}},\wdh{\eta}^{\frac{n}{2}})$ 
as a $\summ$-covariant system, 
in Proposition \ref{isomorphic} 
we show that $C(D_n\setminus \summ)$ is 
isomorphic to $C^*(\eta^{\frac{n}{2}})$ 
as a $\summ$-covariant system. 
As we have said in the previous section, 
we have the isomorphism between 
$\csom$ and $\cso$ as compact quantum groups. 
Hence we can conclude that right coideals are 
quotient by subgroups when 
they are one of type $A_4^*$, $S_4^*$, $A_5^*$, 
$D_m^*$ $(m\geq2)$, $D_\infty^*$, $\T_{2n}$ $(n\geq1)$. 
Hence we have to study other types: $\T_n$ (odd $n\geq3$), 
$A_m'$ ($m\geq3$), $D_n$ (odd $n\geq1$). 

(I) $A_m'\,(m\geq3)$ case. 
The entirely same proof as in the previous section 
derives a contradiction well. 

(II) $D_1$ case. 
If we consider the limit $q_0\rightarrow 1$ in the previous 
section, 
we can show its existence. 
In that procedure we do not need any lemma stated in that section, 
because $\pi_{\frac{1}{2}}$-eigenvector is uniquely determined and 
in 
\cite[Proposition 3.8]{Podles2} 
it has been shown that there exists a right coideal by subgroup 
$D_1$ 
whose $\pi_{\frac{1}{2}}$-multiplicity is one. 

(III) $\T_{n}$ (odd $n\geq3$) case. 
$A$ has the spectral pattern 
\[
\bigoplus_{k\in\Z_{\geq0}}\bigg{(}1+2\bigg{[}\frac{k}{n}\bigg{]}\bigg{)}\pi_k
\oplus 
\bigoplus_{k\in \Z_{\geq0}} 
2\bigg{[}\frac{2k+n+1}{2n}\bigg{]}\pi_{k+\frac{1}{2}}.
\] 
Since its integer part becomes a right coideal in $\csom$, 
it must be quotient by a subgroup $H\subset \som$. 
From the spectral pattern there exist angles $0\leq\chi<2\pi$ and 
$0\leq\psi<\pi$ such that $H=\T_{n}^{\chi,\psi}$. 
We begin to determine a $\pi_1$-eigenvector of 
$C(\T_{n}^{\chi,\psi}\setminus\som)$. 

\begin{lem}\label{pi1 of tn}
The following vector is a $\pi_1$-eigenvector 
for $C(\T_{n}^{\chi,\psi}\setminus\som)$. 
\[\xi^{\chi,\psi}
=\sqrt{2}^{\,-1}
(-ie^{i\chi}\sin\psi,\sqrt{2}\cos\psi,ie^{-i\chi}\sin\psi)w(\pi_1)
.
\]
\end{lem}
\begin{proof}
First we consider a $\pi_1$-eigenvector 
$\zeta^{\chi,\psi}=\sum_{t\in I_1}c_t \bw_t^1$ of 
$C(\T_{n}^{\chi,\psi}\setminus SO(3))$. 
The left action of $\pi_1(g)\in\T_{n}^{\chi,\psi}$ is given by 
multiplication of $w(\pi_1)(g)$ to $w(\pi_1)$ from the left. 
Hence we have to get a vector 
$(c_{-1},c_0,c_1)w(\pi_1)(g)=(c_{-1},c_0,c_1)$ for 
all $g\in SU(2)$ with $\pi_1(g)\in \T_{n}^{\chi,\psi}$. 
If $\chi$ and $\psi$ are equal to $0$, 
we can easily get $c_{\pm1}=0$ and $c_0=1$. 
Since we know 
$\T_{n}^{\chi,\psi}
=\Ad(\pi_1(r^{12}(\chi)r^{13}(\psi)))(\T_{n}^{0,0})$, 
we obtain 
$(c_{-1},c_0,c_1)
=
(0,1,0)
w(\pi_1)(r^{13}(-\psi)r^{12}(-\chi))$. 
Using two matrices 
\[w(\pi_1)(r^{13}(-\psi))
=
\begin{pmatrix}
\cos^2 \frac{\psi}{2} & -\sqrt{2}\sin\frac{\psi}{2}\cos\frac{\psi}{2} &
\sin^2 \frac{\psi}{2}\\
\sqrt{2}\sin\frac{\psi}{2}\cos\frac{\psi}{2} & 
1-2\sin^2\frac{\psi}{2} & 
-\sqrt{2}\sin\frac{\psi}{2}\cos\frac{\psi}{2}\\
\sin^2\psi & -\sqrt{2}\sin\frac{\psi}{2}\cos\frac{\psi}{2} & 
\cos^2 \frac{\psi}{2}
\end{pmatrix}
\]
and 
\[w(\pi_1)(r^{12}(-\chi))
=
\begin{pmatrix}
e^{i\chi}&0&0\\
0&1&0\\
0&0&e^{-i\chi}
\end{pmatrix}
,
\] 
we have $c_{\pm1}=\mp \sqrt{2}^{\,-1}
e^{\mp i\chi}\sin\psi$ and $c_0=\cos\psi$. 
Next we make use of $\Xi_1:C(SO_1(3))\longrightarrow C(SO_{-1}(3))$ 
and get the desired eigenvector. 
\end{proof}

Next we determine the $\pi_{\frac{n}{2}}$-eigenvector space which is 
two-dimensional. 
Let us write $m$ for $\frac{n}{2}$. 
For $\eta=\sum_{t\in I_m}d_t\bw_t^m$ 
we use $q=-1$ version of 
Lemma \ref{product of pi1 negative} and we obtain 
\begin{align*}
\sqrt{2m(2m-1)}\Psi_{m-1}(\xi^{\chi,\psi},\eta)
=
&\,
\sum_{t\in I_m}
\big{\{}
-ie^{i\chi}\sin\psi \sqrt{(m+t)(m+t+1)}d_{t+1}\\
&\qquad\,
-(-1)^{m-t}2\cos\psi\sqrt{(m-t)(m+t)}d_t\\
&\qquad\,
+ie^{-i\chi}\sin\psi \sqrt{(m-t)(m-t+1)}d_{t-1}
\big{\}}
\bw_t^{m}\,.
\end{align*}
Define $\alpha_t$, $\beta_t$ and $\gamma_t$ by
\begin{align*}
\alpha_t
=&\,
-ie^{i\chi}\sin\psi \sqrt{(m+t)(m+t+1)}\,,\\
\beta_t
=&\,
-(-1)^{m-t}2\cos\psi\sqrt{(m-t)(m+t)}\,,\\
\gamma_t
=&\,
ie^{-i\chi}\sin\psi \sqrt{(m-t)(m-t+1)}
\, .
\end{align*}
Then we have $\alpha_{-t}=\ovl{\gamma_t}$ 
and $\beta_{-t}=-\beta_t$. 
If $\eta$ is a $\pi_m$-eigenvector for $A$, 
then we have the following recurrence formula for $t\in I_m$.
\begin{equation}\label{tn01}
\alpha_t d_{t+1}
+\beta_t d_t
+\gamma_t d_{t-1}=0\,.
\end{equation}
If we assume $\eta$ is self-conjugate, we have 
\[
\alpha_t d_{t+1}
-\beta_t d_t
+\gamma_t d_{t-1}=0
\]
where we define $d_t=0$ for $|t|\geq m+1$. 
From (\ref{tn01}) we get $\beta_t d_t=0$ and 
$\alpha_t d_{t+1}+\gamma_t d_{t-1}=0$. 

(a) $\psi\neq 0,\frac{\pi}{2}$ case. 
Since $\beta_t$ is equal to zero if $|t|=m$, 
we have $d_t=0$ for $|t|\leq m-1$. 
We also have $\alpha_{m-1} d_{m}+\gamma_{m-1} d_{m-2}=0$. 
Then we get $d_m=0$ by $\alpha_{m-1}\neq0$ and $m-2\geq-m+1$. 
This shows $d_t$ are all equal to zero. This is not appropriate. 

(b) $\psi=0$ case. 
Then similarly we have $d_t=0$ for $|t|\leq m-1$. 
Since $\pi_m$-eigenvector space is two dimensional, 
it is spanned by $\bw_{\pm m}^m$. 
This shows $A=C(\T_n\setminus \summ)$. 

(c) $\psi=\frac{\pi}{2}$ case. 
We have only a non-trivial equation $\alpha_t d_{t+1}+\gamma_t d_{t-1}=0$. 
Its solution space is two dimensional. 
We give an explicit solution as follows. 
The proof is straightforward. 

\begin{lem}\label{pim-eigenvector}
Let $m$ be a half integer in $\frac{3}{2}+\Z_{\geq0}$ and 
$0\leq\chi<2\pi$. 
Consider the recurrence equation, 
\[
-ie^{i\chi} \sqrt{(m+t)(m+t+1)}d_{t+1}
+ie^{-i\chi}\sqrt{(m-t)(m-t+1)}d_{t-1}
=0
.
\]
Then its solution is a linear combination of 
$d_{t}=e^{i(m-t)\chi}{2m\brack m-t}^{\frac{1}{2}}$ and 
$d_t=(-1)^{m-t}e^{i(m-t)\chi}{2m\brack m-t}^{\frac{1}{2}}$. 
\end{lem}

Define two $\pi_m$-eigenvectors of $\csum$ $\eta_{\chi}^m$ and 
$\widehat{\eta}_{\chi}^m$ by
\begin{align*}\label{eta etahat}
\eta^{m,\chi}
=&\,
\sum_{t\in I_m}e^{i(m-t)\chi}{2m\brack m-t}^{\frac{1}{2}}\bw_t^m\,,
\\
\wdh{\eta}^{m,\chi}
=&\,
\sum_{t\in I_m}(-1)^{m-t}e^{i(m-t)\chi}
{2m\brack m-t}^{\frac{1}{2}}\bw_t^m\,.
\end{align*}
Since the $\pi_m$-eigenvector space of $A$ is two-dimensional, 
it is spanned by the above vectors. 
Conversely we consider the right coideal $B$ which is generated 
by $\eta^{m,\chi}$ and $\wdh{\eta}^{m,\chi}$. 
We want to conclude that $B$ is a right coideal of type $\T_n$ and 
hence $B$ coincides with $A$ as a result. 
For the sake of this, it suffices to show that 
the $\pi_{\ell}$-eigenvector space of $B$ is trivial for all 
$\ell\in\{\frac{1}{2},\ldots,m-1\}$. 
We prepare the element $g^{\theta,\chi}$ of $SU(2)$ defined by 
$g^{\theta,\chi}
=
\begin{pmatrix}
\cos\frac{\theta}{2} & -ie^{-i\chi}\sin\frac{\theta}{2}\\
-ie^{i\chi}\sin\frac{\theta}{2} & \cos\frac{\theta}{2}
\end{pmatrix}
.
$
Recall the $*$-homomorphism 
$\ups_{g^{\theta,\chi}}:\csum\longrightarrow \mathbb{B}(\C^2)$
\[
\begin{pmatrix}
\ups_{g^{\theta,\chi}}(x)&\ups_{g^{\theta,\chi}}(u)\\
\ups_{g^{\theta,\chi}}(v)&\ups_{g^{\theta,\chi}}(y)
\end{pmatrix}
=
\begin{pmatrix}
\cos\frac{\theta}{2}\sigma_1 & ie^{-i\chi}\sin\frac{\theta}{2}\sigma_2\\
-ie^{i\chi}\sin\frac{\theta}{2}\sigma_2 &\cos\frac{\theta}{2} \sigma_1
\end{pmatrix}
.
\]
Define a $*$-homomorphism 
$\rho_{g^{\theta,\chi}}=(\ups_{g^{\theta,\chi}}\otimes\id)\circ\delta$ 
and we have
\[
\begin{pmatrix}
\!\rho_{g^{\theta,\chi}}(x) & \hspace{-3mm}\rho_{g^{\theta,\chi}}(u) \!\\
\!\rho_{g^{\theta,\chi}}(v) & \hspace{-3mm}\rho_{g^{\theta,\chi}}(y)\!
\end{pmatrix}
\hspace{-1.5mm}=\hspace{-1.5mm}
\begin{pmatrix}
\hspace{0mm}\cos\frac{\theta}{2}\sigma_1\!\otimes \!x 
\!+\! ie^{-i\chi}\sin\frac{\theta}{2}\sigma_2\!\otimes \!v 
& \hspace{-2mm}
\cos\frac{\theta}{2}\sigma_1\!\otimes \!u
\!+\!
ie^{-i\chi}\sin\frac{\theta}{2}\sigma_2\!\otimes \!y 
\\
\hspace{0mm}
-\!ie^{i\chi}\sin\frac{\theta}{2} \sigma_2\!\otimes \!x
\!+\!\cos\frac{\theta}{2}\sigma_1\!\otimes \!v 
&\hspace{-2mm}
-\!ie^{i\chi}\sin\frac{\theta}{2} \sigma_2\!\otimes \!u
\!+\!\cos\frac{\theta}{2}\sigma_1\!\otimes \!y 
\end{pmatrix}
.
\]

\begin{lem}\label{rho and eta}
Let $m$ be a half integer in $\frac{1}{2}+\Z_{\geq0}$. 
For all $s\in I_m$ we have the following equalities.
\begin{align*}
&\rho_{g^{\theta,\chi}}((\eta^{m,\chi})_s)
=\hspace{-1.0mm}
\sum_{r\in I_m}
(\cos m\theta\sigma_1+(-1)^{m-r}i\sin m\theta\sigma_2)
\otimes e^{i(m-r)\chi}{2m\brack m-r}^{\frac{1}{2}}w(\pi_m)_{r,s},
\\
&\rho_{g^{\theta,\chi}}
((\wdh{\eta}^{m,\chi})_s)
=\hspace{-1.0mm}
\sum_{r\in I_m}
((-1)^{m-r}\cos m\theta\sigma_1
-i\sin m\theta \sigma_2)
\otimes
e^{i(m-r)\chi}
{2m\brack m\!-\!r}^{\frac{1}{2}}
w(\pi_m)_{r,s} 
.
\end{align*}
\end{lem}
\begin{proof}
Since $\rho_{g^{\theta,\chi}}$ commutes with the lowering operator 
$f$, it is enough to show those equalities with $s=-m$. 
Here it is proved by induction on $m$. 
We can easily show 
$\rho_{g^{\theta,\chi}}\circ\beta_z^L=(\id\otimes\beta_z^L)
\circ
\rho_{g_{\theta,0}}
$ 
where $z$ is equal to $e^{i\frac{\chi}{2}}$. 
Moreover we have 
$\eta^{m,\chi}_{-m}=e^{im\chi}\beta_z^L(\eta^{m,0}_{-m})$ and 
$\wdh{\eta}^{m,\chi}_{-m}=e^{im\chi}\beta_z^L(\wdh{\eta}^{m,0}_{-m})$. 
Hence it yields 
$\rho_{g^{\theta,\chi}}(\eta^{m,\chi}_{-m})
=
e^{im\chi}(\id\otimes\beta_z^L)
\circ
\rho_{g_{\theta,0}}(\eta^{m,0}_{-m})
$ 
and 
$\rho_{g^{\theta,\chi}}(\wdh{\eta}^{m,\chi}_{-m})
=
e^{im\chi}(\id\otimes\beta_z^L)
\circ
\rho_{g_{\theta,0}}(\wdh{\eta}^{m,0}_{-m})
$. 
Hence we may assume $\chi$ is equal to $0$. 
For simplicity of notations we write 
$g^{\theta}$, $\eta^m$ and $\wdh{\eta}^m$ for 
$g_{\theta,0}$, $\eta^{m,0}$ and $\wdh{\eta}^{m,0}$, respectively. 
When $m$ is equal to $\frac{1}{2}$, the desired equalities are 
easily obtained. 
We assume that the desired equalities hold for $m$. 
We can justify an equality about $\rho_{g^{\theta}}(\eta_{-(m+1)}^{m+1})$ 
as follows. 
Making use of the equality 
$\eta_{-(m+1)}^{m+1}=x^2\eta_{-m}^m+2x\eta_{-m}^m v+\eta_{-m}^m v^2$, 
the computation of $\rho_{g^{\theta}}(x^2\eta_{-m}^m)$, 
$\rho_{g^{\theta}}(2x\eta_{-m}^m v)$ and $\rho_{g^{\theta}}(\eta_{-m}^m v^2)$ 
is carried out in the following way.
\begin{align*}
\rho_{g^{\theta}}(x^2\eta_{-m}^m)
=
&\,
(\cos^2\frac{\theta}{2}\otimes x^2
- 2\sin\frac{\theta}{2}\cos\frac{\theta}{2}\sigma_3\otimes xv
-\sin^2\frac{\theta}{2}\otimes v^2)
\\
&\quad\cdot
\sum_{r\in I_m}
(\cos m\theta\sigma_1+(-1)^{m-r}i\sin m\theta\sigma_2)
\otimes {2m\brack m-r} x^{m-r}v^{m+r}\\
=&\,
\sum_{r\in I_{m+1}}
\bigg{(}
\cos^2\frac{\theta}{2}(\cos m\theta\sigma_1 
+ (-1)^{m+1-r}i\sin m\theta\sigma_2){2m \brack m-r-1}\\
&\qquad\hspace{0.5cm}
+2\sin\frac{\theta}{2}\cos\frac{\theta}{2}
(-\sin m\theta\sigma_1+(-1)^{m+1-r}i\cos m\theta\sigma_2)
{2m \brack m-r}
\\
&\qquad\hspace{0.8cm}
-\sin^2\frac{\theta}{2}(\cos m\theta\sigma_1 
+ (-1)^{m+1-r}i\sin m\theta\sigma_2){2m \brack m-r+1}
\bigg{)}\\
&\qquad\hspace{0.5cm}
\otimes x^{m+1-r}v^{m+1+r},
\end{align*}
\begin{align*}
\rho_{g^{\theta}}(x\eta_{-m}^m v)
=&
\sum_{r\in I_{m+1}}
\hspace{-1mm}
\bigg{(}
\sin\frac{\theta}{2}\cos\frac{\theta}{2}
(-\sin m\theta \sigma_1+(-1)^{m+r}i\cos m\theta\sigma_2){2m\brack m-r-1}
\\
&\hspace{1.3cm}
+(\cos^2\frac{\theta}{2}-\sin^2\frac{\theta}{2})
(\cos m\theta\sigma_1-(-1)^{m-r}i\sin m\theta\sigma_2){2m \brack m-r}
\\
&\hspace{1.0cm}
+\sin\frac{\theta}{2}\cos\frac{\theta}{2}
(-\sin m\theta \sigma_1-(-1)^{m-r}i\cos m\theta\sigma_2){2m\brack m-r+1}
\bigg{)}
\\
&
\qquad\hspace{0.5cm}
\otimes x^{m+1-r}v^{m+1+r}
\end{align*}
and 
\begin{align*}
\rho_{g^{\theta}}(\eta_{-m}^m v^2)
=&\,
\sum_{r\in I_{m+1}}
\bigg{(}
-\sin^2\frac{\theta}{2}(\cos m\theta\sigma_1
-(-1)^{m-r}i\sin m\theta\sigma_2){2m \brack m-r-1}\\
&\qquad\hspace{0.8cm}
+2\sin\frac{\theta}{2}\cos\frac{\theta}{2}
(-\sin m\theta\sigma_1+(-1)^{m+r}i\cos m\theta\sigma_2)
{2m \brack m-r}
\\
&
\qquad\hspace{0.8cm}
+\cos^2\frac{\theta}{2}(\cos m\theta\sigma_1
-(-1)^{m-r}i\sin m\theta\sigma_2){2m \brack m-r+1}
\bigg{)}
\\
&
\qquad\hspace{0.5cm}
\otimes x^{m+1-r}v^{m+1+r}. 
\end{align*}
Then we can derive 
$\rho_{g_{\theta}}(\eta_{-(m+1)}^{m+1})
=
\sum_{r\in I_{m+1}}
(\cos(m+1)\theta+(-1)^{m+1-r}i\sin(m+1)\theta)\otimes 
{2m+2\brack m+1-r}x^{m+1-r}v^{m+1+r}
$ where 
we use the formula 
${2m+2\brack m+1-r}={2m\brack m-r-1}+2{2m\brack m-r}+{2m\brack m-r+1}$. 
By induction the assertion about $\rho_{g^{\theta}}(\eta_{-m}^m)$ is 
justified. 
About $\rho_{g^{\theta}}(\wdh{\eta}_{-m}^m)$ we make use of 
$\rho_{g^{\theta}}\circ\beta_i^L
=(\id\otimes\beta_i^L)\circ\rho_{g^{-\theta}}$ 
and 
$\wdh{\eta}_{-m}^m=i^{2m}\beta_i^L(\eta_{-m}^m)$. 
\end{proof}
Especially putting $\theta=\frac{\pi}{m}k$ $(k=0,\ldots,n-1)$ and 
$z=e^{-i\frac{\chi}{2}}$, 
we get
\begin{align*}
&
\rho_{g_{\frac{\pi}{m}k,\chi}}(\eta_s^{m,\chi})
=
(-1)^k\sigma_1\otimes \eta_s^{m,\chi},\\
&
\rho_{g_{\frac{\pi}{m}k,\chi}}(\wdh{\eta}_s^{m,\chi})
=
(-1)^k\sigma_1\otimes \wdh{\eta}_s^{m,\chi}\,.
\end{align*}
Therefore $B$ is contained in a $C^*$-algebra 
$C^{n,\chi}=\{a\in \csum\mid \rho_{g_{\frac{\pi}{m}k,\chi}}(a)\in
(\C+\C\sigma_1)\otimes \C a \mbox{\ for all\ }0\leq k<n\}
$. 
Actually, $C^{n,\chi}$ is a right coideal because of 
$(\rho_{g_{\frac{\pi}{m}k,\chi}}\otimes\id)\circ\delta
=(\id\otimes\delta)\circ\rho_{g_{\frac{\pi}{m}k,\chi}}$. 

\begin{lem}\label{pi-product eta etahat}
The following equalities hold.
\begin{align*}
\Psi_{1}(\eta^{m,\chi},\wdh{\eta}^{m,\chi})
=&\,
-i\sqrt{m}2^{2m}e^{2im\chi}\xi^{\chi,\frac{\pi}{2}},
\\ 
\Psi_{1}(\wdh{\eta}^{m,\chi},\eta^{m,\chi})
=&\,
-i\sqrt{m}2^{2m}e^{2im\chi}\xi^{\chi,\frac{\pi}{2}},
\\
\Psi_{1}(\eta^{m,\chi},\eta^{m,\chi})
=&\,
\Psi_{1}(\wdh{\eta}^{m,\chi},\wdh{\eta}^{m,\chi})
=
0.
\end{align*}
\end{lem}
\begin{proof}
We give a proof for only the first one. 
Others are proved similarly. 
By definition we have 
$\Psi_{1}(\eta^{m,\chi},\wdh{\eta}^{m,\chi})
=
\sum_{r=0}^{2m-1}(C_{m,m}^1)_r\eta^{m,\chi}_{m-r-1}
\wdh{\eta}^{m,\chi}_{-m+r}
$ 
where $(C_{m,m}^1)_r$ is equal to $(-1)^r\sqrt{2m}^{-1}\sqrt{(r+1)(2m-r)}$. 
Define complex numbers $c_{\pm1}$ and $c_0$ by 
$\Psi_{1}(\eta^{m,\chi},\wdh{\eta}^{m,\chi})
=
c_{-1}\bw_{-1}^1+c_0\bw_{0}^1+c_1\bw_{1}^1
$. 
In order to obtain their values we use the restriction homomorphism 
$\pi_\T:\csum\longrightarrow C(\T)$. 
In fact we get 
$\pi_\T(\Psi_{1}(\eta^{m,\chi},\wdh{\eta}^{m,\chi})_s)
=
z^{2s}c_s
$ for all $z\in \T$ and $1\leq s\leq 1$. 
First we have 
\begin{equation}\label{eta-1}
\Psi_{1}(\eta^{m,\chi},\wdh{\eta}^{m,\chi})_{-1}
=
\sum_{r=0}^{2m-1}
\sqrt{2m}^{\,-1}(-1)^r\sqrt{(r+1)(2m-r)}
(\eta^{m,\chi})_{m-1-r}\wdh{\eta}^{m,\chi}_{-m+r}.
\end{equation}
Applying $\pi_\T$ to the both side, we have 
\begin{align*}
c_{-1}
=
&\,
\sum_{r=0}^{2m-1}
\sqrt{2m}^{\,-1}(-1)^r\sqrt{(r+1)(2m-r)}
{2m \brack r+1}^{\frac{1}{2}}{2m \brack 2m-r}^{\frac{1}{2}}\\
&\qquad\cdot
(-1)^{2m-r}
e^{i(r+1)\chi}e^{i(2m-r)\chi}
\\
=
&\,
(-1)^{2m}\sqrt{2m}^{\,-1}
e^{i(2m+1)\chi}\sum_{r=0}^{2m-1}
\frac{2m!}{r!(2m-1-r)!}\\
=&\,
-2^{2m-1}\sqrt{2m}e^{i(2m+1)\chi}.
\end{align*}
If we apply the lowering operator $f$ to (\ref{eta-1}), 
then we can easily obtain 
\begin{equation}\label{eta0}
\Psi_{1}(\eta^{m,\chi},\wdh{\eta}^{m,\chi})_{0}
=
-\sqrt{m}^{\,-1}\sum_{r=0}^{2m}
(m-r)(\eta^{m,\chi})_{m-r}\wdh{\eta}^{m,\chi}_{-m+r}, 
\end{equation}
where we use the fact of 
$f\cdot \xi_p=i^{-2\nu+1}\sqrt{(\nu-p)(\nu+p+1)}\xi_{p+1}$ 
for a $\pi_\nu$-eigenvector $\xi$. 
Applying $\pi_\T$ to (\ref{eta0}), we obtain 
$c_0
=
-\sqrt{m}^{\,-1}\sum_{r=0}^{2m}
(m-r){2m\brack r}(-1)^{2m-r}
=0. 
$
Since we also have 
\begin{equation*}
\Psi_{1}(\eta^{m,\chi},\wdh{\eta}^{m,\chi})_{1}
=
\sum_{r=0}^{2m-1}
\sqrt{2m}^{\,-1}(-1)^r\sqrt{(r+1)(2m-r)}
(\eta^{m,\chi})_{m-r}\wdh{\eta}^{m,\chi}_{-m+r+1},
\end{equation*} 
similarly we can derive 
$c_1=2^{2m-1}\sqrt{2m}e^{i(2m-1)\chi}$. 
\end{proof}

Now we start to prove that $B$ is a right coideal of type $\T_n$. 
Let $n'$ be a smallest odd integer such that the $\pi_{\frac{n'}{2}}$-
eigenvector space of $B$ is not trivial. 
Put $m'=\frac{n'}{2}$. 
We claim $m'\geq \frac{3}{2}$. 
This is because 
there exist no non-zero 
complex numbers $c_{-\frac{1}{2}}$ and $c_{\frac{1}{2}}$ 
such that $c_{-\frac{1}{2}}x+c_{\frac{1}{2}}v$ is an element of 
$C^{m,\chi}$. 
From the previous lemma, $B$ contains $\pi_1$-part and 
it shows $B$ is not of type $D_{\ell}$ for some odd $\ell\geq3$ but 
of $\T_{n'}$. 
Again by Lemma \ref{pim-eigenvector} for $m'$, 
the $\pi_{m'}$-eigenvector space is spanned by 
$\eta^{m',\chi}$ and $\wdh{\eta}^{m',\chi}$. 
In particular, $\eta_{-m'}^{m',\chi}$ is an element of $C^{n,\chi}$. 
Since we have for $0\leq k\leq n-1$ 
\begin{align*}
\rho_{g_{\frac{2\pi}{n}k,\chi}}(\eta_{-m'}^{m',\chi})
=&\,
\sum_{r\in I_{m'}}(\cos \frac{m'}{m}\pi k\, \sigma_1
-(-1)^{m'-r}i\sin\frac{m'}{m}\pi k\,\sigma_2)\\
&\qquad\hspace{0.5cm}\otimes
e^{i(m'-r)\chi}{2m' \brack m'-r}^{\frac{1}{2}}w(\pi_{m'})_{r,-m'}
\end{align*} 
by Lemma \ref{rho and eta}, $\frac{m'}{m}$ must be equal to $1$. 
Hence the $\pi_{\ell}$-eigenvector space is trivial 
for $\ell\in\{\frac{1}{2}, \ldots, m-1\}$. 
This case is only admitted as type $\T_n$. 
As a result we also see that $B$ coincides with $C^{n,\chi}$. 

We have proved that 
a right coideal of type $\T_n$ (odd $n\geq3$) 
arises from the cyclic group in the maximal torus $\T$ or 
the cyclic group $\T_{n}^{\frac{\pi}{2},\chi}$. 
We have to check two right coideals 
$C(D_n\setminus \summ)$ and $C^{n,\chi}$ are not 
$\summ$-isomorphic. 
It suffices to prove it in the case of $\chi=0$ by 
considering $\beta^L$. 
Assume that there exists an $\summ$-isomorphism 
$\vartheta:C(D_n\setminus \summ) \longrightarrow C^{n,0}$. 
It induces the map between eigenvector spaces defined by 
$\vartheta(\xi)=(\vartheta(\xi_r))_{r\in I_\mu}$ for 
$\pi_\mu$-eigenvector $\xi$. 
On the eigenvector spaces, $\vartheta$ preserves 
the inner product, the conjugation operation $T$ 
and moreover the product map $\Psi$. 
Take complex numbers $\lambda$ and $\mu$ which satisfy 
$\vartheta(\bw_{-m}^m)=\lambda \eta^{m}+\mu \wdh{\eta}^{m}$. 
By Lemma \ref{pi-product eta etahat} and 
$\Psi_1(\bw_{-m}^m,\bw_{-m}^m)=0$, 
we have 
$-i\lambda\mu\sqrt{m}2^{2m+1}\xi^{0,\frac{\pi}{2}}=0$. 
Hence $\lambda=0$ or $\mu=0$. 
By using $\beta^L$ we may assume $\mu=0$. 
Now we have $\vartheta(\bw_{-m}^m)=\lambda\eta^m$. 
If we apply the conjugation $T$ to the both side, we obtain 
$\vartheta(\bw_m^m)=\ovl{\lambda}\eta^m$. 
This shows the range of $\vartheta$ onto $\pi_m$-eigenvector space is 
one-dimensional, however, this is contradiction. 
Therefore $C(\T_n\setminus \summ)$ and $C^{n,0}$ are not 
$\summ$-isomorphic. 
Finally we will make a discussion on characterizing 
a right coideal $C^{n,\chi}$ by a quantum subgroup of $\summ$. 
Let us recall a closed subgroup of $\summ$ which 
corresponds to $\T_n^{\frac{\pi}{2},0}$, 
that is, we consider the compact set 
$Z=\pi_1^{-1}(\T_n^{\frac{\pi}{2},0})$ 
and define $C(G_Z)=\pi_Z(\csum)$. 
We denote a subgroup $G_Z$ by $(\T_n^{\frac{\pi}{2},0})_\#$. 
The set $Z$ is actually a binary subgroup $(\T_n^{\frac{\pi}{2},0})^*$ 
and consists of $\{g^{\frac{2\pi}{n}k}\mid 0\leq k \leq 2n-1\}$. 
Hence $C((\T_n^{\frac{\pi}{2},0})_\#)\subset 
\oplus_{0\leq k\leq 2n-1} \mathbb{M}_{s(k)}(\C)$ 
where $s(k)=1$ if $k=0,\ n$ and $s(k)=2$ otherwise. 
We denote the coproduct of $C((\T_n^{\frac{\pi}{2},0})_\#)$ 
by $\delta_{Z}$. 
Now we look at the equality 
$\rho_{g^{\frac{2\pi}{n}k}}(\eta_s^m)=(-1)^k \sigma_1\otimes \eta_s^m$ 
for $0\leq k\leq 2n-1$ and $s\in I_m$. 
Define a self-adjoint unitary in $C((\T_n^{\frac{\pi}{2},0})_\#)$ 
by 
$w=\oplus_{k\neq 0,n} (-1)^k \sigma_1 
    \oplus 1 \oplus -1$. 
Then the above equality is equivalent to 
$(\pi_Z\otimes\id)\circ\delta(\eta_s^m)=w\otimes\eta_s^m$. 
Moreover, we have
\begin{align*}
\delta_{Z}(w)\otimes \eta_s^m
=&\,
(\delta_{Z}\circ\pi_Z\otimes\id)\circ\delta(\eta_s^m)\\
=&\,
(\pi_Z\otimes\pi_Z\otimes\id)\circ(\delta\otimes\id)\delta(\eta_s^m)\\
=&\,
(\pi_Z\otimes\pi_Z\otimes\id)\circ(\id\otimes\delta)\delta(\eta_s^m)\\
=&\,
w\otimes w\otimes \eta_s^m.
\end{align*}
Hence we have proved the following proposition. 

\begin{prop}\label{group like tn}
A self-adjoint unitary $w$ is a group-like element in 
$C((\T_n^{\frac{\pi}{2},0})_\#)$. 
\end{prop}

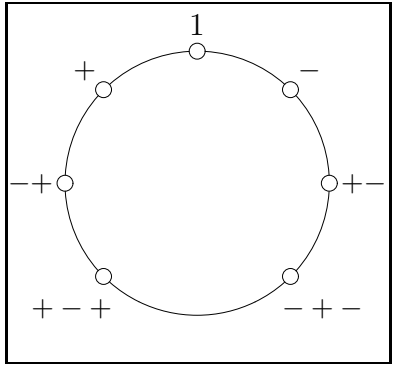
\begin{figure}
\begin{center}
\fbox{
\begin{picture}(131,130)

\put(115,65){\circle{6}}
\put(100.36,100.36){\circle{6}}
\put(65,115){\circle{6}}
\put(29.64,100.36){\circle{6}}
\put(15,65){\circle{6}}
\put(29.64,29.64){\circle{6}}
\put(100.36,29.64){\circle{6}}

\put(128,65){\makebox(0,0){{\small$ +-$}}}
\put(107.43,107.43){\makebox(0,0){{\small$-$}}}
\put(65,125){\makebox(0,0){$1$}}
\put(22.57,107.43){\makebox(0,0){{\small$+$}}}
\put(2,65){\makebox(0,0){{\small$-+$}}}
\put(17.57,17.57){\makebox(0,0){{\small$+-+$}}}
\put(112.43,17.57){\makebox(0,0){{\small$-+-$}}}

\put(65,65){\arc{100}{5.56}{-0.061}}

\put(65,65){\arc{100}{4.774}{-0.846}}

\put(65,65){\arc{100}{3.989}{-1.632}}

\put(65,65){\arc{100}{3.204}{-2.417}}

\put(65,65){\arc{100}{2.419}{-3.202}}


\put(65,65){\arc{100}{2.027}{-3.988}}


\put(65,65){\arc{100}{0.8486}{-5.166}}

\put(65,65){\arc{100}{0.0631}{-5.558}} 

\put(65,65){\arc{100}{1.117}{-4.256}}  


\end{picture}
}
\caption{The McKay diagram of $(\T_n^{0,\frac{\pi}{2}})_\#$, $2n$ nodes}
\label{Tns}
\end{center}
\end{figure}

A $C^*$-algebra generated by $w$ becomes a Hopf $*$-algebra 
which is isomorphic to $C(\Z_2)$. 
With this identification we have 
$C^{n,0}=\{a\in \csum \mid 
(\pi_Z\otimes\id)\circ\delta(a)\in C(\Z_2)\otimes \C a\}$. 
For general $C^{n,\chi}$ we also get similar results 
by using $\beta_z^L$. 
We study all the irreducible representations of 
$(\T_n^{\frac{\pi}{2},0})_\#$ as follows. 
Let us denote the restriction of the fundamental representation to 
$(\T_n^{\frac{\pi}{2},0})_\#$ simply by 
$w(\pi_{\frac{1}{2}})|$. 
Then we have 
$
w(\pi_{\frac{1}{2}})|
=
\big{(}
\begin{smallmatrix}
\pi_Z(x) & \pi_Z(u)\\
\pi_Z(v) & \pi_Z(y)
\end{smallmatrix}
\big{)}
$
 where we have 
\begin{align*}
\pi_Z(x)
=&\,
 \oplus_{k\neq 0,n}\cos\frac{\pi}{n}k \,\sigma_1 
 \oplus 1 \oplus -1,
\\
\pi_Z(u)
=&\,
 \oplus_{k\neq 0,n}i\sin\frac{\pi}{n}k \,\sigma_2 
 \oplus 0 \oplus 0,
\\
\pi_Z(v)
=&\,
 \oplus_{k\neq 0,n}-i\sin\frac{\pi}{n}k \,\sigma_2 
 \oplus 0 \oplus 0,
\\
\pi_Z(y)
=&\,
 \oplus_{k\neq 0,n}\cos\frac{\pi}{n}k \,\sigma_1 
 \oplus 1 \oplus -1.
\end{align*}
The unitary representation $w(\pi_{\frac{1}{2}})|$ invariantly 
acts on two lines 
$
\C\xi_{\pm}
=
\C
\big{(}
\begin{smallmatrix}
1\\
\pm i
\end{smallmatrix}
\big{)}
$ and their one dimensional actions are given by $w_{+}$ and $w_{-}$ 
with 
\[
w_\pm
=
\oplus_{k\neq 0,n}
\begin{pmatrix}
0 & e^{\pm i\frac{\pi}{n}k}\\
e^{\mp i\frac{\pi}{n}k} & 0
\end{pmatrix}
\oplus 1 \oplus -1
.
\]
Hence we have $w(\pi_{\frac{1}{2}})|=w_{+}\oplus w_{-}$. 
Notice that $w_{\pm}$ are the  self-conjugate unitary representations. 
This shows the difference from the fusion rules of the ordinary 
cyclic group $\T_{2n}$. 
And define the tensor product representations 
$w_{+-}=w_{+}w_{-}$ and $w_{-+}=w_{-}w_{+}$. 
Then they are self-conjugate and non-equivalent to each other, 
where non-equivalence comes from the non-commutativity of 
the underlying space $(\T_n^{\frac{\pi}{2},0})_\#$. 
If we want to write the McKay diagram, we compute 
$w(\pi_{\frac{1}{2}})|\cdot w_{+}=1\oplus w_{-+}$ and 
bond the vertices $+$ and $-+$ by a single line. 
We continue these procedures to get all the irreducible representations. 
As a result, they are represented by the reduced words 
($+-+-$ etc.) whose lengths are less than or equal to $n$ 
(Figure \ref{Tns}). 
For length $n$ word, we can easily check 
$w=+-\cdots-+$=$-+\cdots+-$. 
Hence $w$ sits at the bottom vertex. 
The fusion rules are given by 
this equality, $+^2=1$, $-^2=1$ and $\ovl{(\pm)}=\pm$. 
In the above proposition, we get a self-adjoint unitary $w$. 
It is the $n$-th irreducible representation $+-\cdots-+$. 
Let us denote $g=+-$ and $s=+$. 
Then we have $s^2=1$, $g^n=1$ and $sgs=g^{-1}$. 
This shows that the dual group $\widehat{(\T_n^{\frac{\pi}{2},0})_\#}$ 
is isomorphic to $D_n$, and 
hence $C((\T_n^{\frac{\pi}{2},0})_\#) 
\cong C_r^*(D_n)$ as the Hopf-algebras. 

(IV) $D_n$ (odd $n\geq3$) case. 
$A$ has the spectral pattern
\[\bigoplus_{k\in \Z_{\geq0}}
\bigg{(}\bigg{[}\frac{k}{n}\bigg{]}+\frac{1+(-1)^k}{2}\bigg{)}\pi_k
\oplus
\bigoplus_{k\in \Z_{\geq0}}
\bigg{(}\bigg{[}\frac{2k+1}{n}\bigg{]}
-\bigg{[}\frac{k}{n}\bigg{]}\bigg{)}
\pi_{k+\frac{1}{2}}\, . 
\]
Let $B$ be a right coideal of the integer parts of $A$. 
Since $B$ is a right coideal of $\csom$, 
it must be a quotient by subgroup of $\som$ 
for some angles 
$0\leq\phi<\frac{2\pi}{n}$, $0\leq\chi<2\pi$ and $-\pi\leq\psi\leq\pi$. 
Look at the spectral pattern for $B$ and 
we see that $B= C(D_n^{\phi,\chi,\psi} \setminus \som)$ 
We have to clarify its self-conjugate $\pi_2$-eigenvector. 

\begin{lem}
The following vector is a $\pi_2$-eigenvector for 
$C(D_n^{\phi,\chi,\psi}\setminus \summ)$.
\begin{align*}
\zeta^{\phi,\chi,\psi}
=
\Big{(}e^{i2\chi}\frac{\sqrt{6}}{4}\sin^2\psi,
&\,
ie^{i\chi}\frac{\sqrt{6}}{2}\sin\psi \cos\psi, 
1-\frac{6}{4}\sin^2\psi,\\
&
-ie^{-i\chi}\frac{\sqrt{6}}{2}\sin\psi \cos\psi, 
e^{-i2\chi}\frac{\sqrt{6}}{4}\sin^2\psi\Big{)}w(\pi_2)\, .
\end{align*}
\end{lem}
\begin{proof}
The proof is done as Lemma \ref{pi1 of tn} with considering 
$w(\pi_2)$ and $D_{n,\phi,\chi,\psi}
=\Ad(\pi_1(r^{12}(\chi)r^{13}(\psi)r^{12}(\phi)))(D_{n}^{0,0,0})$. 
\end{proof}

Next we use the following lemma in order to get 
a $\pi_{\frac{n}{2}}$-eigenvector of $A$, 
which is proved by 
$q_0\rightarrow1$ in 
Lemma \ref{product of pi2 negative} and an elementary calculation 
for $\bw_{\pm1}^2$. 

\begin{lem}

Let $m$ be a half integer with $m \in \frac{1}{2}+\Z_{\geq0}$. 
Then we have the following equalities for all $t\in I_m$. 

\begin{align*}
&\Psi_{m-1}(\bw_{-2}^2,\bw_t^m)
=
(-1)^{m-t}\frac{1}{2m\!-\!2}
\sqrt{
\frac{4(m\!-\!t\!+\!1)(m\!+\!t)(m\!+\!t\!-\!1)(m\!+\!t\!-\!2)}
{2m(2m\!-\!1)}
}\,\bw_{t-2}^{m-1}
 \!,
\\
&\Psi_{m-1}(\bw_{-1}^2,\bw_t^m)
=
\frac{2(m-2t+2)}{2m-2}\sqrt{\frac{(m+t)(m+t-1)}{2m(2m-1)}}
\bw_{t-1}^{m-1}
 ,
\\
&\Psi_{m-1}(\bw_{0}^2,\bw_t^m)
=
-(-1)^{m-t} 
\frac{2t}{(2m-2)}
\sqrt{
\frac{3!(m+t)(m-t)}
{2m(2m-1)}
}\,\bw_{t}^{m-1}
 ,
\\
&\Psi_{m-1}(\bw_{1}^2,\bw_t^m)
=
-\frac{2(m+2t+2)}{2m-2}\sqrt{\frac{(m-t)(m-t-1)}{2m(2m-1)}}
\bw_{t+1}^{m-1}
,
\\
&\Psi_{m-1}(\bw_{2}^2,\bw_t^m)
=\!
-(
\hspace{-0.1mm}-1)^{m-t}\frac{1}{2m\!-\!2}
\sqrt{
\frac{4(m\!+\!t\!+\!1)(m\!-\!t)(m\!-\!t\!-\!1)(m\!-\!t\!-\!2)}
{2m(2m\!-\!1)}
}\,\bw_{t+2}^{m-1}
\!.
\end{align*}
\end{lem}

Let us define $\alpha_t,\beta_t,\gamma_t,\delta_t$ and $\varepsilon_t$ 
for $t\in I_m$ by
\begin{align*}
\alpha_t
=&\, 
\frac{\sqrt{6}}{4}e^{i2\chi}\sin^2\psi\, (-1)^{m-t}
\sqrt{4(m-t-1)(m+t+2)(m+t+1)(m+t)}
\, ,
\\
\beta_t
=&\,
\sqrt{6}e^{i\chi}\sin\psi\cos\psi \,(m-2t) \sqrt{(m+t+1)(m+t)}
\, ,
\\
\gamma_t
=&\,
2(-1+\frac{6}{4}\sin^2\psi)(-1)^{m-t}t\sqrt{3!(m+t)(m-t)}
\, ,
\\
\delta_t
=&\,
-\sqrt{6}e^{-i\chi}\sin\psi\cos\psi \,(m+2t) \sqrt{(m-t+1)(m-t)}
\, ,
\\
\varepsilon_t
=&\,
-\frac{\sqrt{6}}{4}e^{-i2\chi}\sin^2\psi\, (-1)^{m-t}
\sqrt{4(m+t-1)(m-t+2)(m-t+1)(m-t)}
\, ,
\end{align*}
with $\alpha_m=0=\varepsilon_{-m}$. 
Now we put $m=\frac{n}{2}$. 
Take a self-conjugate $\pi_{m}$-eigenvector 
$\eta=\sum_{t\in I_m}d_t \bw_{t}^{m}$ of $A$. 
The self-conjugacy means $d_{-t}=\ovl{d_t}$. 
Since the $\pi_{m-1}$-part of $A$ is zero, we obtain 
$\Psi_{m-1}(\zeta_\psi,\eta)=0$ 
and this is equivalent to the following recurrence formula: 
\begin{equation}\label{dn01}
\alpha_t d_{t+2}+\beta_t d_{t+1}+\gamma_t d_{t} 
+\delta_t d_{t-1}+\varepsilon_t d_{t-2}=0
\end{equation}
for $t\in I_m$, 
where we define $d_t=0$ for $|t|\geq m+1$ as usual. 
Making use of $d_{-t}=\ovl{d_t}$, 
$\alpha_{-t}=\ovl{\varepsilon_t}$, $\beta_{-t}=-\ovl{\delta_t}$ and 
$\gamma_{-t}=\gamma_t$, where we use $(-1)^{2t}=-1$, 
we also have 
\[
\alpha_t d_{t+2}-\beta_t d_{t+1}+\gamma_t d_{t} 
-\delta_t d_{t-1}+\varepsilon_t d_{t-2}=0\, .
\]
From (\ref{dn01}) and this we get 
\begin{align}
&\,
\alpha_t d_{t+2}+\gamma_t d_{t}+\varepsilon_t d_{t-2}=0
\, ,\label{dn02}\\
&\,
\beta_t d_{t+1}+\delta_t d_{t-1}=0\, 
\label{dn03}.
\end{align}
We analyze them as follows. 

(1) $0<\psi<\frac{\pi}{2}$ case. 
We want to derive a contradiction to non-triviality of $\eta$. 
If we give a number $d_{-m}$, then by (\ref{dn02}) or 
(\ref{dn03}) 
we can inductively determine $d_{-m+2},\cdots, d_{m-1}$. 
By the self-conjugacy, we obtain the whole numbers $d_t$. 
Hence $\eta$ is uniquely determined by $d_{-m}$ or $d_m$. 
Because of non-triviality of solutions $\{d_t\}_{t\in I_m}$, 
the determinant of the following matrix must be zero for 
all $t\in I_m$: 
\[
\begin{pmatrix}
\alpha_t & \gamma_t & \varepsilon_t\\
\beta_{t+1} & \delta_{t+1} & 0 \\
0 & \beta_{t-1} & \delta_{t-1}
\end{pmatrix}
.
\]
Hence we have the equation for $t,m$ and $\psi$
\begin{equation}\label{dn04}
\alpha_t\delta_{t-1}\delta_{t+1}-\beta_{t+1}\gamma_t\delta_{t-1}
+\beta_{t-1}\beta_{t+1}\varepsilon_{t}=0
\end{equation}
for $t\in I_m$. 
If we put $t=m-1$, then we easily get 
$\sin^2\psi=\frac{3m-4}{5m-8}$. 

(1a) $n=3$ case. 
Then $m$ is equal to $\frac{3}{2}$ and 
we have $\sin^2\psi=-1<0$. 
This is a contradiction. 

(1b) $n\geq5$ case. 
If we put $t=m-2(\geq -m+1)$, then we easily get 
$\sin^2\psi=\frac{12(m-2)^2}{5m^2-24m+24}$. 
This is contradiction to $\sin^2\psi=\frac{3m-4}{5m-8}$. 

Therefore $\psi$ must be equal to $0$ or $\frac{\pi}{2}$. 

(2) $\psi=0$ case.
We have $\eta=\bw_0^m$ and 
$\gamma_t=-2(-1)^{m-t}t\sqrt{3!(m+t)(m-t)}$ 
and $\alpha_t=\beta_t=\delta_t=\varepsilon_t=0$. 
Then we get $d_{t}=0$ for $|t|\leq m-1$. 
This case $A$ is $C(D_n\setminus\summ)$. 

(3) $\psi=\frac{\pi}{2}$ case. 
We have 
$\zeta^{\phi,\chi,\frac{\pi}{2}}=\frac{\sqrt{6}}{4}e^{i2\chi}\bw_{-2}^2
-\frac{1}{2}\bw_0^2+\frac{\sqrt{6}}{4}e^{-i2\chi}\bw_{2}^2$ 
and $\beta_t=\delta_t=0$ and others are 
\begin{align*}
\alpha_t
=&\, 
\frac{\sqrt{6}}{4}e^{i2\chi}
(-1)^{m-t}
\sqrt{4(m-t-1)(m+t+2)(m+t+1)(m+t)}
\, ,
\\
\gamma_t
=&\,
(-1)^{m-t}t\sqrt{3!(m+t)(m-t)}
\, ,
\\
\varepsilon_t
=&\,
-\frac{\sqrt{6}}{4}e^{-i2\chi}(-1)^{m-t}
\sqrt{4(m+t-1)(m-t+2)(m-t+1)(m-t)}
\, .
\end{align*}
Recall vectors $\eta^{m,\chi}$ and $\wdh{\eta}^{m,\chi}$ 
which are introduced in the previous case $\T_n$. 
We easily confirm that they are independent solutions of 
$\alpha_t d_{t+2}+\gamma_t d_t+\varepsilon_t d_{t-2}=0$. 
Let $\eta=\lambda_0\eta^{m,\chi}+\lambda_1\wdh{\eta}^{m,\chi}$ be 
a self-conjugate $\pi_m$-eigenvector of $A$ where 
$\lambda_0$ and $\lambda_1$ are complex numbers. 
Since the conjugation operation $T$ satisfies 
$T\eta^{m,\chi}=e^{-i2m\chi}\eta^{m,\chi}$ 
and 
$T\wdh{\eta}^{m,\chi}=-e^{-i2m\chi}\wdh{\eta}^{m,\chi}$, 
there exist real numbers $\mu_0$ and $\mu_1$ with 
$\lambda_0=\mu_0 e^{-im\chi}$ 
and $\lambda_1=i\mu_1e^{-im\chi}$. 
The $\pi_1$-eigenvector $\Psi_1(\eta,\eta)$ must be zero 
because of the absence of $\pi_1$-part in $A$. 
By Lemma \ref{pi-product eta etahat} we have 
$\Psi_1(\eta,\eta)=2\mu_0\mu_1\sqrt{2m}2^{2m-1}\xi^{\chi,\frac{\pi}{2}}$. 
It shows $\mu_0$ or $\mu_1$ is equal to zero. 
Hence the $\pi_m$-eigenvector space of $B$ consists of scalar multiples 
of $\eta^{m,\chi}$ or $\wdh{\eta}^{m,\chi}$. 
We can easily check 
$\beta_{i}^L(\eta_r^{m,\chi})=i^{-2m}\wdh{\eta}_r^{m,\chi}$ for 
all $r\in I_m$ and 
hence may assume that the $\pi_m$-eigenvector space of $A$ is 
spanned by $\eta^{m,\chi}$ if necessary by applying $\beta_i^L$. 
Moreover we may also assume $\chi=0$ by applying 
$\beta_z^L$ where $z=e^{i\frac{\chi}{2}}$. 
Let $B$ be a right coideal generated by $\eta_r^{m}=\eta_r^{m,0}$ 
with all $r\in I_m$. 
Notice that $B$ is contained in $C^{n,0}$ which is 
defined in the discussion on the type $\T_n$. 
Hence $B$ is of type $\T_n$ or $D_n$. 
We want to show $B$ is actually of type $D_n$. 
This is proved by Proposition \ref{isomorphic} or 
the following direct calculation 
which takes us the similar characterization on the right coideal 
as Proposition \ref{group like tn}. 
In order to study it we have to clarify the type of 
the right coideal $C$ which is the integer part of $B$. 
Whether $C$ is of type $\T_n$ or $D_n$, its $\pi_2$-multiplicity 
is one. 
Hence the $\pi_2$-eigenvector of $C$ is a scalar multiple of 
$\zeta^{\phi,0,\frac{\pi}{2}}$. 
We show that there exists an angle $0\leq\phi<\frac{2\pi}{n}$ 
such that $C=C(D_n^{\phi,0,\frac{\pi}{2}}\setminus SO_{-1}(3))$. 
We already know all $\eta_r^{m}$ for $r\in I_m$ are fixed by the 
rotation of $\T_n^{0,\frac{\pi}{2}}$. 
It suffices to show $\eta_r^m \eta_s^m$ is fixed by 
the rotation of angle $\pi$ around the axis 
$\pi_1(r^{23}(-\phi))\cdot \sigma_3$ for any $r,s\in I_m$. 
This rotation is obtained by the matrix $k^{\phi}\in SU(2)$
$
k^\phi
=
\big{(}
\begin{smallmatrix}
-i\cos\phi & -\sin\phi\\
\sin\phi & i\cos\phi
\end{smallmatrix}
\big{)}
$. 
Define the $*$-homomorphism 
$\rho_{k^\phi}=(\ups_{k^\phi}\otimes\id)\circ\delta$. 

\begin{lem}
For any $m\in \frac{1}{2}+\Z_{\geq0}$ and $s\in I_m$ 
the following equalities holds.
\begin{align*}
\rho_{k^{\phi}}(\eta_{s}^m)
=
&\,
i^{2m}\sum_{r\in I_m}
(\cos2m\phi\sigma_1
\!-\!(-1)^{m-r}i\sin2m\phi\sigma_2)\!\otimes\!
(\hspace{-0.1mm}-1)^{m-r}
{2m\brack m\!-\!r}^{\!\frac{1}{2}}
w(\pi_m)_{r,s},
\\
\rho_{k^{\phi}}(\wdh{\eta}_{s}^m)
=
&\,
i^{2m}\sum_{r\in I_m}
(\cos2m\phi\sigma_1+(-1)^{m-r}i\sin2m\phi\sigma_2)\otimes
{2m\brack m-r}^{\frac{1}{2}}
w(\pi_m)_{r,s}.
\end{align*}
\end{lem}
\begin{proof}
$k^{\phi}$ is equal to $r^{12}(\pi)r^{23}(2\phi)=r^{12}(\pi)g^{2\phi,0}$ 
where $g^{\theta,\chi}$ is a matrix defined in the study of $\T_n$ case. 
Since we have 
$\rho_{k^\phi}(x)=-i\rho_{g^{-2\phi,0}}(x)$ and 
$\rho_{k^\phi}(v)=i\rho_{g^{-2\phi,0}}(v)$, 
the equalities 
$\rho_{k^\phi}(\eta_{-m}^m)
=
i^{2m}\rho_{g^{2\phi,0}}(\wdh{\eta}_{-m}^m)
$ 
and 
$\rho_{k^\phi}(\wdh{\eta}_{-m}^m)
=
i^{2m}\rho_{g^{2\phi,0}}(\eta_{-m}^m)
$ 
hold. 
By Lemma \ref{rho and eta} we obtain the desired equalities. 
\end{proof}
Now we solve the equation for $\phi$; 
$\rho_{k^\phi}(\eta_s^m\eta_t^m)=1\otimes\eta_s^m\eta_t^m$ 
for all $s,t\in I_m$. 
First we consider $s=t=-m$. 
We focus on the term of $x^{4m}$ in the both sides. 
In the left hand side, we have 
$i^{4m}(\cos 2m\phi\sigma_1-i\sin2m\phi\sigma_2)^2\otimes x^{4m}$. 
In the right one, we have 
$1\otimes x^4$. 
Hence we have 
$i^{4m}(\cos 2m\phi\sigma_1-i\sin2m\phi\sigma_2)^2=1$. 
It yields $\cos4m\phi=-1$ and hence $\phi$ must be equal to $\frac{\pi}{2n}$ 
or $\frac{3\pi}{2n}$. 
Then we obtain 
$\rho_{k^{\phi}}(\eta_{s}^m)=\pm i^{2m+1}\sigma_2\otimes \eta_{s}^m$ 
for all $s\in I_m$ with respect to 
$\phi=\frac{\pi}{2n}$ 
or $\frac{3\pi}{2n}$. 
It shows the desired fixed element property of $\eta_s^m\eta_t^m$ 
by $\rho_{k^\phi}$ for all $s,t\in I_m$ 
because of $i^{2m+1}\in\{-1,1\}$. 
Therefore we have proved $B$ is a right coideal of type $D_n$. 
Note that $B$ does not depend on the choice of 
$\phi=\frac{\pi}{2n}$ 
or $\frac{3\pi}{2n}$. 
In summary, a right coideal of type $D_n$ is made from 
the quotient $D_n\setminus \summ$ or generated by a 
$\pi_m$-eigenvector $\eta^m$. 

As in the previous case of $\T_n$ (odd $n\geq3$) 
we can also observe the similar results about a group-like unitary.
Here we treat only the case $\chi=0,\phi=\frac{\pi}{2n}$, that is, 
$Z$ is $\pi_1^{-1}(D_n^{\frac{\pi}{2n},0,\frac{\pi}{2}})$ 
in order to avoid similar arguments. 
Denote the subgroup $G_Z$ by 
$(D_n^{\frac{\pi}{2n},0,\frac{\pi}{2}})_\#$. 
Now $Z$ consists of 
$\{g^{\frac{2\pi}{n}k}\}_{0\leq k\leq 2n-1} 
\cup 
\{k^{\phi}g^{\frac{2\pi}{n}\ell}\}_{0\leq \ell\leq 2n-1}$. 
Hence the subgroup $(D_n^{\frac{\pi}{2n},0,\frac{\pi}{2}})_\#$ 
contains $(\T_n^{0,\frac{\pi}{2}})_\#$. 
We call the subsets $\{g^{\frac{2\pi}{n}k}\}_{0\leq k\leq 2n-1}$ 
and 
$\{k^{\phi}g^{\frac{2\pi}{n}\ell}\}_{0\leq \ell\leq 2n-1}$ 
the cyclic component and the reflective component, respectively. 
By definition of the restriction map, we have 
$
C((D_n^{\frac{\pi}{2n},0,\frac{\pi}{2}})_\#)
\subset 
\oplus_{0\leq k\leq 2n-1} \mathbb{M}_{s(k)}(\C)
\oplus
\oplus_{0\leq \ell\leq 2n-1} \mathbb{M}_{t(\ell)}(\C)$ 
where $s(k)=2$ except for $k=0,2$ and 
$t(\ell)=2$ except for $\ell=\frac{n-1}{2}, \frac{3n-1}{2}$. 
First note the equality $k^\phi g^\theta=k^{\phi+\frac{\theta}{2}}$. 
Then we get 
\[\rho_{k^\phi g^\theta}(\eta_s^m)
=
i^{2m}\sum_{r\in I_m}(-\sin m\theta\sigma_1-(-1)^{m-r}i\cos m\theta\sigma_2)
\otimes 
(-1)^{m-r}{2m \brack m\!-\!r}^{\!\frac{1}{2}}w(\pi_m)_{r,s}
.
\] 
Hence it immediately derives 
$\rho_{k^\phi g^{\frac{2\pi}{n}k}}(\eta_s^m)
=i^{2m+1}(-1)^{k+1}\sigma_2\otimes \eta_s^m$ for $0\leq k \leq n-1$. 
Define a self-adjoint unitary 
$\wdt{w}=w  
\oplus \oplus_{\ell\neq \frac{n-1}{2},\frac{3n-1}{2}} 
(-1)^{\ell+1} i^{2m+1}\sigma_2
\oplus
1
\oplus
-1
$ 
where $w$ is a self-adjoint group-like unitary defined 
in the previous $\T_n$ case. 

\begin{prop}
A self-adjoint unitary $\wdt{w}$ is a group-like element in 
$C((D_n^{\frac{\pi}{2n},0,\frac{\pi}{2}})_\#)$. 
\end{prop}

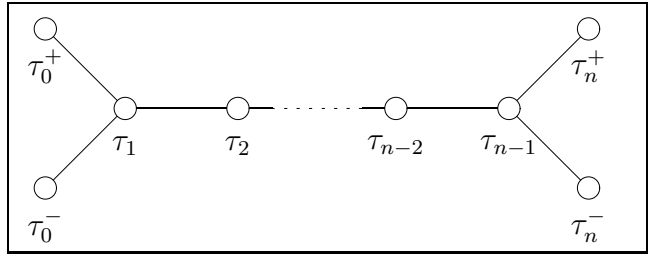
\begin{figure}
\begin{center}
\setlength{\unitlength}{0.5mm}

\fbox{
\begin{picture}(160,62)
\put(5,15){\circle{6}}
\put(7.12,17.12){\line(1,1){17.0}}
\put(26.22,36.22){\circle{6}}
\put(24.1,38.34){\line(-1,1){17.0}}
\put(5,57.4){\circle{6}}
\put(29.22,36.22){\line(1,0){24}}
\put(56.22,36.22){\circle{6}}
\put(59.22,36.22){\line(1,0){6}}

\dashline{1}(65.22,36.22)(89.22,36.22)

\put(89.22,36.22){\line(1,0){6}}

\put(98.22,36.22){\circle{6}}
\put(101.22,36.22){\line(1,0){24}}

\put(128.22,36.22){\circle{6}}

\put(149.44,15){\circle{6}}
\put(149.44,57.4){\circle{6}}

\put(147.32,17.12){\line(-1,1){17.0}}
\put(130.34,38.34){\line(1,1){17.0}}

\put(5,5){\makebox(0,0){$\tau_0^-$}}
\put(26.22,26.22){\makebox(0,0){$\tau_1$}}
\put(5,47.4){\makebox(0,0){$\tau_0^+$}}
\put(56.22,26.22){\makebox(0,0){$\tau_2$}}
\put(98.22,26.22){\makebox(0,0){$\tau_{n-2}$}}
\put(128.22,26.22){\makebox(0,0){$\tau_{n-1}$}}
\put(149.44,5){\makebox(0,0){$\tau_n^-$}}
\put(149.44,47.4){\makebox(0,0){$\tau_n^+$}}
\end{picture}
}
\caption{The McKay diagram for 
$(D_n^{\frac{\pi}{2n},0,\frac{\pi}{2}})_\#$, $n+3$ nodes}
\label{D_ns}
\end{center}
\end{figure}

Let $C(\Z_2)$ be a Hopf $*$-subalgebra $\C+\C \wdt{w}$ in 
$C((D_n^{\frac{\pi}{2n},0,\frac{\pi}{2}})_\#)$. 
Then we get $B=\{a\in \csum \mid 
(\pi_Z\otimes \id)\circ \delta(a)\in C(\Z_2)\otimes \C a\}$. 
We study all the irreducible representations of 
$(D_n^{\frac{\pi}{2n},0,\frac{\pi}{2}})_\#$. 
The McKay diagram of 
$(D_n^{\frac{\pi}{2n},0,\frac{\pi}{2}})_\#
\subset \summ
$ about $\pi_{\frac{1}{2}}$ 
is shown in Figure \ref{D_ns}. 
In Figure \ref{D_ns}, $\tau_0^+$ is the trivial representation 
and $\tau_1$ is the restriction of $w(\pi_{\frac{1}{2}})$ to 
$D_n^{\frac{\pi}{2n},0,\frac{\pi}{2}}$. 
We study where the above group-like unitary $\wdt{w}$ sits. 
The restriction of $w(\pi_{\frac{1}{2}})$ to 
$(D_n^{\frac{\pi}{2n},0,\frac{\pi}{2}})_\#$ is denoted simply by 
$w(\pi_{\frac{1}{2}})|$. 
For the $\pi_{\frac{1}{2}}$-module $W_1$ we use an orthonormal basis 
$\xi_\pm=\C
\big{(}
\begin{smallmatrix}
1\\
\pm i
\end{smallmatrix}
\big{)}
$. 
Then we obtain 
$w(\pi_{\frac{1}{2}})|\cdot\xi_{\pm}
=
\xi_{\pm}\otimes w_{\pm} 
$ 
for the cyclic component 
and 
$w(\pi_{\frac{1}{2}})|\cdot\xi_{\pm}
=
\xi_{\mp}\otimes v_{\pm} 
$ 
for the the reflective component, 
where $w_{\pm}$ have been already defined in the previous 
$\T_n$ case and 
\[
v_{\pm}= 
\oplus_{\ell\neq \frac{n-1}{2},\frac{3n-1}{2}}
\Big{(}-i\cos\Big{(}\phi+\frac{\pi}{n}\ell\Big{)}\,\sigma_1
\pm i\sin\Big{(}\phi+\frac{\pi}{n}\ell\Big{)}\,\sigma_2\
\Big{)}
\oplus
\pm i
\oplus 
\mp i
.\] 
Hence $v_{\pm}^*=-v_{\pm}$ and $v_{\pm}^2=-1$. 
Let us consider the tensor product 
$W_1\otimes W_1$. 
Its irreducible submodules are 
$W_0^{\pm}=\C(\xi_{+}\otimes\xi_{+}\mp\xi_{-}\otimes\xi_{-})$ and 
$W_2=\C\xi_{+}\otimes\xi_{-}+\C\xi_{-}\otimes\xi_{+}$ which give 
$\tau_0^{\pm}$ and $\tau_2$, respectively. 
Hence for $\tau_0^{-}$, the cyclic component acts trivially and 
the reflective component acts as the multiplication of $-1$. 
Proceeding tensor products and decompositions, 
we get the $\tau_j$-module 
$W_j=\C\xi_{+}\otimes\xi_{-}\otimes\cdots\xi_{-}\otimes\xi_{+}
+\C\xi_{-}\otimes\xi_{+}\otimes\cdots\xi_{+}\otimes\xi_{-}$ 
for $1\leq j \leq n-1$, 
where the length of the words $\xi_{\pm}$ is $j$. 
The cyclic component acts on $W_j$ 
as the direct sum module of 
$w_{+}w_{-}\cdots w_{-}w_{+}$ and $w_{-}w_{+}\cdots w_{+}w_{-}$. 
The reflective component 
acts there 
$v_{+}v_{-}\cdots v_{-}v_{+}$ and $v_{-}v_{+}\cdots v_{+}v_{-}$. 
When $j$ is odd, we have 
\[
v_{\pm}v_{\mp}\cdots v_{\mp}v_{\pm}
=
\oplus_{\ell\neq \frac{n-1}{2},\frac{3n-1}{2}}
i^j\Big{(}
-\cos j\Big{(}\phi+\frac{\pi}{n}\ell\Big{)}\,\sigma_1
\pm\sin j\Big{(}\phi+\frac{\pi}{n}\ell\Big{)}\,\sigma_2
\Big{)}
\oplus
\pm i
\oplus
\mp i
.
\] 
Finally considering the tensor products module $W_1\otimes W_{n-1}$, 
we get its one-dimensional submodules 
$W_n^{\pm}=\C
(\xi_{\pm}\otimes\xi_{\mp}\otimes\cdots\xi_{\mp}\otimes\xi_{\pm}
\pm i
\xi_{\mp}\otimes\xi_{\pm}\otimes\cdots\xi_{\pm}\otimes\xi_{\mp})
$. 
We investigate the action of 
$(D_n^{\frac{\pi}{2n},0,\frac{\pi}{2}})_\#$ on them. 
The cyclic component acts by $w_{+}w_{-}\cdots w_{-}w_{+}=w$. 
Since we have the equality 
$v_{\pm}v_{\mp}\cdots v_{\mp}v_{\pm}
=
\oplus_{k\neq \frac{n-1}{2},\frac{3n-1}{2}}
i^n(-1)^k (-1)^{\pm1+1}\sigma_2
\oplus
\pm i
\oplus
\mp i
$, 
the reflection component acts on $W_n^{\pm}$ by 
$
\oplus_{k\neq \frac{n-1}{2},\frac{3n-1}{2}}
i^{n+1}(-1)^k (-1)^{\frac{\pm1+1}{2}}\sigma_2
\oplus
\pm 1
\oplus
\mp 1
$. 
Hence the self-adjoint group-like unitary $\wdt{w}$ is 
the unitary representation on $W_n^{+}$. 
Especially we get $\ovl{\tau_n^{\pm}}=\tau_n^{\pm}$. 
Recall the classical $q=1$ case $D_n^*\subset SU(2)$. 
Then it has the same McKay diagram, however, 
we know $\ovl{\tau_n^{\pm}}=\tau_n^{\mp}$. 
This shows the difference of the representation theory of 
$(D_n^{\frac{\pi}{2n},0,\frac{\pi}{2}})_\#\subset\summ$ and 
$D_n^*\subset SU(2)$. 

We have classified all the right coideals of type $D_n$ (odd $n\geq3$). 
They are $C(D_n\setminus \summ)$ or $C^*(\eta^m)$ 
up to conjugacy by the maximal torus action $\beta^L$, 
where $C^*(\eta^m)$ is a right coideal generated by $\eta_r^m$ 
for all $r\in I_m$. 
Finally we end this section with the following result. 

\begin{prop}\label{isomorphic}
As $\summ$-covariant systems, 
$C(D_n\setminus \summ)$ and $C^*(\eta^m)$ are isomorphic.
\end{prop}

\begin{proof}
Let $g$ be a matrix in $SU(2)$ 
$
r^{13}(\frac{\pi}{2})
=
\begin{pmatrix}
\sqrt{2}^{-1} & \sqrt{2}^{-1}\\
-\sqrt{2}^{-1} & \sqrt{2}^{-1} 
\end{pmatrix}
$. 
Recall the $\summ$-homomorphism 
$\rho_g=(\upsilon_g\otimes\id)\circ\delta:
\csum\longrightarrow \mathbb{B}(\C^2)\otimes \csum$. 
$\upsilon_g$ satisfies
$
\begin{pmatrix}
\upsilon_g(x) & \upsilon_g(u) \\
\upsilon_g(v) & \upsilon_g(y)
\end{pmatrix}
=
\begin{pmatrix}
\sqrt{2}^{-1} \sigma_1 & -\sqrt{2}^{-1}\sigma_2 \\ 
-\sqrt{2}^{-1} \sigma_2 & \sqrt{2}^{-1}\sigma_1
\end{pmatrix}
$. 
Then we obtain
\begin{align*}
\upsilon_g(w(\pi_m)_{-m,r})
=
&\,
\sqrt{2}^{-2m}(-1)^{m+r}{2m \brack m-r}^{\frac{1}{2}}
\sigma_1^{m-r}\sigma_2^{m+r}
, \\
\upsilon_g(w(\pi_m)_{m,r})
=
&\,
\sqrt{2}^{-2m}(-1)^{m-r}{2m \brack m-r}^{\frac{1}{2}}
\sigma_2^{m-r}\sigma_1^{m+r}
\end{align*}
for all $r\in I_m$. 
We also have the equality
\[
(-1)^{m+r}\sigma_1^{m-r}\sigma_2^{m+r}
+(-1)^{m-r}\sigma_2^{m-r}\sigma_1^{m+r}
=
(-1)^{2m}(\sigma_2-\sigma_1)\,
\]
for all $r\in I_m$. 
Then we have 
\begin{align*}
\rho_g(x^n+v^n)
=
&\,
\rho_g(w(\pi_m)_{-m,-m}+w(\pi_m)_{m,-m})
\\
=
&\,
\sum_{r\in I_m}
(\upsilon_g(w(\pi_m)_{-m,r})+\upsilon_g(w(\pi_m)_{m,r}))
\otimes
w(\pi_m)_{r,-m}
\\
=
&\,
\sqrt{2}^{-2m}
\sum_{r\in I_m}
{2m \brack m-r}^{\frac{1}{2}}
((-1)^{m+r}\sigma_1^{m-r}\sigma_2^{m+r}
+(-1)^{m-r}\sigma_2^{m-r}\sigma_1^{m+r})\\
&\qquad\qquad\qquad\otimes
w(\pi_m)_{r,-m}
\\
=
&\,
-
\sqrt{2}^{-2m}
\sum_{r\in I_m}
{2m \brack m-r}^{\frac{1}{2}}
(\sigma_2-\sigma_1)
\otimes
w(\pi_m)_{r,-m}\\
=
&\,
-
\sqrt{2}^{-2m}
(\sigma_2-\sigma_1)
\otimes
\eta_{-m}^m
\,.
\end{align*}
Let us define the self-adjoint unitary 
$\nu=\sqrt{2}^{-1}(\sigma_2-\sigma_1)$. 
Then we have 
$
\rho_g(w(\pi_m)_{-m,s}+w(\pi_m)_{m,s})
=
-
\sqrt{2}^{1-2m}
\nu\otimes \eta_s^m
$ 
for all $s\in I_m$. 
Hence we get 
$\rho_g(C(D_n\setminus\summ))\subset C^*(\nu)\otimes C^*(\eta^m)$. 
Take a $*$-homomorphism $\omega:C^*(\nu)\longrightarrow \C$ and 
we have an $\summ$-homomorphism 
$(\omega\otimes\id)\circ\rho_g:
C(D_n\setminus\summ)\longrightarrow C^*(\eta^m)$. 
This is clearly surjective and the injectivity follows from 
Remark \ref{faithfulness}. 

\end{proof}

\section{Appendix}

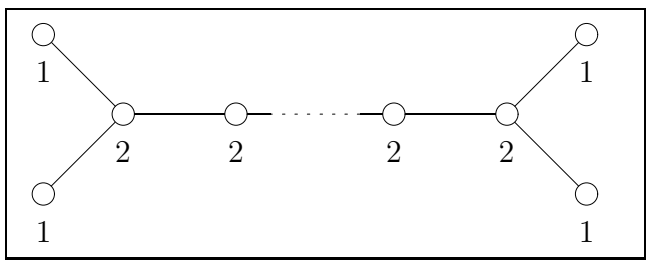
\begin{figure}
\begin{center}
\setlength{\unitlength}{0.5mm}

\fbox{
\begin{picture}(160,62)
\put(5,15){\circle{6}}
\put(7.12,17.12){\line(1,1){17.0}}
\put(26.22,36.22){\circle{6}}
\put(24.1,38.34){\line(-1,1){17.0}}
\put(5,57.4){\circle{6}}
\put(29.22,36.22){\line(1,0){24}}
\put(56.22,36.22){\circle{6}}
\put(59.22,36.22){\line(1,0){6}}

\dashline{1}(65.22,36.22)(89.22,36.22)

\put(89.22,36.22){\line(1,0){6}}

\put(98.22,36.22){\circle{6}}
\put(101.22,36.22){\line(1,0){24}}

\put(128.22,36.22){\circle{6}}

\put(149.44,15){\circle{6}}
\put(149.44,57.4){\circle{6}}

\put(147.32,17.12){\line(-1,1){17.0}}
\put(130.34,38.34){\line(1,1){17.0}}

\put(5,5){\makebox(0,0){1}}
\put(26.22,26.22){\makebox(0,0){2}}
\put(5,47.4){\makebox(0,0){1}}
\put(56.22,26.22){\makebox(0,0){2}}
\put(98.22,26.22){\makebox(0,0){2}}
\put(128.22,26.22){\makebox(0,0){2}}
\put(149.44,5){\makebox(0,0){1}}
\put(149.44,47.4){\makebox(0,0){1}}

\end{picture}
}
\caption{$D_n^*$ ($n\geq 2$), $n+3$ nodes}
\label{D_n*}
\end{center}
\end{figure}

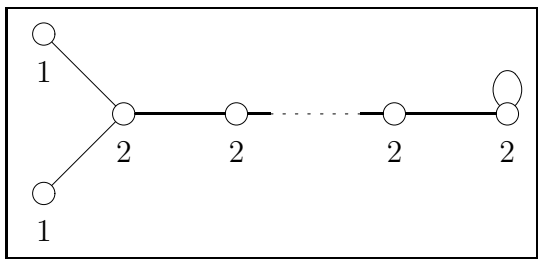
\begin{figure}
\begin{center}
\setlength{\unitlength}{0.5mm}

\fbox{
\begin{picture}(131,62)
\put(5,15){\circle{6}}
\put(7.12,17.12){\line(1,1){17.0}}
\put(26.22,36.22){\circle{6}}
\put(24.1,38.34){\line(-1,1){17.0}}
\put(5,57.4){\circle{6}}
\put(29.22,36.22){\line(1,0){24}}
\put(56.22,36.22){\circle{6}}
\put(59.22,36.22){\line(1,0){6}}

\dashline{1}(65.22,36.22)(89.22,36.22)

\put(89.22,36.22){\line(1,0){6}}

\put(98.22,36.22){\circle{6}}
\put(101.22,36.22){\line(1,0){24}}

\put(128.22,36.22){\circle{6}}
\put(130.34,42.58){\arc{12}{2.36}{-2.36}}
\put(128.22,44.71){\arc{6}{3.93}{-0.785}}
\put(126.10,42.58){\arc{12}{5.50}{0.79}}

\put(5,5){\makebox(0,0){1}}
\put(26.22,26.22){\makebox(0,0){2}}
\put(5,47.4){\makebox(0,0){1}}
\put(56.22,26.22){\makebox(0,0){2}}
\put(98.22,26.22){\makebox(0,0){2}}
\put(128.22,26.22){\makebox(0,0){2}}

\end{picture}
}
\caption{$D_n$ (odd $n\geq 3$), $(n+3)\slash 2$ nodes}
\label{D_n}
\end{center}
\end{figure}

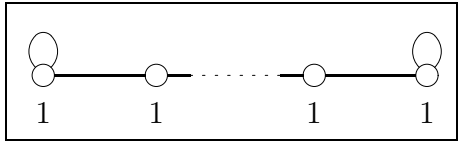
\begin{figure}
\begin{center}
\setlength{\unitlength}{0.5mm}

\fbox{
\begin{picture}(110,32)
\put(5,15){\circle{6}}
\put(7.12,21.36){\arc{12}{2.36}{-2.36}}
\put(5,23.49){\arc{6}{3.93}{-0.785}}
\put(2.88,21.36){\arc{12}{5.50}{0.79}}
\put(8,15){\line(1,0){24}}
\put(35,15){\circle{6}}
\put(38,15){\line(1,0){6}}
\dashline{1}(44,15)(68,15)

\put(68,15){\line(1,0){6}}
\put(77,15){\circle{6}}
\put(80,15){\line(1,0){24}}
\put(107,15){\circle{6}}
\put(109.12,21.36){\arc{12}{2.36}{-2.36}}
\put(107,23.49){\arc{6}{3.93}{-0.785}}
\put(104.88,21.36){\arc{12}{5.50}{0.79}}

\put(5,5){\makebox(0,0){1}}
\put(35,5){\makebox(0,0){1}}
\put(77,5){\makebox(0,0){1}}
\put(107,5){\makebox(0,0){1}}

\end{picture}
}
\caption{$A_m'$, ($m\geq3$) nodes}
\label{A_m'}
\end{center}
\end{figure}

We list all the connected graphs of norm $2$ from 
Figure \ref{1} to Figure \ref{A_infty'} 
and the spectral patterns of ergodic systems 
which sit at the vertex of $1$ (the entry of 
the Perron-Frobenius vector). 
Because we need not to specify the detailed pattern of type $A_m'$ 
in our classification program, 
we do not list it in the $A_m'$ case.


\begin{figure}[h]
\begin{center}
\setlength{\unitlength}{0.5mm}
\fbox{
\begin{picture}(26,23)

\put(12.73,15){\circle{6}}
\put(6.36,17.12){\arc{12}{0.79}{-3.93}}
\put(6.36,13.12){\arc{12}{3.93}{-0.79}}
\put(4.24,15){\arc{6}{2.36}{-2.36}}

\put(21.21,15){\arc{6}{5.5}{0.79}}
\put(19.09,13.12){\arc{12}{3.93}{-0.79}}
\put(19.09,17.12){\arc{12}{0.79}{-3.93}}

\put(12.73,5){\makebox(0,0){1}}

\end{picture}
}
\caption{$1$}
\label{1}
\end{center}
\end{figure}

\begin{figure}
\begin{center}
\fbox{
\begin{picture}(130,130)

\put(115,65){\circle{6}}
\put(100.36,100.36){\circle{6}}
\put(65,115){\circle{6}}
\put(29.64,100.36){\circle{6}}
\put(15,65){\circle{6}}
\put(29.64,29.64){\circle{6}}
\put(100.36,29.64){\circle{6}}

\put(125,65){\makebox(0,0){1}}
\put(107.43,107.43){\makebox(0,0){1}}
\put(65,125){\makebox(0,0){1}}
\put(22.57,107.43){\makebox(0,0){1}}
\put(5,65){\makebox(0,0){1}}
\put(22.57,22.57){\makebox(0,0){1}}
\put(107.43,22.57){\makebox(0,0){1}}

\put(65,65){\arc{100}{5.56}{-0.061}}

\put(65,65){\arc{100}{4.774}{-0.846}}

\put(65,65){\arc{100}{3.989}{-1.632}}

\put(65,65){\arc{100}{3.204}{-2.417}}

\put(65,65){\arc{100}{2.419}{-3.202}}


\put(65,65){\arc{100}{2.027}{-3.988}}


\put(65,65){\arc{100}{0.8486}{-5.166}}

\put(65,65){\arc{100}{0.0631}{-5.558}} 

\put(65,65){\arc{100}{1.117}{-4.256}}  


\end{picture}
}
\caption{$\T_m$, $m$ nodes}
\label{Torus}
\end{center}
\end{figure}
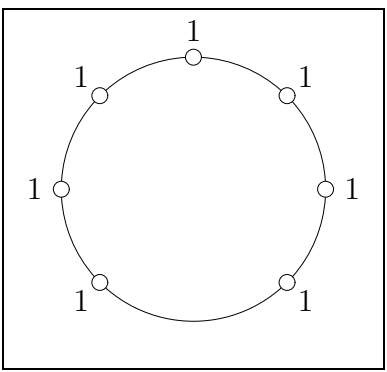


\begin{figure}
\begin{center}
\setlength{\unitlength}{0.5mm}

\fbox{
\begin{picture}(183,20)
\dashline{1}(4,15)(23,15)

\put(23,15){\line(1,0){6}}
\put(32,15){\circle{6}}
\put(35,15){\line(1,0){24}}
\put(62,15){\circle{6}}
\put(65,15){\line(1,0){24}}
\put(92,15){\circle{6}}
\put(95,15){\line(1,0){24}}
\put(122,15){\circle{6}}
\put(125,15){\line(1,0){24}}
\put(152,15){\circle{6}}
\put(155,15){\line(1,0){6}}

\dashline{1}(161,15)(179,15)

\put(32,5){\makebox(0,0){1}}
\put(62,5){\makebox(0,0){1}}
\put(92,5){\makebox(0,0){1}}
\put(122,5){\makebox(0,0){1}}
\put(152,5){\makebox(0,0){1}}

\end{picture}
}
\caption{$\T$}
\label{Torus}
\end{center}
\end{figure}


\begin{figure}
\begin{center}
\setlength{\unitlength}{0.5mm}

\fbox{
\begin{picture}(145,21)
\put(5,15){\circle{6}}
\put(8,15){\line(1,0){24}}
\put(35,15){\circle{6}}
\put(38,15){\line(1,0){24}}
\put(65,15){\circle{6}}
\put(68,15){\line(1,0){24}}
\put(95,15){\circle{6}}
\put(98,15){\line(1,0){6}}

\dashline{1}(104,15)(140,15)

\put(5,5){\makebox(0,0){1}}
\put(35,5){\makebox(0,0){2}}
\put(65,5){\makebox(0,0){3}}
\put(95,5){\makebox(0,0){4}}

\end{picture}
}
\caption{$SU(2)$ }
\label{SU(2)}
\end{center}
\end{figure}


\begin{figure}
\begin{center}
\setlength{\unitlength}{0.5mm}

\fbox{
\begin{picture}(160,62)
\put(5,15){\circle{6}}
\put(7.12,17.12){\line(1,1){17.0}}
\put(26.22,36.22){\circle{6}}
\put(24.1,38.34){\line(-1,1){17.0}}
\put(5,57.4){\circle{6}}
\put(29.22,36.22){\line(1,0){24}}
\put(56.22,36.22){\circle{6}}
\put(59.22,36.22){\line(1,0){24}}
\put(86.22,36.22){\circle{6}}
\put(89.22,36.22){\line(1,0){24}}
\put(116.22,36.22){\circle{6}}
\put(119.22,36.22){\line(1,0){6}}

\dashline{1}(125.22,36.22)(155.22,36.22)

\put(5,5){\makebox(0,0){1}}
\put(26.22,26.22){\makebox(0,0){2}}
\put(5,47.4){\makebox(0,0){1}}
\put(56.22,26.22){\makebox(0,0){2}}
\put(86.22,26.22){\makebox(0,0){2}}
\put(116.22,26.22){\makebox(0,0){2}}

\end{picture}
}
\caption{$D_\infty^*$}
\label{D_infty*}
\end{center}
\end{figure}


\begin{figure}
\begin{center}
\setlength{\unitlength}{0.5mm}

\fbox{
\begin{picture}(130,80)

\put(5,75){\circle{6}}
\put(8,75){\line(1,0){24}}
\put(35,75){\circle{6}}
\put(38,75){\line(1,0){24}}
\put(65,75){\circle{6}}
\put(68,75){\line(1,0){24}}
\put(95,75){\circle{6}}
\put(98,75){\line(1,0){24}}
\put(125,75){\circle{6}}

\put(65,72){\line(0,-1){24}}
\put(65,45){\circle{6}}
\put(65,42){\line(0,-1){24}}
\put(65,15){\circle{6}}

\put(5,65){\makebox(0,0){1}}
\put(35,65){\makebox(0,0){2}}
\put(70,65){\makebox(0,0){3}}
\put(95,65){\makebox(0,0){2}}
\put(125,65){\makebox(0,0){1}}
\put(75,45){\makebox(0,0){2}}
\put(75,15){\makebox(0,0){1}}

\end{picture}
}
\caption{$A_4^*$}
\label{A_4}
\end{center}
\end{figure}


\begin{figure}
\begin{center}
\setlength{\unitlength}{0.5mm}

\fbox{
\begin{picture}(190,50)

\put(5,45){\circle{6}}
\put(8,45){\line(1,0){24}}
\put(35,45){\circle{6}}
\put(38,45){\line(1,0){24}}
\put(65,45){\circle{6}}
\put(68,45){\line(1,0){24}}
\put(95,45){\circle{6}}
\put(98,45){\line(1,0){24}}
\put(125,45){\circle{6}}
\put(128,45){\line(1,0){24}}
\put(155,45){\circle{6}}
\put(158,45){\line(1,0){24}}
\put(185,45){\circle{6}}

\put(95,42){\line(0,-1){24}}
\put(95,15){\circle{6}}

\put(5,35){\makebox(0,0){1}}
\put(35,35){\makebox(0,0){2}}
\put(70,35){\makebox(0,0){3}}
\put(100,35){\makebox(0,0){4}}
\put(125,35){\makebox(0,0){3}}
\put(155,35){\makebox(0,0){2}}
\put(185,35){\makebox(0,0){1}}

\put(105,15){\makebox(0,0){2}}

\end{picture}
}
\caption{$S_4^*$}
\label{S_4}
\end{center}
\end{figure}


\begin{figure}
\begin{center}
\setlength{\unitlength}{0.5mm}

\fbox{
\begin{picture}(220,50)

\put(5,45){\circle{6}}
\put(8,45){\line(1,0){24}}
\put(35,45){\circle{6}}
\put(38,45){\line(1,0){24}}
\put(65,45){\circle{6}}
\put(68,45){\line(1,0){24}}
\put(95,45){\circle{6}}
\put(98,45){\line(1,0){24}}
\put(125,45){\circle{6}}
\put(128,45){\line(1,0){24}}
\put(155,45){\circle{6}}
\put(158,45){\line(1,0){24}}
\put(185,45){\circle{6}}
\put(188,45){\line(1,0){24}}
\put(215,45){\circle{6}}

\put(65,42){\line(0,-1){24}}
\put(65,15){\circle{6}}

\put(5,35){\makebox(0,0){2}}
\put(35,35){\makebox(0,0){4}}
\put(70,35){\makebox(0,0){6}}
\put(100,35){\makebox(0,0){5}}
\put(125,35){\makebox(0,0){4}}
\put(155,35){\makebox(0,0){3}}
\put(185,35){\makebox(0,0){2}}
\put(215,35){\makebox(0,0){1}}

\put(75,15){\makebox(0,0){3}}

\end{picture}
}
\caption{$A_5^*$}
\label{A_5}
\end{center}
\end{figure}

\begin{figure}
\begin{center}
\setlength{\unitlength}{0.5mm}

\fbox{
\begin{picture}(50,30)
\put(10,15){\circle{6}}
\put(12.12,21.36){\arc{12}{2.36}{-2.36}}
\put(10,23.49){\arc{6}{3.93}{-0.785}}
\put(7.88,21.36){\arc{12}{5.50}{0.79}}
\put(13,15){\line(1,0){24}}
\put(40,15){\circle{6}}
\put(42.12,21.36){\arc{12}{2.36}{-2.36}}
\put(40,23.49){\arc{6}{3.93}{-0.785}}
\put(37.8,21.36){\arc{12}{5.50}{0.79}}

\put(10,5){\makebox(0,0){1}}
\put(40,5){\makebox(0,0){1}}

\end{picture}
}
\caption{$D_1$}
\label{D_1}
\end{center}
\end{figure}

\begin{figure}
\begin{center}
\setlength{\unitlength}{0.5mm}

\fbox{
\begin{picture}(139,32)
\put(5,15){\circle{6}}
\put(7.12,21.36){\arc{12}{2.36}{-2.36}}
\put(5,23.49){\arc{6}{3.93}{-0.785}}
\put(2.88,21.36){\arc{12}{5.50}{0.79}}
\put(8,15){\line(1,0){24}}
\put(35,15){\circle{6}}
\put(38,15){\line(1,0){24}}
\put(65,15){\circle{6}}
\put(68,15){\line(1,0){24}}
\put(95,15){\circle{6}}
\put(98,15){\line(1,0){6}}

\dashline{1}(104,15)(134,15)

\put(5,5){\makebox(0,0){1}}
\put(35,5){\makebox(0,0){1}}
\put(65,5){\makebox(0,0){1}}
\put(95,5){\makebox(0,0){1}}
\end{picture}
}
\caption{$A_\infty'$ }
\label{A_infty'}
\end{center}
\end{figure}
 
\begin{tabular}{ll}
$1$ & $\oplus_{\nu\in \frac{1}{2}}(2\nu+1)\pi_\nu$,\\
$\T_{n}$\, (even $n\geq2$) &
$\oplus_{k\in \Z_{\geq0}}(1+2\big{[}\frac{2k}{n}\big{]})\pi_k$, \\
$\T_{n}$\, (odd $n\geq3$)  &
$\oplus_{k\in\Z_{\geq0}}\big{(}1+2\big{[}\frac{k}{n}\big{]}\big{)}\pi_k
\oplus 
\oplus_{k\in \Z_{\geq0}} 
2\big{[}\frac{2k+n+1}{2n}\big{]}\pi_{k+\frac{1}{2}}$,\\
$\T$ &
$\oplus_{k\in \Z_{\geq0}}\pi_k$, \\
$SU(2)$ & $\pi_0$,\\
$D_m^*$\, $(m\geq 2)$ & 
$\oplus_{k\in\Z_{\geq0}}\big{(}\frac{1+(-1)^k}{2}
+\big{[}\frac{k}{m}\big{]}\big{)} \pi_k$,\\
$D_\infty^*$ & 
$\oplus_{k\in \Z_{\geq0}}\pi_{2k}$, \\
$A_4^*$ & 
$\pi_0\oplus \pi_3 \oplus \pi_4\oplus 2\pi_6\oplus \pi_7 \oplus \cdots$,\\
$S_4^*$ &  
$\pi_0 \oplus \pi_4 \oplus \pi_6 \oplus \pi_8 \oplus \pi_9 
\oplus \pi_{10} \oplus \cdots$, \\
$A_5^*$ & 
$\pi_0 \oplus \pi_6 \oplus \pi_{10} \oplus \pi_{12} \oplus \cdots$, \\
$D_n$\, (odd $n\geq 1$) & 
$\oplus_{k\in\Z_{\geq0}}
\big{(}\frac{1+(-1)^k}{2}+\big{[}\frac{k}{n}\big{]}\big{)}\pi_k 
\oplus 
\oplus_{k\in\Z_{\geq0}} 
\big{(}\big{[}\frac{2k+1}{n}\big{]}- \big{[}\frac{k}{n}\big{]}\big{)}
\pi_{k+\frac{1}{2}}.
$
\end{tabular}

\end{document}